\documentclass[10pt]{amsart}
      \usepackage[mathscr]{eucal}
      \usepackage{amsmath,amsfonts}
      \parskip\smallskipamount
          \newtheorem{theorem}{Theorem}[section]
      
      \newtheorem{proposition}[theorem]{Proposition}
      \newtheorem{corollary}[theorem]{Corollary}
      \newtheorem{lemma}[theorem]{Lemma}

      \newtheorem{remark}[theorem]{Remark}
      \makeatletter
      \@addtoreset{equation}{section}
      \makeatother

      \newcommand{\BB}{{\mathbb B}}
      \newcommand{\CC}{{\mathbb C}}
      \newcommand{\NN}{{\mathbb N}}
      
      \newcommand{\ZZ}{{\mathbb Z}}
      \newcommand{\DD}{{\mathbb D}}
      \newcommand{\RR}{{\mathbb R}}
      \newcommand{\FF}{{\mathbb F}}

      \newcommand{\cA}{{\mathcal A}}
      
      \newcommand{\cC}{{\mathcal C}}

      \newcommand{\cF}{{\mathcal F}}
      
      \newcommand{\cH}{{\mathcal H}}
      \newcommand{\cK}{{\mathcal K}}

      \newcommand{\cP}{{\mathcal P}}

      \newcommand{\cU}{{\mathcal U}}

      \newcommand{\cX}{{\mathcal X}}

      \newdimen\expt
      \def\boxit#1{\setbox0\hbox{$\displaystyle{#1}$}
            \hbox{\lower.4\expt
       \hbox{\lower3\expt\hbox{\lower\dp0
            \hbox{\vbox{\hrule height.4\expt
       \hbox{\vrule width.4\expt\hskip3\expt
            \vbox{\vskip3\expt\box0\vskip2\expt}%
       \hskip3\expt\vrule width.4\expt}\hrule height.4\expt}}}}}}
      
\begin{document}
       \pagestyle{myheadings}
      \markboth{ Gelu Popescu}{ Free holomorphic functions on the unit ball of $B(\cH)^n$ }

      \title [ Free holomorphic functions on the unit ball of $B(\cH)^n$   ]
      {  Free holomorphic functions on the unit ball of $B(\cH)^n$}
        \author{Gelu Popescu}
\date{December 27, 2005}
      \thanks{Research supported in part by an NSF grant}

      \address{Department of Mathematics, The University of Texas
      at San Antonio \\ San Antonio, TX 78249, USA}
      \email{\tt gelu.popescu@utsa.edu}

      %
      \maketitle

\bigskip

\section*{Contents}
{\it

\quad Introduction

\begin{enumerate}
   \item[1.]  Free holomorphic functions and Hausdorff derivations
 \item[2.]  Cauchy, Liouville, and Schwartz type results for free holomorphic functions
\item[ 3.]  Algebras of free holomorphic functions
\item[ 4.]   Free analytic functional calculus  and noncommutative Cauchy transforms
 \item[ 5.]  Weierstrass and Montel theorems for free holomorphic functions
  \item[6.]  Free pluriharmonic functions and noncommutative Poisson transforms
\item[7.] Hardy spaces of free holomorphic functions
   \end{enumerate}

\quad References

}

\bigskip

\section*{Introduction}

The Shilov-Arens-Calderon theorem (\cite{S}, \cite{AC})  states that
if $a_1,\ldots, a_n$ are elements of a commutative Banach algebra $A$ with the joint spectrum  included in  a domain $\Omega\subset \CC^n$,
then the algebra homomorphism
$$\CC[z_1,\ldots,z_n]\ni p\mapsto p(a_1,\ldots, a_n)\in A
$$
extends to a continuous homomorphism from the algebra $Hol(\Omega)$, of holomorphic functions on $\Omega$, to the algebra $A$.
This result was greatly improved by the pioneering work of J.L.~Taylor (\cite{T1}, \cite{T2}, \cite{T3}) who introduced  a ``better'' notion of joint spectrum for $n$-tuples of commuting operators, which is now called Taylor spectrum, and developed an analytic functional calculus.
Stated for  the open unit ball of $\CC^n$,
$$
\BB_n:=\{(\lambda_1,\ldots, \lambda_n)\in \CC^n: \ |\lambda|^2+\cdots +|\lambda_n|^2<1\},
$$
his result states that
 if $[T_1,\ldots, T_n]\in B(\cH)^n$
is an $n$-tuple of commuting bounded linear operators on a Hilbert space $\cH$ with Taylor spectrum $\sigma(T_1,\ldots, T_n)\subset \BB_n$, then there is a unique continuous unital algebra homomorphism
$$
Hol(\BB_n)\ni f\mapsto f(T_1,\ldots, T_n)\in B(\cH)
$$
such that $z_i\mapsto T_i$, \ $i=1,\ldots,n$.
Due to a result of V.~M\" uller \cite{M}, the condition that
$\sigma(T_1,\ldots, T_n)\subset \BB_n$ is equivalent to the fact that the joint spectral radius
$$r(T_1,\ldots, T_n):=\lim_{k\to\infty}\left\|\sum_{|\alpha|=k} T_\alpha T_\alpha^*\right\|^{1/2k}<1.
$$

 F.H.~Vasilescu introduced and studied,  in \cite{Va} and  a joint paper with R.E.~Curto \cite{CV},  operator-valued  Cauchy and Poisson transforms on the unit ball $\BB_n$  associated with commuting operators with $r(T_1,\ldots, T_n)<1$,
in connection with commutative multivariable dilation theory.

In recent years, there has been exciting progress in noncommutative multivariable operator theory regarding  noncommutative dilation theory (\cite{F}, \cite{B}, \cite{Po-models}, \cite{Po-isometric}, \cite{Po-charact}, \cite{DKS}, \cite{BB}, \cite{BBD}, \cite{Po-unitary},  \cite{Po-varieties}, etc.) and its applications concerning  interpolation in several variables
(\cite{Po-charact}, \cite{Po-analytic}, \cite{Po-interpo}, \cite{ArPo2}, \cite{DP}, \cite{BTV},  \cite{Po-entropy}, etc.) and unitary invariants
for $n$-tuples of operators (\cite{Po-charact}, \cite{Arv}, \cite{Arv2}, \cite{Po-curvature}, \cite{Kr}, \cite{Po-similarity}, \cite{BT},
 \cite{Po-entropy}, \cite{Po-unitary}, etc.).

Our program to develop a {\it free} analogue of Sz.-Nagy--Foia\c s theory
\cite{SzF-book}, for row contractions, fits perfectly  the setting of the present paper, which includes that of free holomorphic  functions on  the open
operatorial unit ball
$$
[B(\cH)^n]_1:=\left\{ [X_1,\ldots, X_n]\in B(\cH)^n: \ \|X_1X_n^*+\cdots + X_nX_n^*\|<1\right\}.
$$
The present work is an attempt  to develop a theory of holomorphic functions in several noncommuting (free) variables and thus provide a framework for  the study of arbitrary $n$-tuples of operators,  and to introduce
and study a free analytic functional calculus  in connection with Hausdorff derivations,  noncommutative Cauchy and Poisson transforms, and von Neumann inequalities.

In Section 1, we introduce a notion of radius of convergence  for formal power series in $n$ noncommuting indeterminates $Z_1,\ldots, Z_n$ and prove noncommutative multivariable analogues of Abel theorem and Hadamard formula from complex analysis (\cite{Co}, \cite{R}). This enables us to define, in particular,
the algebra $Hol(B(\cX)^n_1)$~ of free holomorphic functions on the open operatorial unit $n$-ball, as the set of all power series  $\sum_{\alpha\in \FF_n^+}a_\alpha Z_\alpha$ with radius of convergence $\geq 1$. When $n=1$, \, $Hol(B(\cX)^1_1)$  coincides with the algebra of all analytic functions on the open unit disc ~$\DD:=\{z\in \CC:\ |z|<1\}$. The algebra of free holomorphic functions  ~$Hol(B(\cX)^n_1)$~ has the following universal property.

{\it Any   strictly contractive representation  $\pi: \CC[Z_1,\ldots, Z_n]\to B(\cH)$, i.e., $\|[\pi(Z_1),\ldots, \pi(Z_n)]\|<1$,  extends uniquely
to  a representation of $Hol(B(\cX)^n_1)$.}

A free  holomorphic function on the open operatorial unit ball ~$[B(\cH)^n]_1$~
is the representation of an element $F\in Hol(B(\cX)^n_1)$ on the Hilbert space $\cH$, that is,  the mapping
$$
[B(\cH)^n]_1\ni (X_1,\ldots, X_n)\mapsto F(X_1,\ldots, X_n)\in B(\cH).
$$
As expected, we prove that  any  free holomorphic function  is continuous on $[B(\cH)^n]_1$ in the operator norm topology. In the last part of this section, we show that  the Hausdorff derivations $\frac{\partial}{\partial Z_i}$, \, $i=1,\ldots, n$, on the algebra of noncommutative polynomials  $\CC[Z_1,\ldots, Z_n]$ (\cite{MKS}, \cite{RSS}) can be  extended to the algebra of free holomorphic functions.

In Section 2, we obtain Cauchy type estimates for the coefficients of free holomorphic functions  and a Liouville type theorem for free entire functions. Based on a noncommutative version of Gleason's problem
 \cite{R2}, which is obtained here,  and the noncommutative von Neumann inequality \cite{Po-von}, we provide a free analogue of Schwartz lemma from complex analysis (\cite{Co}, \cite{R}).
In particular, we prove that if $f$ is a free holomorphic function on
$[B(\cH)^n]_1$ such that $\|f\|_\infty\leq 1$ and  $f(0)=0$, then
$$
\|f(X_1,\ldots, X_n)\|\leq \|[X_1,\ldots, X_n]\|,\qquad r(f(X_1,\ldots, X_n))\leq r(X_1,\ldots, X_n),
$$
and \ $\sum_{i=1}^n \left|\frac{\partial f}{\partial X_i}(0)\right|^2\leq 1$.

In Section 3, following the classical case (\cite{H}, \cite{RR}), we
introduce two  Banach algebras of free holomorphic functions,
$H^\infty(B(\cX)^n_1)$ and $A(B(\cX)^n_1)$, and prove that, together
with a natural operator space structure,  they are completely
isometrically isomorphic to the noncommutative analytic Toeplitz
algebra $F_n^\infty$ and the noncommutative disc algebra $\cA_n$,
respectively, which were introduced in \cite{Po-von} in connection
with a multivariable von Neumann inequality. We recall that the
algebra $F_n^\infty$ (resp. $\cA_n$) is the weakly (resp. norm)
closed algebra generated by the left creation  operators
$S_1,\ldots, S_n$ on the full Fock space with $n$ generators,
$F^2(H_n)$,  and the identity. These algebras have been intensively
studied in recent years by many authors (\cite{Po-charact},
\cite{Po-multi}, \cite{Po-von}, \cite{Po-funct}, \cite{Po-analytic},
\cite{Po-disc}, \cite{Po-poisson}, \cite{Po-curvature},
\cite{Po-similarity}, \cite{ArPo}, \cite{ArPo2}, \cite{DP1},
\cite{DP2}, \cite{DKP}, \cite{PPoS}, \cite{Po-unitary}). The results
of this section are used to obtain a maximum principle for free
holomorphic functions.

In Section 4,  we provide a  free analytic functional calculus
for $n$-tuples $T:=[T_1,\ldots, T_n]\in B(\cH)^n$  of operators
with joint spectral radius
$r(T_1,\ldots, T_n)<1$. We show that there is a continuous  unital  algebra homomorphism
$$
\Phi_T:Hol(B(\cX)^n_1)\to B(\cH), \quad \Phi_T(f)=f(T_1,\ldots, T_n),
$$
  which is  uniquely determined by the mapping $z_i\mapsto T_i$, \ $i=1,\ldots,n$. (The continuity  and uniqueness of $\Phi_T$ are  proved in Section 5.)  We introduce a noncommutative Cauchy transform
$\cC_T:B(F^2(H_n))\to B(\cH)$ associated with any $n$-tuple of operators with joint spectral radius $r(T_1,\ldots, T_n)<1$.
The definition is  based on the {\it  reconstruction operator}
$$
S_1\otimes T_1^*+\cdots + S_n\otimes T_n^*,
$$
which has played an important role in noncommutative multivariable operator
theory (\cite{Po-entropy}, \cite{Po-unitary}, \cite{Po-varieties}).
We prove  that
$$
f(T_1,\ldots, T_n)=C_T(f(S_1,\ldots, S_n)),\quad f\in H^\infty (B(\cX)^n_1),
$$
where  $f(S_1,\ldots, S_n)$ is the boundary function of $f$.
Hence, we deduce that
$$
\|f(T_1,\ldots, T_n)\|\leq M \|f\|_\infty
$$
where $M:=\sum_{k=0}^\infty \left\|\sum_{|\alpha|=k} T_\alpha T_\alpha^*\right\|^{1/2}$.
Similar Cauchy representations are obtained for the $k$-order Hausdorff derivations of $f$. Finally, we show  that the noncommutative Cauchy transform commutes with the action of the unitary group $\cU(\CC^n)$. More precisely, we prove that
$$
\cC_T(\beta_U(f))=\cC_{\beta_U(T)}(f)\quad \text{ for any } \ U\in \cU(\CC^n), \, f\in \cA_n,
$$
where  $\beta_U$ denotes  a natural isometric automorphism (generated by $U$) of  the noncommutative disc algebra $\cA_n$, or the open unit ball $[B(\cH)^n]_1$.

In Section 5, we obtain Weierstrass and Montel type theorems \cite{Co} for the algebra  of free holomorphic functions on the open operatorial unit $n$-ball.
This enables us to introduce a metric on $Hol(B(\cX)^n_1)$ with respect to which it becomes a complete metric space, and
the Hausdorff derivations
$$
\frac{\partial}{\partial Z_i}:Hol(B(\cX)^n_1)\to Hol(B(\cX)^n_1),\quad i=1,\ldots,n,
$$
are continuous. In the end of this section, we prove  the continuity  and uniqueness of the free functional calculus  for $n$-tuples of operators with joint spectral radius $r(T_1,\ldots, T_n)<1$.
Connections with  the $F_n^\infty$-functional calculus for row contractions \cite{Po-funct} and, in the commutative case, with Taylor's functional calculus \cite{T2} are  also discussed.

Given an operator $A\in B(F^2(H_n))$, the noncommutative Poisson transform \cite{Po-poisson} generates a function
$$
P[A]: [B(\cH)^n]_1\to B(\cH).
$$
In Section 6,  we provide classes of operators $A\in B(F^2(H_n))$
such that $P[A]$ is a free holomorphic (resp. pluriharmonic)
function on $[B(\cH)^n]_1$. In this case, the operator $A$ can be
regarded as the boundary function of the Poisson extension $P[A]$.
Using some results from \cite{Po-von}, \cite{Po-funct},  and
\cite{Po-poisson}, we characterize the free holomorphic functions
$u$ on  the open unit ball $[B(\cH)^n]_1$ such that $u=P[f]$ for
some boundary function $f$ in the noncommutative analytic Toeplitz
algebra $F_n^\infty$, or the noncommutative disc algebra $\cA_n$.
For example, we prove that there exists $f\in F_n^\infty$ such that
$u=P[f]$ if and only if
$$
\sup_{0\leq r<1}\|u(rS_1,\ldots, rS_n)\|<\infty.
$$
We also  obtain noncommutative multivariable versions of Herglotz
theorem  and Dirichlet extension problem (\cite{Co}, \cite{H}) for
free pluriharmonic functions.

In Section 7, we define the radial  maximal Hardy space
$H^p(B(\cX)^n_1)$, $p\geq 1$,  as the set of all free holomorphic function $F$ such that
$$
\|F\|_p:=\left( \int_0^1\|F(rS_1,\ldots, rS_n)\|^p dr \right)^{1/p}<\infty,
$$
and  prove that it is a Banach space. Moreover,   we show that
$$
\|f(T_1,\ldots, T_n)\|\leq \frac{1}{(1-\|[T_1,\ldots, T_n]\|)^{1/p}} \|f\|_p
$$
for any $[T_1,\ldots, T_n]\in [B(\cH)^n]_1$ and  $f\in H^p(B(\cX)^n_1)$.

Finally,  we introduce  the symmetrized Hardy space
$H^\infty_{\text{\rm sym}}(\BB_n)$ as the set of all holomorphic function on $\BB_n$ such that
 $
\|f\|_{\text{\rm sym}}:= \|f_{\text{\rm sym}}\|_\infty<\infty,
$
where $f_{\text{\rm sym}}\in Hol(B(\cX)^n_1)$ is the symmetrized functional calculus of $f\in Hol(\BB_n)$.
We prove that $H^\infty_{\text{\rm sym}}(\BB_n)$ is a Banach space and
$$
\|f(T_1,\ldots, T_n)\|\leq M \|f_{\text{\rm sym}}\|_\infty,
$$
for any commuting $n$-tuple of operators with $r(T_1,\ldots, T_n)<1$.

Several classical results from complex analysis are extended to our noncommutative multivariable setting. The present paper exhibits, in particular,  a ``very good''  free analogue  of the algebra of analytic functions on the open unit disc
$\DD$. This claim is also supported by the fact that numerous results  in noncommutative multivariable  operator theory (\cite{Po-von}, \cite{Po-funct}, \cite{Po-disc},
\cite{Po-poisson}, \cite{Po-unitary})
fit perfectly our setting and can be seen in a new light.
We strongly believe that many other results in the theory of analytic functions have free analogues in our noncommutative multivariable setting.

In a forthcoming paper \cite{Po-Bohr}, we  consider  operator-valued
 Wiener and Bohr type inequalities for   free
holomorphic  (resp. pluriharmonic) functions on the open operatorial
unit $n$-ball. As consequences, we obtain operator-valued   Bohr
inequalities for the noncommutative Hardy algebra
$H^\infty(B(\cX)^n_1)$ and the symmetrized Hardy space
$H^\infty_{\text{\rm sym}}(\BB_n)$.

\bigskip

      \section{Free holomorphic functions }
\label{free holomorphic}

 We introduce a notion of radius of convergence  for formal power series in $n$ noncommuting indeterminates $Z_1,\ldots, Z_n$ and prove noncommutative multivariable analogues of Abel theorem and Hadamard formula. This enables us to define
 algebras  of free holomorphic functions on open operatorial  $n$-balls.
We show that the Hausdorff derivations $\frac{\partial}{\partial Z_i}$, \, $i=1,\ldots, n$, on the algebra of noncommutative polynomials  $\CC[Z_1,\ldots, Z_n]$ (see \cite{MKS}, \cite{RSS}) can be  extended to  algebras of free holomorphic functions.

Let $\FF_n^+$ be the unital free semigroup on $n$ generators $g_1,\ldots, g_n$ and the identity $g_0$.  The length of $\alpha\in \FF_n^+$ is defined by $|\alpha|=0$ if $\alpha=g_0$  and $|\alpha|:=k$ if $\alpha=g_{i_1}\cdots g_{i_k}$, where $i_1,\ldots, i_k\in \{1,\ldots, n\}$.
  We consider formal power series in $n$ noncommuting indeterminates $Z_1,\ldots, Z_n$ and coefficients in $B(\cK)$,
 the algebra of all bounded linear operators on the Hilbert space $\cK$,
 of the  form
\begin{equation*}
\sum_{\alpha\in \FF_n^+} A_{(\alpha)}\otimes Z_\alpha,\quad A_{(\alpha)}\in B(\cK),
\end{equation*}
where $Z_\alpha:=Z_{i_1}\cdots Z_{i_k}$ if $\alpha=g_{i_1}\cdots g_{i_k}$ and $Z_{g_0}:=I$.
If $F=\sum_{\alpha\in \FF_n^+} A_{(\alpha)}\otimes Z_\alpha $
and
$G=\sum_{\alpha\in \FF_n^+} B_{(\alpha)}\otimes Z_\alpha $
are such formal power series, we define their sum and product
by   setting
$$
F+G:=\sum_{\alpha\in \FF_n^+} (A_{(\alpha)}+B_{(\alpha)})\otimes Z_\alpha \quad
\text{ and }\quad
FG:=\sum_{\alpha\in \FF_n^+} C_{(\alpha)}\otimes Z_\alpha,
$$
respectively,
where $C_{(\alpha)}:=\sum\limits_{\sigma, \beta\in \FF_n^+:\  \alpha=\sigma \beta} A_{(\sigma)} B_{(\beta)}$.

By abuse of notation, throughout this paper, we will denote by $[T_1,\ldots,T_n]$ either  the $n$-tuple of operators
$(T_1,\ldots, T_n)\in B(\cH)^n$ or the row operator matrix
$[T_1\,\cdots \,T_n]\in B(\cH^{(n)}, \cH)$ acting as an operator  from $\cH^{(n)}$, the direct sum of $n$ copies of the Hilbert space  $\cH$, to $\cH$. We also denote
by $[T_\alpha:\ |\alpha|=k]$ the row operator matrix acting from $\cH^{n^k}$ to $\cH$, where the entries are arranged  in the lexicographic order of the free semigroup $\FF_n^+$.

 In what follows  we show that given a sequence of operators
$A_{(\alpha)}\in B(\cK)$, $\alpha\in \FF_n^+$, there is a unique $R\in [0,\infty]$ such that the series
\begin{equation*}
\sum_{k=0}^\infty\sum_{|\alpha|=k} A_{(\alpha)}\otimes X_\alpha
\end{equation*}
converges in the operator norm of $B(\cK\otimes \cH)$ ($\cK\otimes \cH $ is the Hilbert tensor product) for any Hilbert space $\cH$ and any $n$-tuple  $[X_1,\ldots, X_n]\in B(\cH)^n$ with $\|[X_1,\ldots, X_n]\|<R$, and it is divergent for some   $n$-tuples $[Y_1,\ldots, Y_n]$ of operators  with
$\|[Y_1,\ldots, Y_n]\|>R$.

The result can be regarded as a noncommutative multivariable analogue of Abel theorem and Hadamard's formula from complex analysis.

\begin{theorem}\label{Abel} Let $\cH$, $\cK$ be Hilbert spaces and let $A_{(\alpha)}\in B(\cK)$, $\alpha\in \FF_n^+$, be a sequence of operators.
 Define
 $R\in [0,\infty]$ by setting
$$
\frac {1} {R}:=
\limsup_{k\to\infty}
\left\|\sum_{|\alpha|=k} A^*_{(\alpha)} A_{(\alpha)}\right\|^{\frac{1} {2k}}.
$$
Then the following properties hold:
\begin{enumerate}
\item [(i)]
For any  $n$-tuple of operators
$[X_1,\ldots, X_n]\in B(\cH)^n$, the series
$\sum\limits_{k=0}^\infty \left\| \sum\limits_{|\alpha|=k}A_{(\alpha)}\otimes X_\alpha\right\|
$ converges if \ $\|[X_1,\ldots,X_n]\|<R$. Moreover, if $0\leq \rho<R$, then the convergence is uniform  for
$[X_1,\ldots, X_n]$ with $\|[X_1,\ldots,X_n]\|\leq \rho$.
\item [(ii)]
If  $R<R'<\infty$ and $\cH$ is infinite dimensional, then there is
an  $n$-tuple  $[X_1,\ldots, X_n]\in B(\cH)^n$ of operators with
$$\|X_1X_1^*+\cdots +X_nX_n^*\|^{1/2}=R'
$$
such that
$\sum\limits_{k=0}^\infty \left( \sum\limits_{|\alpha|=k}A_{(\alpha)}\otimes X_\alpha\right)$
is divergent in the operator norm of $B(\cK\otimes \cH)$.
\end{enumerate}

Moreover, the number $R$ satisfying properties (i) and (ii) is unique.
\end{theorem}
\begin{proof}
Assume that  $R>0$ and  $[X_1,\ldots, X_n]$ is an $n$-tuple of operators on $\cH$ such that $\|[X_1,\ldots,X_n]\|<R$.
Let $\rho',\rho>0$ be such that
$\|[X_1,\ldots,X_n]\|<\rho'<\rho<R$.
 Since $\frac{1}{\rho}> \frac{1}{R}$, we can find $m_0\in\NN:=\{1,2,\ldots\}$ such that
\begin{equation*}
\left\|\sum\limits_{|\alpha|=k}A_{(\alpha)}^*A_{(\alpha)}\right\|^{1/2k}< \frac{1}{\rho}\quad  \text{ for any }\ k\geq m_0.
\end{equation*}
Hence, we deduce that
 \begin{equation*}
\begin{split}
\left\|\sum\limits_{|\alpha|=k}A_{(\alpha)}\otimes X_\alpha \right\|&=
\left\|\left[ I\otimes X_\alpha:\ |\alpha|=k\right]
\left[\begin{matrix}
A_{(\alpha)}\otimes I\\
:\\|\alpha|=k
\end{matrix}\right]\right\|\\
&=\left\|\sum\limits_{|\alpha|=k}X_{\alpha}X_\alpha^*\right\|^{1/2}\left\|\sum\limits_{|\alpha|=k}A_{(\alpha)}^*A_{(\alpha)}\right\|^{1/2}\\
&\leq \left\|\sum_{i=1}^nX_iX_i^*\right\|^{k/2}\left\|\sum\limits_{|\alpha|=k}A_{(\alpha)}^*A_{(\alpha)}\right\|^{1/2}\\
&\leq \left(\frac{\rho'}{\rho}\right)^k
\end{split}
\end{equation*}
for any $k\geq m_0$.
This proves the convergence of the series
$\sum\limits_{k=0}^\infty \left\| \sum\limits_{|\alpha|=k}A_{(\alpha)}\otimes X_\alpha\right\|$.
Assume now  that $0\leq \rho<R$ and $\|[X_1,\ldots, X_n]\|\leq \rho$.
 Choose $\gamma$ such that $0\leq \rho< \gamma<R$ and notice that, due to similar calculations as above,
   there exists $n_0\in \NN$ such that
\begin{equation*}
\left\|\sum\limits_{|\alpha|=k}A_{(\alpha)}\otimes X_\alpha \right\|
\leq \left(\frac{\rho}{\gamma}\right)^k
\end{equation*}
for any  $(X_1,\ldots, X_n)$ with $\|[X_1,\ldots, X_n]\|\leq \rho$,
and $k\geq n_0$, which proves the uniform convergence of the above
series. The case $R=\infty$, can be treated in a similar manner.

To prove part (ii), assume that $R<\infty$ and  $\cH$ is infinite dimensional. Let
$R', \rho>0$ be such that
$R<\rho< R'$ and define the operators $X_i:= R' V_i$, \ $i=1,\ldots, n$, where $V_1,\ldots, V_n$   are isometries with orthogonal ranges.
Notice that $\|[X_1,\ldots, X_n]\|=R'$ and

\begin{equation*}
\begin{split}
\left\|\sum\limits_{|\alpha|=k}A_{(\alpha)}\otimes X_\alpha \right\|&={R'}^k
\left\|\left(\sum\limits_{|\alpha|=k}A_{(\alpha)}^*\otimes V_\alpha^*\right) \left(\sum\limits_{|\alpha|=k}A_{(\alpha)}\otimes V_\alpha\right)\right\|^{1/2}\\
&=
{R'}^k
\left\|\sum\limits_{|\alpha|=k}A_{(\alpha)}^*A_{(\alpha)} \otimes I\right\|^{1/2}\\
&=
{R'}^k
\left\|\sum\limits_{|\alpha|=k}A_{(\alpha)}^* A_{(\alpha)} \right\|^{1/2}.
\end{split}
\end{equation*}
On the other hand, since $\frac{1}{\rho}<\frac{1}{R}$, there are arbitrarily large $k\in \NN$ such that
$$
\left\|\sum\limits_{|\alpha|=k}A_{(\alpha)}^* A_{(\alpha)} \right\|^{1/2}>\left(\frac{1}{\rho}\right)^k.
$$
Consequently, we deduce that
$$
\left\|\sum\limits_{|\alpha|=k}A_{(\alpha)}\otimes X_\alpha \right\|>\left(\frac{R'}{\rho}\right)^k,
$$
which proves that
 the series
$\sum\limits_{k=0}^\infty \left( \sum\limits_{|\alpha|=k}A_{(\alpha)}\otimes X_\alpha\right)$
is divergent in the operator norm. The uniqueness of the number $R$ satisfying properties (i) and (ii) is now obvious.
\end{proof}

   As expected, the number $R$ in the above theorem is called the radius of convergence of the power series
$\sum_{\alpha\in \FF_n^+} A_{(\alpha)}\otimes Z_\alpha.$

Let us consider the full Fock space
$$
F^2(H_n)=\CC 1\oplus\ \oplus_{m\ge1}H_n^{\otimes m}
$$
where $H_n$ is an $n$-dimensional complex Hilbert space with orthonormal basis
$\{e_1,\dots,e_n\}$.
Setting  $e_\alpha:=e_{i_1}\otimes\cdots e_{i_k}$ if $\alpha=g_{i_1}\cdots g_{i_k}$,  and $e_{g_0}=1$, it is clear that $\{ e_\alpha:\ \alpha\in \FF_n^+\}$ is an orthonormal basis of the full Fock space $F^2(H_n)$.
For each $i=1,2,\dots$,  we define  the left creation operator
$\ S_i\in B(F^2(H_n))$  by
 $$S_i\xi=e_i\otimes\xi,\qquad  \xi\in F^2(H_n).
$$

We can now obtain the following characterization of the radius of convergence, which will be useful later.

\begin{corollary}\label{Cs} Let
  $\sum\limits_{\alpha\in \FF_n^+}A_{(\alpha)}\otimes Z_\alpha$ be a formal power series with radius of convergence $R$.

\begin{enumerate}
\item[(i)]
If  $R>0$ and $0<r<R$, then  there exists $C>0$ such that
$$
\left\|\sum\limits_{|\alpha|=k}A_{(\alpha)}^* A_{(\alpha)} \right\|^{1/2}\leq \frac {C}{r^k}\quad
\text{for any } \ k=0,1,\ldots.
$$
\item[(ii)] The radius of convergence of the power series satisfies the relations
$$
R=\sup\left\{ r\geq 0:\ \text{ the sequence }\ \left\{r^k\left\|\sum\limits_{|\alpha|=k}A_{(\alpha)}^* A_{(\alpha)} \right\|^{1/2} \right\}_{k=0}^\infty  \text{ is bounded }\right\}
$$
 and
$$
R=\sup\left\{ r\geq 0: \ \sum_{k=0}^\infty \sum_{|\alpha|=k}r^{|\alpha|} A_{(\alpha)}\otimes S_\alpha\
 \text{ is convergent in the operator norm }\right\}.
$$
\end{enumerate}
\end{corollary}
\begin{proof}
Setting $X_i:= rS_i$, \ $i=1,\ldots, n$, where $S_1,\ldots, S_n$ are the left creation operators on the full Fock space, we have
$\|[X_1,\ldots, X_n]\|=r<R$. According to Theorem \ref{Abel},
the series
$
\sum_{k=0}^\infty \left\|r^k\sum_{|\alpha|=k} A_{(\alpha)} \otimes S_\alpha\right\|$ is convergent.
Since $S_1,\ldots, S_n$ are isometries with orthogonal ranges, the above series is equal to
$\sum_{k=0}^\infty r^k \left\|\sum_{|\alpha|=k}A_{(\alpha)}^* A_\alpha \right\|^{1/2}.
$
Consequently, there is a constant $C>0$ such that
$$
r^k \left\|\sum\limits_{|\alpha|=k}A_{(\alpha)}^* A_\alpha \right\|^{1/2}\leq C\ \text{ for any } k=0,1,\ldots.
$$
Now, the second part of this corollary follows easily from part (i) and Theorem \ref{Abel}.
This completes the proof.
\end{proof}

We establish terminology which will be used throughout the paper.
    Denote by $[B(\cH)^n]_{\gamma}$ the open ball of $B(\cH)^n$  of radius $\gamma> 0$, i.e.,
$$
[B(\cH)^n]_{\gamma}:=\{[X_1,\ldots, X_n]:\
\|X_1X_1^*+\cdots +X_nX_n^*\|^{1/2}<\gamma\}.
$$
We also use the notation $[B(\cH)^n]_1^-$ for the closed ball.
A formal power series $F:=\sum\limits_{\alpha\in \FF_n^+} A_{(\alpha)}\otimes Z_\alpha$
represents a free holomorphic function
 on  the open  operatorial $n$-ball of radius $\gamma$ with coefficients in $B(\cK)$,   if  for any Hilbert space $\cH$ and  any representation
$$
\pi:\CC[Z_1,\ldots, Z_n]\to B(\cH)\quad \text{  such that } \quad
[\pi(Z_1),\ldots, \pi(Z_n)]\in [B(\cH)^n]_{\gamma}
$$
 the series
$$F(\pi(Z_1),\ldots, \pi(Z_n)):=\sum\limits_{k=0}^\infty \sum\limits_{|\alpha|=k} A_{(\alpha)} \otimes \pi(Z_\alpha)$$
 converges in the operator norm of $B(\cK\otimes \cH)$.
Due to  Theorem \ref{Abel}, we must have
$\gamma\leq R$, where $R$ is the radius of convergence of $F$.
The mapping
$$
[B(\cH)^n]_{\gamma}\ni [X_1,\ldots, X_n]\mapsto F(X_1,\ldots X_n)\in B(\cK\otimes \cH).
$$
is called the representation of $F$ on the Hilbert space $\cH$.
 Given a Hilbert space $\cH$,
we say that a function  $G:[B(\cH)^n]_{\gamma}\to B(\cK  \otimes \cH)$ is a {\it free holomorphic  function} on  $[B(\cH)^n]_{\gamma
}$ with coefficients in $B(\cK)$ if there exist operators  $A_{(\alpha)}\in B(\cK)$, $\alpha\in \FF_n^+$, such that
the power series $\sum_{\alpha\in \FF_n^+} A_{(\alpha)}\otimes Z_\alpha$ has radius of convergence
 $\geq \gamma$ and
$$
G(X_1,\ldots, X_n)=\sum\limits_{k=0}^\infty \left( \sum\limits_{|\alpha|=k}A_{(\alpha)}\otimes X_\alpha\right),
$$
where the series converges in the operator  norm for any
$[X_1,\ldots, X_n]\in [B(\cH)^n]_{\gamma}$.

We remark  that  the coefficients  of a  free holomorphic
function  are uniquely determined by its representation
on
   an infinite dimensional   Hilbert space. Indeed,  let $0<r<\gamma$ and assume $F(rS_1,\ldots, rS_n)=0$, where $S_1,\ldots, S_n$ are the left creation operators on the full Fock space $F^2(H_n)$. Taking into account that   $S_i^* S_j=\delta_{ij} I$, we have
\begin{equation*}
\left< F(rS_1,\ldots, rS_n)(x\otimes 1), (I_\cK\otimes S_\alpha)(y\otimes 1)\right>=\left<A_{(\alpha)}x,y\right>=0
\end{equation*}
for any $x,y\in \cK$ and $\alpha\in \FF_n^+$. Therefore $A_{(\alpha)}=0$ for any $\alpha\in \FF_n^+$.

 We establish now the continuity of free holomorphic functions on the open  ball  $[B(\cH)^n]_{\gamma}$.

\begin{theorem}\label{continuous}
Let  $ f(X_1,\ldots, X_n)=\sum\limits_{k=0}^\infty \left( \sum\limits_{|\alpha|=k}A_{(\alpha)}\otimes X_\alpha\right)
$
 be a free holomorphic  function on $[B(\cH)^n]_{\gamma}$ with coefficients in $B(\cK)$.
If $X:=[X_1,\ldots, X_n]$, $Y:=[Y_1,\ldots, Y_n]$ are in the closed ball \, $[B(\cH)^n]_r^-$, $0<r<\gamma$, then
$$
\|f(X)-f(Y)\|\leq
\|X-Y\|\sum _{k=1}^\infty kr^{k-1}\left\|\sum_{|\alpha|=k} A_{(\alpha)}^* A_{(\alpha)}\right\|^{1/2}.
$$
In particular,  $f$ is continuous on
$[B(\cH)^n]_{\gamma}$  and uniformly continuous on $[B(\cH)^n]_r^-$ in the operator  norm topology.
\end{theorem}

\begin{proof}
Let $X^{[k]}:=[X_\alpha:\ \alpha\in \FF_n^+, \ |\alpha|=k]$, \ $k=1,2,\ldots$, be the row operator matrix with entries arranged in the lexicographic order of the free semigroup $\FF_n^+$.
First, we prove that if $\|X\|\neq \|Y\|$, then
\begin{equation}
\label{[k]}
\frac{\|X^{[k]}-Y^{[k]}\|}{\|X-Y\|}\leq
\frac{\|X\|^k-\|Y\|^k}{\|X\|-\|Y\|}.
\end{equation}
Notice that
\begin{equation*}
\begin{split}
X^{[k]}-Y^{[k]}&=
\left[(X_1-Y_1)X^{[k-1]},\ldots, (X_n-Y_n) X^{[k-1]}\right]\\
&\qquad +
\left[ Y_1(X^{[k-1]}-Y^{[k-1]}),\ldots, Y_n(X^{[k-1]}-Y^{[k-1]})\right]\\
&=
(X-Y)\text{\rm diag}_n(X^{[k-1]})+Y\text{\rm diag}_n(X^{[k-1]}-Y^{[k-1]}),
\end{split}
\end{equation*}
where $\text{\rm diag}_n(A)$  is the $n\times n$ block diagonal operator matrix with $A$ on the diagonal and $0$ otherwise.
Hence, we deduce that
$$
\|X^{[k]}-Y^{[k]}\|\leq \|X-Y\|\|X^{[k-1]}\|+\|Y\|
\|X^{[k-1]}-Y^{[k-1]}\|
$$
for any $k\geq 2$. Iterating this relation and taking into account that $\|X^{[k]}\|\leq \|X\|^k$ for  $k=1,2,\ldots$,  we obtain
\begin{equation*}
\begin{split}
\|X^{[k]}-Y^{[k]}\|&\leq \|X-Y\|\left(\|X^{[k-1]}\|
+\|\|X^{[k-2]}\|\|Y^{[1]}\|+\cdots + \|Y^{[k-1]}\|\right)\\
&\leq \|X-Y\|\left(\|X\|^{k-1}+\|X\|^{k-2}\|Y\|+\cdots+ \|Y\|^{k-1}\right),
\end{split}
\end{equation*}
which proves  inequality \eqref{[k]}.
Assuming that $\|X\|\leq r$ and $\|Y\|\leq r$, we deduce that
\begin{equation*}
\|X^{[k]}-Y^{[k]}\|\leq kr^{k-1}\|X-Y\|,\quad k=1,2,\ldots.
\end{equation*}
Hence, we obtain
\begin{equation*}
\begin{split}
\|f(X)-f(Y)\|&\leq \sum_{k=1}^\infty\left\|\sum_{|\alpha|=k} A_{(\alpha)}\otimes (X_\alpha-Y_\alpha)\right\|\\
&\leq \sum_{k=1}^\infty\left\|\sum_{|\alpha|=k} A_{(\alpha)}^* A_{(\alpha)}\right\|^{1/2} \|X^{[k]}-Y^{[k]}\|\\
&\leq \|X-Y\|
\sum_{k=1}^\infty kr^{k-1} \left\|\sum_{|\alpha|=k} A_{(\alpha)}^* A_{(\alpha)}\right\|^{1/2}.
\end{split}
\end{equation*}
Let $\rho$ be a constant such that $r<\rho<\gamma$.
Since $\gamma\leq R$ ($R$ is the radius of convergence of $f$) and $\frac{1}{\rho}>\frac{1}{\gamma}\geq \frac{1}{R}$,
we can find $m_0\in \NN$, such that
$$
 \left\|\sum_{|\alpha|=k} A_{(\alpha)}^* A_{(\alpha)}\right\|^{1/2k}<\frac{1}{\rho}\quad \text{ for any
} \ k\geq m_0.
$$
Combining this with the above inequality, we deduce that
$$
 \|f(X)-f(Y)\|\leq \|X-Y\|\left(\sum_{k=1}^{m_0-1} kr^{k-1}
\left\|\sum_{|\alpha|=k} A_{(\alpha)}^* A_{(\alpha)}\right\|^{1/2}+ \sum_{k=m_0}^\infty \frac{k}{r} \left(\frac{r}{\rho}\right)^k\right).
$$
Since $r<\rho$, the above series is convergent. Consequently, there exists a constant $M>0$ such that
$$
\|f(X)-f(Y)\|\leq M \|X-Y\|\qquad \text{ for any }\ X,Y\in [B(\cH)^n]_r^-.
$$
This  implies the  uniform continuity of $f$  on  any closed ball $[B(\cH)^n]_r^-$, $0<r<\gamma$,  in the norm topology and,
consequently,   the continuity of $f$  on  $[B(\cH)^n]_{\gamma}$.
\end{proof}

\begin{theorem}\label{operations}
Let $F$ and $G$ be formal power series such that
\begin{equation*}
\begin{split}
F(X_1,\ldots, X_n)&=\sum\limits_{k=0}^\infty \left( \sum\limits_{|\alpha|=k}A_{(\alpha)}\otimes X_\alpha\right)\text{ and }\\
G(X_1,\ldots, X_n)&=\sum\limits_{k=0}^\infty \left( \sum\limits_{|\alpha|=k}B_{(\alpha)}\otimes X_\alpha\right)
\end{split}
\end{equation*}
are free holomorphic functions on  $[B(\cH)^n]_{\gamma}$, and let $a, b\in \CC$. Then  the power series  $aF+bG$, and $FG$ generate  free holomorphic functions
on $[B(\cH)^n]_{\gamma}$. Moreover,
\begin{equation*}
\begin{split}
 aF(X_1,\ldots, X_n)+bG(X_1,\ldots, X_n)&=\sum\limits_{k=0}^\infty \left( \sum\limits_{|\alpha|=k}(aA_{(\alpha)}+bB_{(\alpha)})\otimes X_\alpha\right) \text{ and }\\
F(X_1,\ldots, X_n)G(X_1,\ldots, X_n)&=\sum\limits_{k=0}^\infty \left( \sum\limits_{|\alpha|=k} C_{(\alpha)}\otimes X_\alpha \right)
\end{split}
\end{equation*}
for any $[X_1,\ldots, X_n]\in [B(\cH)^n]_\gamma$,
where $C_{(\alpha)}:= \sum\limits_{\alpha=\sigma\beta}A_{(\sigma)} B_{(\beta)}$, \ $\alpha\in \FF_n^+$.
\end{theorem}
\begin{proof}
According to the hypotheses, both power series $F$ and $G$
have radius of convergence $\geq \gamma$. Due to Theorem \ref{Abel},  we deduce that, given  any $\epsilon>0$, there exists $k_0\in \NN$ such that
$$
\left\|\sum_{|\alpha|=k} A_{(\alpha)}^* A_{(\alpha)}\right\|^{1/2k}\leq \frac{1}{\gamma} +\epsilon
\ \text{ and }\
\left\|\sum_{|\alpha|=k} B_{(\alpha)}^* B_{(\alpha)}\right\|^{1/2k}\leq \frac{1}{\gamma} +\epsilon
$$
for any $k\geq k_0$.
Assume that $|a|+|b|\neq 0$.
Since the left creation operators $S_1,\ldots, S_n$ are isometries with orthogonal ranges, we have
  \begin{equation*}
\begin{split}
\Biggl\|\sum_{|\alpha|=k} (aA_{(\alpha)}+bB_{(\alpha)})^* &(aA_{(\alpha)}+ bB_{(\alpha)})\Biggr\|^{1/2}\\
&=
\left\|\sum_{|\alpha|=k} (aA_{(\alpha)} +bB_{(\alpha)})\otimes S_\alpha \right\|\\
&\leq
\left\|\sum_{|\alpha|=k} aA_{(\alpha)}  \otimes S_\alpha \right\|+\left\|\sum_{|\alpha|=k} bB_{(\alpha)}\otimes S_\alpha\right\|\\
&=\left\|\sum_{|\alpha|=k} |a|^2A_{(\alpha)}^* A_{(\alpha)}\right\|^{1/2}+\left\|\sum_{|\alpha|=k} |b|^2B_{(\alpha)}^* B_{(\alpha)}\right\|^{1/2}\\
&=(|a|+|b|)\left( \frac{1}{\gamma}+\epsilon\right)^k
\end{split}
\end{equation*}
for any $k\geq k_0$. Hence, we deduce that
$$
\limsup_{k\to\infty}
\left\|\sum_{|\alpha|=k} (aA_{(\alpha)}+bB_{(\alpha)})^* (aA_{(\alpha)}+ bB_{(\alpha)})\right\|^{1/2k}\leq \frac{1}{\gamma}+\epsilon
$$
for any $\epsilon>0$. Taking $\epsilon\to 0$, we  deduce that
the power series $aF+bG$ has the radius of convergence $\geq \gamma$.
Now, we prove that the power series $FG$ has radius of
convergence $\geq \gamma$. If $0<r<\gamma$, then, due to
 Corollary \ref{Cs}, there is a constant $M>0$ such that
\begin{equation*}
\begin{split}
\left\|\sum_{|\sigma|=k} C_{(\sigma)}^* C_{(\sigma)}\right\|^{1/2}&= \left\| \sum_{|\sigma|=k}C_{(\sigma)}\otimes S_\sigma\right\|\\
&=
\left\|
\sum_{p+q=k} \left( \sum_{|\alpha|=p}A_{(\alpha)}\otimes S_\alpha\right) \left( \sum_{|\beta|=q}B_{(\beta)}\otimes S_\beta\right)\right\|\\
&\leq
\sum_{p+q=k}\left\| \sum_{|\alpha|=p}A_{(\alpha)}^* A_{(\alpha)}  \right\|^{1/2} \left\|\sum_{|\beta|=q}B_{(\beta)}^*B_{(\beta)} \right\|^{1/2}\\
&\leq
\sum_{p+q=k} \frac{M}{r^p}\cdot\frac{M}{r^q}\\
&= (k+1) \frac{M^2}{r^k}
\end{split}
\end{equation*}
for any $k=0,1,\ldots$.
Hence, we obtain
$$
\limsup_{k\to\infty}
\left\|\sum_{|\sigma|=k}  C_{(\sigma)}^* C_{(\sigma)}\right\|^{1/2k}\leq \frac {1}{r}
$$
for any $r$ such that $0<r<\gamma$.
Consequently, the radius of convergence of the power series $FG$ is $\geq \gamma$.
The last part of the theorem follows easily  using Theorem
\ref{Abel}.
\end{proof}

We are in position to give  a characterization  as well as models for free holomorphic functions on the open operatorial $n$-ball of radius $\gamma$.

\begin{theorem} \label{caract-shifts}
A power series $F:=\sum\limits_{\alpha\in \FF_n^+} A_{(\alpha)}\otimes Z_\alpha$
represents a free holomorphic function on the open  operatorial $n$-ball of radius $\gamma$ with coefficients in $B(\cK)$      if and only if
the series
$$
\sum\limits_{k=0}^\infty  \sum\limits_{|\alpha|=k} r^{|\alpha|}  A_{(\alpha)}\otimes S_\alpha
$$
is convergent for any $r\in [0,\gamma)$, where
$S_1,\ldots, S_n$ are the left creation operators on the Fock space  $F^2(H_n)$.
 Moreover, in this case, the series
\begin{equation}\label{cre-seri}
\sum\limits_{k=0}^\infty\left\| \sum\limits_{|\alpha|=k} r^{|\alpha|}  A_{(\alpha)}\otimes S_\alpha \right\|=\sum_{k=0}^\infty r^k\left\|
\sum_{|\alpha|=k} A_{(\alpha)}^*A_{(\alpha)}\right\|^{1/2}
\end{equation}
are convergent for any $r\in [0,\gamma)$.
\end{theorem}

\begin{proof} Assume that $F$ represents a free holomorphic function on the open  operatorial $n$-ball of radius $\gamma$.
According to Theorem \ref{Abel}, $\gamma\leq R$, where $R$ is the radius of convergence of  $F$, and $\sum\limits_{k=0}^\infty \left\|\sum\limits_{|\alpha|=k} A_{(\alpha)}\otimes X_\alpha\right\|$ converges for any $n$-tuple $[X_1,\ldots, X_n]$ with $\|[X_1,\ldots, X_n]\|=r<\gamma$. Since $\|[rS_1,\ldots, rS_n]\|=r<\gamma$, we deduce that the series
\eqref{cre-seri} is convergent for any $r\in [0,\gamma)$.

Now, assume that  the series
\eqref{cre-seri} is convergent for any $r\in [0,\gamma)$.
According to the noncommutative von Neumann inequality \cite{Po-von}, we have
$$
\sum\limits_{k=0}^\infty\left\| \sum\limits_{|\alpha|=k} r^{|\alpha|}  A_{(\alpha)}\otimes T_\alpha \right\|\leq
\sum\limits_{k=0}^\infty\left\| \sum\limits_{|\alpha|=k} r^{|\alpha|}  A_{(\alpha)}\otimes S_\alpha \right\|
$$
for any $n$-tuple $[T_1,\ldots, T_n]\in B(\cH)^n$ with $T_1T_1^*+\cdots T_nT_n^*\leq I$ and any $r\in [0,\gamma)$.
Hence, we deduce that the series
$$
\sum_{k=0}^\infty \left\|\sum_{|\alpha|=k} A_{(\alpha)}\otimes X_\alpha\right\|
$$
 converges for any $n$-tuple  of operators $[X_1,\ldots, X_n]$ with $\|[X_1,\ldots, X_n]\|<\gamma$.
Due to Theorem  \ref{Abel}, the power series
$F=\sum\limits_{\alpha\in \FF_n^+} A_{(\alpha)}\otimes Z_\alpha$
represents a free holomorphic function on the open  operatorial $n$-ball of radius $\gamma$. This completes the proof.
 \end{proof}

\begin{corollary}  Let $\{a_k\}_{k=0}^\infty$ be a sequence of complex numbers. Then the following statements are equivalent:
\begin{enumerate}
\item[(i)]$f(z):=\sum_{k=0}^\infty a_k z^k$ is an analytic function on the open unit disc $\DD:=\{z\in \CC:\ |z|<1\}$.
\item[(ii)] $f_r(S):=\sum_{k=0}^\infty r^ka_k S^k$ is convergent in the operator norm for each $r\in [0,1)$, where $S$ is the unilateral shift on the Hardy space $H^2$.
\item[(iii)]
$f(Z):=\sum_{k=0}^\infty a_k Z^k$ is a free holomorphic function  on the open  operatorial unit $1$-ball.
\end{enumerate}
\end{corollary}
\begin{proof}
If $f(z)=\sum\limits_{k=0}^\infty a_k z^k$ is  an analytic  function on the open unit disc, then Hadamard's theorem implies
$\limsup\limits_{k\to\infty} |a_k|^{1/k}\leq 1$. Hence $\sum\limits_{k=0}^\infty r^k|a_k|<\infty$ for any $r\in [0,1)$ and, consequently, the series  $\sum\limits_{k=0}^\infty r^k a_k S^k$ is convergent in the operator norm. Conversely, if the latter series is norm convergent, then, due to von Neumann inequality \cite{vN},  the series $\sum\limits_{k=0}^\infty r^k a_k z$ converges  for any $r\in[0,1)$ and $z\in \DD$. Hence, we deduce (i). The equivalence
(ii)$\Longleftrightarrow$ (iii) is a particular case of
Theorem \ref{caract-shifts}.
\end{proof}

If $\lambda:=(\lambda_1,\ldots, \lambda_n)\in\CC^n$ and $\alpha=g_{i_1}\cdots g_{i_k}\in \FF_n^+$, then we set $\lambda_\alpha:=\lambda_{i_1}\cdots \lambda_{i_k}$ and $\lambda_0=1$.

\begin{corollary}\label{part-case}
If $f =\sum_{\alpha\in \FF_n^+} a_\alpha Z_\alpha$, \, $a_\alpha\in \CC$,  is a free holomorphic function on the open  operatorial unit $n$-ball, then its representation on $\CC$,
$$f(\lambda_1,\ldots, \lambda_n)=\sum_{k=0}^\infty \sum_{|\alpha|=k} a_\alpha \lambda_\alpha,$$  is a holomorphic function on  $\BB_n$,  the open unit ball of $\CC^n$.
\end{corollary}

\begin{proof}
Due to Theorem  \ref{caract-shifts},  we have
\begin{equation*}
\begin{split}
\sum_{k=0}^\infty\sum_{|\alpha|=k}|a_\alpha||\lambda_\alpha|&\leq
\sum_{k=0}^\infty\left(\sum_{|\alpha|=k}|a_\alpha|^2\right)^{1/2}\left(\sum_{|\alpha|=k}|\lambda_\alpha|^2\right)^{1/2}\\
&\leq
\sum_{k=0}^\infty\left(\sum_{|\alpha|=k}|a_\alpha|^2\right)^{1/2}\left( \sum_{i=1}^n |\lambda_i|^2\right)^{k/2}<\infty
\end{split}
\end{equation*}
for any $(\lambda_1,\ldots, \lambda_n)\in \BB_n$.
Hence, the result follows.
\end{proof}

\smallskip

In the last part of this section, we show that the Hausdorff derivations on the algebra of noncommutative polynomials $\CC[Z_1,\ldots, Z_n]$ (see \cite{MKS}, \cite{RSS}) can be extended to the algebra of free holomorphic functions.
 For each $i=1,\ldots, n$,  we define the free  partial derivation  $\frac{\partial } {\partial Z_i}$  on $\CC[Z_1,\ldots, Z_n]$ as the unique linear operator  on this algebra, satisfying the conditions
$$
\frac{\partial I} {\partial Z_i}=0, \quad  \frac{\partial Z_i} {\partial Z_i}=I, \quad  \frac{\partial Z_j} {\partial Z_i}=0\ \text{ if }  \ i\neq j,
$$
and
$$
\frac{\partial (fg)} {\partial Z_i}=\frac{\partial f} {\partial Z_i} g +f\frac{\partial g} {\partial Z_i}
$$
for any  $f,g\in \CC[Z_1,\ldots, Z_n]$ and $i,j=1,\ldots n$.
The same definition extends to formal power series in the noncommuting indeterminates $Z_1,\ldots, Z_n$.

Notice that if $\alpha=g_{i_1}\cdots g_{i_p}$, $|\alpha|=p$, and $q$ of the  $g_{i_1},\ldots, g_{i_p}$ are equal to $g_j$, then  $\frac{\partial Z_\alpha} {\partial Z_j}$ is the sum
of the $q$ words obtained by deleting each occurence of $Z_j$ in $Z_\alpha:=Z_{i_1}\cdots Z_{i_p}$. For example,
$$
\frac{\partial (Z_1 Z_2 Z_1^2)} {\partial Z_1}=
Z_2 Z_1^2+ Z_1Z_2Z_1+ Z_1Z_2Z_1.
$$
One can easily show that $\frac{\partial } {\partial Z_i}$
coincides with the Hausdorff  derivative.
If $\beta:=g_{i_1}\cdots g_{i_k}\in \FF_n^+$,\ $i_1,\ldots, i_k\in \{1,2,\ldots, n\}$,
we denote $Z_\beta:=Z_{i_1}\cdots Z_{i_k}$ and  define the $k$-order free partial derivative of $G\in \CC[Z_1,\ldots, Z_n]$ with respect to $Z_{i_1},\ldots, Z_{i_k}$ by
$$
\frac {\partial^k G}{\partial Z_{i_1}\cdots \partial Z_{i_k}}:=
\frac{\partial} {\partial Z_{i_1}}\left(\frac{\partial} {\partial Z_{i_2}}\cdots \left( \frac{\partial G} {\partial Z_{i_k}}\right)\cdots \right).
$$
  These definitions can easily be  extended to formal power series.
  If $F:=\sum\limits_{\alpha\in \FF_n^+} A_{(\alpha)} \otimes Z_\alpha$ is a  power series with  operator-valued coefficients, then   we define the $k$-order  free partial derivative  of $F$ with respect to $Z_{i_1}, \ldots, Z_{i_k}$ to be the power series
$$
\frac {\partial^k F}{\partial Z_{i_1}\cdots \partial Z_{i_k}}
:=
\sum_{\alpha\in \FF_n^+} A_{(\alpha)} \otimes
\frac {\partial^k Z_\alpha}{\partial Z_{i_1}\cdots \partial Z_{i_k}}.
$$

\begin{proposition}\label{deriv-comu}
If $i,j\in \{1,\ldots, n\}$, then
$$\frac{\partial^2 F}{\partial Z_i \partial Z_j}=
\frac{\partial^2 F}{\partial Z_j \partial Z_i}
$$
for any formal power series $F$.
\end{proposition}
\begin{proof}
Due to linearity, it is enough to prove the result  for monomials.
Let $\alpha:=g_{i_1}\cdots g_{i_k}$ be a word in $ \FF_n^+$ and $Z_\alpha:=Z_{i_1}\cdots Z_{i_k}$. Let $i,j\in \{1,\ldots, n\}$ be such that  $i\neq j$. Assume that $Z_i$ occurs $q$ times in $Z_\alpha$, and $Z_j$ occurs $p$ times in $Z_\alpha$. Then
$\frac{\partial Z_\alpha}{\partial Z_i}$ is the sum of the $q$ words obtained by deleting each occurence of $Z_i$ in $Z_\alpha$. Notice that $Z_j$ occurs $p$ times in each of these $q$ words.
Therefore, $\frac{\partial ^2 Z_\alpha}{\partial Z_j \partial Z_i}$ is the sum of the $qp$ words  obtained by deleting each
occurence of $Z_i$ in $Z_\alpha$ and then deleting each  occurence of $Z_j$ in the resulting words.
Similarly, $\frac{\partial ^2 Z_\alpha}{\partial Z_i \partial Z_j}$ is the sum of the $qp$ words  obtained by deleting each
occurence of $Z_j$ in $Z_\alpha$ and then deleting each  occurence of $Z_i$ in the resulting words. Hence, it is clear that
$$\frac{\partial^2 Z_\alpha}{\partial Z_i \partial Z_j}=
\frac{\partial^2 Z_\alpha}{\partial Z_j \partial Z_i}.
$$
This completes the proof.
\end{proof}

\begin{theorem} \label{derivation}
Let $F=\sum\limits_{\alpha\in \FF_n^+} A_{(\alpha)} \otimes Z_\alpha$
be  a power series with radius of convergence $R$  and let $R'$ be the radius  of convergence of the power series \, $  \frac{\partial^k F} {\partial Z_{j_1}\cdots \partial Z_{j_k}}$, where $j_1,\ldots, j_k\in \{1,\ldots, n\}$. Then $R'\geq R$ and, in general, the inequality is strict.
\end{theorem}
\begin{proof}
It is enough to prove the result for first order  free partial derivatives.
For any word $\omega:=g_{i_1}\cdots g_{i_k}$, $|\omega|=k\geq 1$, and $0\leq m\leq k$, we define the insertion mapping of
$g_j$, $j=1,\ldots, n$, on the $m$ position  of $\omega$
by setting
$$
\chi(g_j, m,\omega):=
\begin{cases}
g_j\omega & \text{ if } m=0,\\
g_{i_1}\cdots g_{i_m} g_j g_{i_{m+1}}\cdots g_{i_k} & \text{ if } 1\leq m\leq k-1,\\
\omega g_j & \text{ if } m=|\omega|=k,
\end{cases}
$$
and $\chi(g_j, 0, g_0):=g_j$.
Let
$$
 \frac{\partial F} {\partial Z_j}=
\sum_{\beta\in \FF_n^+} B_{(\beta)} \otimes Z_\beta.
$$
Using the definition of the Hausdorff derivation   and the insertion mapping, we deduce that
$$
B_{(\beta)}=\sum_{m=0}^k A_{(\chi(g_j,m,\beta))}
$$
 for any $\beta\in \FF_n^+$ with $|\beta|=k$. This is the case, since the monomial $Z_\beta$ comes from free differentiation with respect to $Z_j$ of  the monomials $Z_{\chi(g_j,m,\beta)}$,\ $m=0,1,\ldots, |\beta|$.
Therefore, we have
\begin{equation*}
\begin{split}
\sum_{|\beta|=k} B_{(\beta)}^* B_{(\beta)}
&=
\sum_{|\beta|=k} \left( \sum_{m=0}^k A_{(\chi(g_j,m,\beta))}^*\right)\left( \sum_{m=0}^k A_{(\chi(g_j,m,\beta))}\right)\\
&\leq (k+1)\sum_{|\beta|=k} \sum_{m=0}^k A_{(\chi(g_j,m,\beta))}^*A_{(\chi(g_j,m,\beta))}
\\
&\leq (k+1)^2 \sum_{|\alpha|=k+1} A_{(\alpha)}^* A_{(\alpha)}.
\end{split}
\end{equation*}
The last inequality holds since, for each $j=1,\ldots,n$,
each $\alpha\in \FF_n^+$ with $|\alpha|=k+1$, and each $\beta\in \FF_n^+$ with $|\beta|=k$, the cardinal of the set
$$
\{(g_j,m,\beta):\ \chi(g_j,m,\beta)=\alpha, \text{ where } m=0,1,\ldots, k\}
$$
is $\leq k+1$.
Hence, we deduce that
$$
\left(\sum_{|\beta|=k} B_{(\beta)}^* B_{(\beta)} \right)^{1/2k}\leq (k+1)^{1/k} \left( \sum_{|\alpha|=k+1} A_{(\alpha)}^* A_{(\alpha)}\right)^{1/2k}.
$$
Consequently,
due to Theorem \ref{Abel}, we have  $\frac {1}{R'}\leq \frac {1}{R}$. Therefore, $R'\geq R$.

To prove the last part of the theorem, let $R_1, R_2>0$ be  such that $R_1<R_2$. Let us consider two power series
$$
F=\sum_{k=0}^\infty a_k Z_1^k \ \text { and } \
G=\sum_{k=0}^\infty b_k Z_2^k
$$
 with radius of convergence
$R_1$ and $R_2$, respectively.
We shall show that the  power series
$$
F+G=\sum_{k=0}^\infty (a_kZ_1^k+b_k Z_2^k)
$$
has  the radius of convergence equal to $R_1$.
First, since
$$
\sup_k\left( |a_k|^2+|b_k|^2\right)^{1/2k}\geq \sup|a_k|^{1/k}=\frac{1}{R_1},
$$
we deduce that the radius of convergence of $F+G$ is $\leq R_1$. On the other hand,
 if $r<R_1$, Corollary \ref{Cs} shows that both sequences $\{r^k |a_k|\}_{k=0}^\infty$ and $\{r^k |b_k|\}_{k=0}^\infty$
are bounded.  This implies that the sequence
$\{r^k \left(|a_k|^2+|b_k|^2\right)^{1/2}\}_{k=0}^\infty$
is bounded. Applying again Corollary \ref{Cs}, we can conclude that $F+G$ has radius of convergence $R_1$.
Since
$$
\frac{\partial (F+G)}{\partial Z_2}=\sum_{k=1}^\infty k b_kZ_2^{k-1},
$$
the power series
$ \frac{\partial (F+G)}{\partial Z_2}$ has radius of convergence $R_2$, which is strictly larger than the radius of convergence of $F+G$.
This completes the proof.
 \end{proof}

\bigskip

\section{Cauchy, Liouville, and Schwartz type results for free holomorphic functions} \label{Liouville}

In this section , we obtain Cauchy type estimates for the coefficients of free holomorphic functions  and a Liouville type theorem for free entire functions. Based on a noncommutative version of Gleason's problem
 \cite{R2}  and the noncommutative von Neumann inequality \cite{Po-von}, we provide a free analogue of Schwartz lemma.

First, we  obtain  Cauchy type estimates for the coefficients of  free holomorphic functions on the open ball $[B(\cH)^n]_{\gamma}$ with coefficients in $B(\cK)$.

\begin{theorem} \label{Cauchy-est} Let
$F:[B(\cH)^n]_{\gamma}\to B(\cK)\bar\otimes B(\cH)$   be a free holomorphic function  on $[B(\cH)^n]_{\gamma}$ with  the representation
 $$
F(X_1,\ldots, X_n)=\sum\limits_{k=0}^\infty
\left( \sum\limits_{|\alpha|=k} A_{(\alpha)}\otimes X_\alpha\right),
$$
and
  define
$$
M(\rho):=  \|F(\rho S_1,\ldots, \rho S_n)\|\quad
\text{for  any } \ \rho\in (0,\gamma),
$$
where $S_1,\ldots, S_n$ are the left creation operators on the full Fock space.
 Then, for each $k=0,1,\ldots,$
$$
\left\|\sum_{|\alpha|=k} A_\alpha^* A_\alpha\right\|^{1/2}\leq \frac{1} {\rho^k} M(\rho).
$$
\end{theorem}

\begin{proof}
Let $\{Y_{(\alpha)}\}_{|\alpha|=k}$ be an arbitrary sequence of operators in $B(\cK)$. Using Theorem \ref{caract-shifts}, we have
\begin{equation*}
\begin{split}
\left|\left<\left(\sum_{|\alpha|=k} Y_{(\alpha)}^*\otimes S_\alpha^*\right)F(\rho S_1,\ldots, \rho S_n) h\otimes 1, h\otimes 1\right>\right|
&\leq  \left\|\sum_{|\alpha|=k} Y_{(\alpha)}^*\otimes S_\alpha^*\right\| M(\rho) \|h\|^2\\
&=\left\|\sum_{|\alpha|=k} Y_{(\alpha)}^* Y_{(\alpha)}\right\|^{1/2} M(\rho) \|h\|^2
\end{split}
\end{equation*}
for any $h\in \cK$. On the other hand, since $S_1,\ldots, S_n$ are isometries with orthogonal ranges, we have
\begin{equation*}
\begin{split}
\Bigl|\Bigl<\Bigl(\sum_{|\alpha|=k} Y_{(\alpha)}^*\otimes S_\alpha^*\Bigr)F(\rho S_1,\ldots, &\rho S_n) h\otimes 1, h\otimes 1\Bigr>\Bigr|\\
&=
 \rho^k\left|\left<\left(\sum_{|\alpha|=k} Y_{(\alpha)}^* A_{(\alpha)}\otimes I\right) h\otimes 1, h\otimes 1\right>\right|\\
&=
\rho^k\left|\left<[Y_{(\alpha)}^*: |\alpha|=k]
\left[\begin{matrix}A_{(\alpha)}\\:\\|\alpha|=k\end{matrix} \right]h,h\right>\right|.
\end{split}
\end{equation*}
Combining these relations and taking $Y_{(\alpha)}:=A_{(\alpha)}$, $|\alpha|=k$, we deduce that
$$
\rho^k\left\|\left[\begin{matrix}A_{(\alpha)}\\:\\|\alpha|=k\end{matrix} \right]h\right\|^2\leq \left\|\left[\begin{matrix}A_{(\alpha)}\\:\\|\alpha|=k\end{matrix} \right]\right\| M(\rho) \|h\|^2
$$
for any $h\in \cK$.
Therefore,
$$\left\|\sum_{|\alpha|=k} A_\alpha^* A_\alpha\right\|^{1/2}=\|[A_{(\alpha)}^*:\ |\alpha|=k]\|\leq \frac{1} {\rho^k} M(\rho),
$$
which completes  the proof.
\end{proof}

A free holomorphic function with radius of convergence $R=\infty$ is called  free entire function. We can prove now the following noncommutative multivariable generalization of Liouville's theorem.

 \begin{theorem}\label{Liou}
Let $F$   be an entire function and let
$$
F(X_1,\ldots, X_n)=\sum_{k=0}^\infty\sum_{|\alpha|=k} A_{(\alpha)}\otimes X_\alpha
$$
 be its representation   on an infinite dimensional Hilbert space $\cH$.   Then  $F$ is a polynomial of degree $\leq m$, \, $m=0, 1, \ldots$,
if and only if  there are  constants $M>0$
 and  $C>1$ such that
$$
\|F(X_1,\ldots, X_n)\|\leq M\|[X_1,\ldots, X_n]\|^m
$$
for any  $[X_1,\ldots, X_n]\in B(\cH)^n$ such that
 $\|[X_1,\ldots, X_n]\|\geq C$.
 \end{theorem}

\begin{proof}
If  $F=\sum_{|\alpha|\leq m} A_{(\alpha)}\otimes X_\alpha$
is a polynomial, then
\begin{equation*}
\begin{split}
\|F\|&\leq \sum_{k=0}^m \left\|\sum_{|\alpha|=m} A_{(\alpha)}\otimes X_\alpha \right\|\\
&\leq \sum_{k=0}^m \left\| \sum_{|\alpha|=k} A_{(\alpha)}^* A_{\alpha)} \right\|^{1/2}\|[X_1,\ldots, X_n]\|^k
\end{split}
\end{equation*}
if $\|[X_1,\ldots, X_n]\|\geq 1$. Therefore, there exists $M>0$
and $R>1$ such that
\begin{equation}\label{f-norm}
\|F(X_1,\ldots, X_n)\|\leq M \|[X_1,\ldots, X_n]\|^k
\end{equation}
for  any  $n$-tuple of operators $[X_1,\ldots, X_n]$ with
$\|[X_1,\ldots, X_n]\|\geq R$.

Conversely, if the   inequality \eqref{f-norm}
 holds, then
 $$
\|F(\rho S_1,\ldots, \rho S_n)\|\leq M \rho^m,\quad \text{ as }\ \rho\to\infty.
$$
 According to Theorem \ref{Cauchy-est}, we have
$$
\left\|\sum_{|\alpha|=k} A_\alpha^* A_\alpha\right\|^{1/2}\leq \frac{1} {\rho^k} M(\rho),
$$
where $M(\rho):=\|F(\rho S_1,\ldots, \rho S_n)\|$.
Combining these inequalities, we deduce that
$$
\left\|\sum_{|\alpha|=k} A_\alpha^* A_\alpha\right\|^{1/2}\leq M\frac{1} {\rho^{k-m}}.
$$
Consequently, if $k>m$ and $\rho\to\infty$, we obtain
$\sum_{|\alpha|=k} A_\alpha^* A_\alpha=0$.
This shows that $A_{(\alpha)}=0$ for any $\alpha\in \FF_n^+$ with $|\alpha|>m$.
 \end{proof}

We say that a free holomorphic function $F$ on the open  operatorial $n$-ball of radius $\gamma$ is bounded if
$$
\|F\|_\infty:=\sup  \|F(X_1,\ldots, X_n)\|<\infty,
$$
where the supremum is taken over all $n$-tuples  of operators
$[X_1,\ldots, X_n]\in (B(\cH)^n)_\gamma$ and any Hilbert space $\cH$.
In the particular case when $m=0$, Theorem \ref{Liou} implies   the following free analogue
   of Liouville's theorem from  complex analysis (see \cite{R}, \cite{Co}).

\begin{corollary}  If  $F$ is a bounded    free entire function,    then  it  is    constant.
 \end{corollary}

We recall   that the joint spectral radius of an $n$-tuple of operators $[T_1,\ldots, T_n]\in B(\cH)^n$,
$$r(T_1,\ldots,T_n):=\lim_{k\to\infty}
\left\|\sum_{|\alpha|=k} T_\alpha T_\alpha^*\right\|^{\frac{1} {2k}},
$$
  is also equal to the spectral radius of the reconstruction operator
 $S_1\otimes T_1^*+\cdots + S_n\otimes T_n^*$ (see \cite{Po-unitary}).
Consequently, $r(T_1,\ldots,T_n)<1$ if and only if
$$
\sigma(S_1\otimes T_1^*+\cdots + S_n\otimes T_n^*)\subset \DD.
$$
Moreover, the joint right spectrum
~$\sigma_r(T_1,\ldots, T_n)$~ is included in the closed ball of $\CC^n$ of radius equal to $r(T_1,\ldots,T_n)$.
We recall that
  $\sigma_r(T_1,\ldots, T_n)$  is the set of all $n$-tuples $(\lambda_1,\ldots, \lambda_n)\in \CC^n$ such that the right ideal of $B(\cH)$ generated  by $\lambda_1 I-T_1,\ldots, \lambda_n I-T_n$ does not contain the identity.

Now, we  prove an analogue of Schwartz lemma, in our multivariable operatorial setting.

\begin{theorem}\label{Schwartz}
Let $F(X_1,\ldots, X_n)=\sum_{\alpha\in \FF_n^+} A_{(\alpha)}\otimes X_\alpha$, \ $A_{(\alpha)}\in B(\cK)$,  be a free holomorphic function on $[B(\cH)^n]_1$  with the properties:
\begin{enumerate}
\item[(i)] $\|F\|_\infty\leq 1$ and
\item[(ii)] $A_{(\beta)}=0$ for any $\beta\in \FF_n^+$ with $|\beta|\leq m-1$, where $m=1,2,\ldots$.
\end{enumerate}
Then
$$
\|F(X_1,\ldots, X_n)\| \leq \|[X_1,\ldots, X_n]\|^m \quad
\text{
and }\quad
r(F(X_1,\ldots, X_n))\leq r(X_1,\ldots, X_n)^m
$$
for any \, $[X_1,\ldots, X_n]\in [B(\cH)^n]_1$.
Moreover,
$$
\left\|\sum_{|\alpha|=k} A_\alpha A_\alpha^*\right\|^{1/2} \leq 1\quad  \text{ for any } \  k\geq m.
$$
\end{theorem}

\begin{proof}
For each $\beta\in \FF_n^+$ with $|\beta|\leq m$, define the formal power series
$$
\Phi_{(\beta)}(Z_1,\ldots, Z_n):=\sum_{\alpha\in \FF_n^+} A_{(\beta\alpha)}\otimes Z_\alpha.
$$
Since
$$
\left\|\sum_{|\alpha|=k} A_{(\beta \alpha)}^* A_{(\beta\alpha)}^*\right\|\leq
\left\| \sum_{|\gamma|=m+k} A_{(\gamma)}^* A_{(\gamma)}
\right\|,
$$
we deduce that
$$
\limsup_{k\to\infty}\left\|\sum_{|\alpha|=k} A_{(\beta \alpha)}^* A_{(\beta\alpha)}^*\right\|^{1/2k}
\leq \limsup_{k\to\infty} \left\| \sum_{|\gamma|=m+k} A_{(\gamma)}^* A_{(\gamma)}
\right\|^{\frac{1}{2(m+k)}}.
$$
Consequently, due to Theorem \ref{Abel}, the radius of convergence of $\Phi_{(\beta)}$ is greater than  the radius of convergence of $F$. Therefore, $\Phi_{(\beta)}$ represents a free holomorphic function on the open  operatorial unit $n$-ball.
Since $A_{(\beta)}=0$ for any $\beta\in \FF_n^+$ with $|\beta|\leq m-1$,  and due to Theorem \ref{operations}, we have the following Gleason type decomposition
$$
F(Z_1,\ldots, Z_n)=\sum_{|\beta|=m}\left[(I_\cK\otimes Z_\beta)\sum_{\alpha\in \FF_n^+} A_{\beta \alpha)} \otimes Z_\alpha\right]= \sum_{|\beta|=m}(I_\cK\otimes Z_\beta)\Phi_{(\beta)}(Z_1,\ldots, Z_n).
$$
Therefore,
\begin{equation}\label{F-Phi}
F(rS_1,\ldots, rS_n)
=\sum_{|\beta|=m} (I_\cK\otimes r^{|\beta|} S_\beta )\Phi_{(\beta)}(rS_1,\ldots, rS_n)
\end{equation}
for any $r\in [0,1)$.
  Since $S_1,\ldots, S_n$ are isometries with orthogonal ranges,  $S_\beta$, $|\beta|=m$, are  also isometries with orthogonal ranges  and we have
$$
F(rS_1,\ldots, rS_n)^*F(rS_1,\ldots, rS_n)
=r^{2m}\sum_{|\beta|=m} \Phi_{(\beta)}(rS_1,\ldots, rS_n)^*\Phi_{(\beta)}(rS_1,\ldots, rS_n).
$$
Now, due to  the noncommutative von Neumann inequality \cite{Po-von} and Theorem  \ref{caract-shifts}, we deduce that
\begin{equation}\label{M-P}
\left\|\left[\begin{matrix}
\Phi_{(\beta)}(rX_1,\ldots, rX_n)\\
:\\
|\beta|=m
\end{matrix}
\right]\right\|\leq \left\|\left[\begin{matrix}
\Phi_{(\beta)}(rS_1,\ldots, rS_n)\\
:\\
|\beta|=m
\end{matrix}
\right]\right\|
\end{equation}
for any $[X_1,\ldots, X_n]\in [B(\cH)^n]_1$.
Consequently, using relations \eqref{F-Phi} and \eqref{M-P}, we obtain
\begin{equation*}
\begin{split}
\|F(rX_1,\ldots, rX_n)\|
&=
\left\|\sum_{|\beta|=m} (I_\cK\otimes r^{|\beta|} X_\beta )\Phi_{(\beta)}(rX_1,\ldots, rX_n)\right\| \\
&\leq
\left\|[r^m X_\beta:\ |\beta|=m]\right\|
\left\|\left[\begin{matrix}
\Phi_{(\beta)}(rX_1,\ldots, rX_n)\\
:\\
|\beta|=m
\end{matrix}
\right]\right\| \\
&\leq
r^m\left\|\sum_{|\beta|=m} X_\beta X_\beta^*\right\|^{1/2}\left\|\left[\begin{matrix}
\Phi_{(\beta)}(rS_1,\ldots, rS_n)\\
:\\
|\beta|=m
\end{matrix}
\right]\right\| \\
&=
r^m\left\|\sum_{|\beta|=m} X_\beta X_\beta^*\right\|^{1/2}
\left\|
\sum_{|\beta|=m} \Phi_{(\beta)}(rS_1,\ldots, rS_n)^*\Phi_{(\beta)}(rS_1,\ldots, rS_n)
\right\|^{1/2}\\
&=
r^m\left\|\sum_{|\beta|=m} X_\beta X_\beta^*\right\|^{1/2}
\left\|F(rS_1,\ldots, rS_n)^*F(rS_1,\ldots, rS_n)\right\|\\
&\leq
r^m\left\|\sum_{|\beta|=m} X_\beta X_\beta^*\right\|^{1/2}
\|F\|_\infty\\
&\leq r^m\left\|\sum_{|\beta|=m} X_\beta X_\beta^*\right\|^{1/2}
\leq r^m\|[X_1,\ldots, X_n]\|^m.
 \end{split}
\end{equation*}
Taking $r\to 1$ and using the continuity of the free holomorphic function  $F$  on $[B(\cH)^n]_1$ (see Theorem \ref{continuous}), we infer that
$$
\|F(X_1,\ldots, X_n)\|\leq \left\|\sum_{|\beta|=m} X_\beta X_\beta^*\right\|^{1/2}
\leq \|[X_1,\ldots, X_n]\|^m
$$
for any $[X_1,\ldots, X_n]\in [B(\cH)^n]_1$.

Due to Theorem \ref{operations}, the power series  $F^k =\sum_{\alpha\in \FF_n^+} B_{(\alpha)}\otimes Z_\alpha$  represents a free holomorphic function on the open operatorial unit $n$-ball, with
$B_{(\alpha)}=0$ for any $\alpha\in \FF_n^+$ with $|\alpha|\leq mk$. Applying the above inequality to $F^k$, we obtain
$$
\|F(X_1,\ldots,  X_n)^k\|\leq \left\|\sum_{|\beta|=mk} X_\beta X_\beta^*\right\|^{1/2}\leq \left\|\sum_{|\beta|=k} X_\beta X_\beta^*\right\|^{m/2}.
$$
Hence, and using the definition of the joint spectral radius,  we deduce that $r(F(X_1,\ldots, X_n))\leq r(X_1,\ldots, X_n)^m.
$

To prove the last part of the theorem, notice that, according to Theorem
\ref{Cauchy-est}, we have
$$
\left\|\sum_{|\alpha|=k} A_{(\alpha)} A_{(\alpha)}^*\right\|^{1/2}\leq \frac{1}{\rho^k} M(\rho)
$$
for any $\rho\in (0,1)$, where
$M(\rho)=\|F(\rho S_1,\ldots, \rho S_n)\|$.
Since $M(\rho)\leq \|F\|_\infty\leq 1$, we take $\rho\to 1$ and deduce that
$\left\|\sum_{|\alpha|=k} A_{(\alpha)} A_{(\alpha)}^*\right\|^{1/2}\leq 1
$
for any $k\geq m$. The proof is complete.
 \end{proof}

In the scalar case we get a little bit more.

\begin{corollary}
Let $f(X_1,\ldots, X_n)=\sum_{\alpha\in \FF_n^+} a_\alpha X_\alpha$, $a_\alpha\in \CC$,   be a free holomorphic function  on  $[B(\cH)^n]_1$ with scalar coefficients and the properties:
\begin{enumerate}
\item[(i)] $\|f\|_\infty\leq 1$ and
\item[(ii)] $f(0)=0$.
\end{enumerate}
Then
\begin{enumerate}
\item[(iii)]
$\|f(X_1,\ldots, X_n)\|\leq \left\|[X_1,\ldots, X_n]\right\| $ \
and  \ $r(f(X_1,\ldots, X_n))\leq r(X_1,\ldots, X_n)
$
\quad
for any  $n$-tuple $ \ [X_1,\ldots, X_n]\in [B(\cH)^n]_1$;
\item[(iv)]
$
\sum\limits_{i=1}^n \left|\frac{\partial f}{\partial X_i}(0)\right|^2\leq 1.
$
\end{enumerate}
Moreover, if $
\sum\limits_{i=1}^n \left|\frac{\partial f}{\partial X_i}(0)\right|^2= 1,
$
then $\|f\|_\infty=1$.
\end{corollary}

\begin{proof} The first part of this corollary is a particular case of Theorem \ref{Schwartz}, when $m=1$ and $\cK=\CC$. To prove the second part, assume that
 $
\sum\limits_{i=1}^n \left|\frac{\partial f}{\partial X_i}(0)\right|^2= 1
$
Consequently, we have $\sum_{i=1}^n |a_i|^2=1$.
Hence, and due to Theorem \ref{Cauchy-est}, we have
$$
1\leq \sum_{i=1}^n |a_i|^2\leq \frac{1}{\rho} \|f\|_\infty $$
for any $0<\rho<1$.
Therefore, $\|f\|_\infty=1$. This completes the proof.
 \end{proof}

\bigskip
\smallskip

\section{ Algebras of free holomorphic functions }
\label{algebras}

In this section,   we introduce two  Banach algebras of free holomorphic functions,
$H^\infty(B(\cX)^n_1)$ and $A(B(\cX)^n_1)$, and prove that they are  isometrically isomorphic to the the noncommutative analytic Toeplitz algebra $F_n^\infty$ and the noncommutative disc algebra $\cA_n$, respectively.
  The results of this section are used to obtain a maximum principle for free holomorphic functions.

We denote by $Hol(B(\cX)^n_\gamma)$ the set of all free holomorphic functions with scalar coefficients on the open  operatorial  $n$-ball of radius $\gamma$.
Due to Theorem \ref{operations} and Theorem \ref{Abel}, $Hol(B(\cX)^n_\gamma)$ is an algebra and
 an element
$F=\sum\limits_{\alpha\in \FF_n^+} a_\alpha Z_\alpha$ is in
$Hol(B(\cX)^n_\gamma)$ if and only if
$$\limsup_{k\to\infty}\left(\sum_{|\alpha|=k}
|a_\alpha|^2\right)^{1/2k}\leq 1.
$$
 Let $H^\infty(B(\cX)^n_1)$  denote the set of  all elements $F$ in
$Hol(B(\cX)^n_1)$     such that
$$\|F\|_\infty:= \sup \|F(X_1,\ldots, X_n)\|<\infty,
$$
where the supremum is taken over all $n$-tuples $[X_1,\ldots, X_n]\in [B(\cH)^n]_1$ and any Hilbert space $\cH$.
 We  denote by  $A(B(\cX)^n_1)$  be the set of all  elements $F$  in $Hol(B(\cX)^n_1)$   such that,  for any Hilbert space $\cH$, the
mapping
$$[B(\cH)^n]_1\ni (X_1,\ldots, X_n)\mapsto F(X_1,\ldots, X_n)\in B(\cH)$$
 has a continuous extension to the closed unit ball $[B(\cH)^n]^-_1$.
In this section, we will show that
$H^\infty(B(\cX)^n_1)$  and
$A(B(\cX)^n_1)$ are Banach algebras under pointwise multiplication and the norm $\|\cdot \|_\infty$, which can be identified with the noncommutative analytic Toeplitz algebra $F_n^\infty$ and the noncommutative disc algebra $\cA_n$, respectively.

Let us recall (see \cite{Po-von}, \cite{Po-funct}, \cite{Po-disc},
\cite{Po-poisson}) a few facts about the  Banach algebras $\cA_n$ and $F_n^\infty$.
Any element $~f$ in the full Fock space $ F^2(H_n)$
has the form
$$
f=\sum\limits_{\alpha\in \FF_n^+} a_\alpha e_\alpha,\quad\text{\ with\ }a_\alpha\in\CC, \quad\text
{\ such that\ }\
\|f\|_2:=\left(\sum\limits_{\alpha\in\FF_n^+}|a_\alpha|^2\right)^{1/2}<\infty.
$$
 If
$g=\sum\limits_{\alpha\in \FF_n^+} b_\alpha e_\alpha\in F^2(H_n)$,   we define the product $f\otimes g$ to be the formal
power series
$$
f\otimes g:=\sum_{\gamma\in \FF_n^+} c_\gamma e_\gamma,\quad
\text{ where }\quad  c_\gamma :=\sum_{\stackrel{\alpha,\beta\in \FF_n^+}{ \alpha\beta=\gamma}} a_\alpha b_\beta, \quad  \gamma\in \FF_n^+.
$$
We also make the natural identification of $e_\alpha\otimes 1$ and
$1\otimes e_\alpha$ with $e_\alpha$.
Let $\cP$  denote the set of all  polynomials $p\in F^2(H_n)$, i.e., elements of the form $p=\sum_{|\alpha|\leq m} a_\alpha e_\alpha$,
where $m=0,1,\ldots$.

 In \cite{Po-von}, we introduced
the noncommutative Hardy algebra $F_n^\infty$ as the set of all $f\in F^2(H_n)$ such that
\begin{equation}\label{norm}
\|f\|_\infty:=\sup\{\|f\otimes p\|_2:p\in\cP, \ \|p\|_2\le 1\}<\infty.
\end{equation}
If $f\in F^2(H_n)$, then  $f\in F_n^\infty$ if and only if
   $f\otimes g\in F^2(H_n)$ for any $g\in F^2(H_n)$.
Moreover, if $f\in F_n^\infty$, then the left multiplication mapping
$L_f:F^2(H_n)\to F^2(H_n)$ defined by
$$L_fg:=f\otimes g, \quad g\in F^2(H_n),$$
is a bounded linear operator with
$\|L_f\|=\|f\|_\infty$.
The noncommutative Hardy algebra $F_n^\infty$  is isometrically isomorphic  to the left
multiplier algebra of the full Fock space $F^2(H_n)$, which is also called the noncommutative Toeplitz algebra.
Under this identification,  $F_n^\infty$ is the weakly closed algebra generated by the left creation operators $S_1,\ldots, S_n$ and the identity. The noncommutative disc algebra $\cA_n$ was introduced in \cite{Po-von} as is the norm closed algebra generated by the left creation operators   and the identity.

  Let
$
f=\sum\limits_{\alpha\in \FF_n^+} a_\alpha e_\alpha$  be an element in  $F^2(H_n)$ and
  define
$$
f_r:=\sum\limits_{\alpha\in \FF_n^+} r^{|\alpha|} a_\alpha e_\alpha
\quad \text{
for  } \ ~0<r<1.
$$
 In \cite{Po-funct}, \cite{Po-poisson},  we proved that if $~f\in F_n^\infty~$ then   $\|f_r\|_\infty\leq \|f\|_\infty$  for  $ 0\leq r<1$, and
\begin{equation}\label{So}
  L_f=\text{\rm{SOT-}}\lim\limits_{r\to1}f_r(S_1,\dots,S_n),
\end{equation}
where $f_r(S_1,\ldots, S_n):=\sum_{k=0}^\infty\sum_{|\alpha|=k} r^{|\alpha|} a_\alpha S_\alpha$.
 Moreover, if  $f\in \cA_n$ then the above limit exists in
  the operator norm topology.

  We identify $M_m(B(\cH))$, the set of
$m\times m$ matrices with entries from $B(\cH)$, with
$B( \cH^{(m)})$, where $\cH^{(m)}$ is the direct sum of $m$ copies
of $\cH$.
Thus we have a natural $C^*$-norm on
$M_m(B(\cH))$. If $X$ is an operator space, i.e., a closed subspace
of $B(\cH)$, we consider $M_m(X)$ as a subspace of $M_m(B(\cH))$
with the induced norm.
Let $X, Y$ be operator spaces and $u:X\to Y$ be a linear map. Define
the map
$u_m:M_m(X)\to M_m(Y)$ by
$$
u_m ([x_{ij}]):=[u(x_{ij})].
$$
We say that $u$ is completely bounded ($cb$ in short) if
$$
\|u\|_{cb}:=\sup_{m\ge1}\|u_m\|<\infty.
$$
If $\|u\|_{cb}\leq1$
(resp. $u_m$ is an isometry for any $m\geq1$) then $u$ is completely
contractive (resp. isometric),
 and if $u_m$ is positive for all $m$, then $u$ is called
 completely positive.

For each $m=1,2,\ldots$, we define the norms $\|\cdot
\|_m:M_m\left(H^\infty(B(\cX)^n_1)\right)\to [0,\infty)$ by setting
$$
\|[F_{ij}]_m\|_m:= \sup \|[F_{ij}(X_1,\ldots, X_n)]_m\|,
$$
where the supremum is taken over all $n$-tuples $[X_1,\ldots,
X_n]\in [B(\cH)^n]_1$ and any Hilbert space $\cH$. It is easy to see
that the norms  $\|\cdot\|_m$, $m=1,2,\ldots$, determine  an
operator space structure  on $H^\infty(B(\cX)^n_1)$,
 in the sense of Ruan (see \cite{ER}).

\begin{theorem}\label{f-infty}
Let $F:=\sum_{\alpha\in \FF_n^+} a_\alpha Z_\alpha$ be a free holomorphic function on the open  operatorial unit $n$-ball. Then the following statements are equivalent:
\begin{enumerate}
\item[(i)] $F$ is in $H^\infty(B(\cX)^n_1)$;
\item[(ii)] $f:=\sum_{\alpha\in \FF_n^+} a_\alpha e_\alpha$ is in $F_n^\infty$;
\item[(iii)]$\sup\limits_{0\leq r<1}\|F(rS_1,\ldots, rS_n)\|<\infty$;
\item[(iv)]
The map $\varphi:[0,1)\to B(F^2(H_n))$
defined by
$$
\varphi(r):=F(rS_1,\ldots, rS_n)
\quad
\text{
for any } \ r\in [0,1)
$$
has a continuous extension to $[0,1]$ with respect to the strong
operator topology of $B(F^2(H_n))$.
\end{enumerate}
In this case,  we have
\begin{equation}\label{many eq}
\|L_f\|=\|f\|_\infty=\sup_{0\leq r<1}\|F(rS_1,\ldots, rS_n)\|=
\lim_{r\to 1}\|F(rS_1,\ldots, rS_n)\|=\|F\|_\infty.
\end{equation}
 Moreover, the map
$$
\Phi:H^\infty((B(\cX)^n_1)\to F_n^\infty\quad \text{ defined by } \quad \Phi(F):=f
$$ is a completely isometric isomorphism of operator algebras.
\end{theorem}

\begin{proof} Assume (ii) holds.
 Since $f\in F_n^\infty$, we have
\begin{equation}\label{f-inf}
\|F(rS_1,\ldots, rS_n)\|=\|f(rS_1,\ldots, rS_n)\|=\|L_{f_r}\|=\|f_r\|\leq  \|f\|_\infty
\end{equation}
for any $r\in [0,1)$. Therefore, (ii)$\implies$(iii).
To prove that  (iii)$\implies$(ii),
assume that  (iii) holds. Consequently, we have
\begin{equation*}
\begin{split}
\sum_{\alpha\in \FF_n^+} r^{2|\alpha|} |a_\alpha|^2&=
\left\|\sum_{\alpha\in \FF_n^+} r^{|\alpha|} a_\alpha S_\alpha(1)\right\|\\
&\leq \sup\limits_{0\leq r<1}\|F(rS_1,\ldots, rS_n)\|<\infty
\end{split}
\end{equation*}
for any $0\leq r<1$. Hence, $\sum_{\alpha\in \FF_n^+} |a_\alpha|^2<\infty$, which shows that $f:=\sum_{\alpha\in \FF_n^+} a_\alpha e_\alpha$ is in $F^2(H_n)$.
Now assume that
  $f\notin F_n^\infty$.
 Due to the definition of $F_n^\infty$, given an arbitrary positive number $M$, there exists a polynomial $q\in \cP$ with $\|q\|_2=1$ such that
$$
\|f\otimes q\|_2>M.
$$
Since $\|f_r-f\|_2\to 0$ as $r\to 1$, we have
$$\|f\otimes q-f_r\otimes q\|_2=\|(f-f_r)\otimes q\|_2\to 0,
\quad
\text{ as }\ r\to 1.
$$
Therefore, there is $r_0\in (0,1)$ such that
$
\|f_{r_0}\otimes q\|_2> M.
$
Hence,
$$
\|f_{r_0}(S_1,\ldots, S_n)\|=\|L_{f_{r_0}}\|=\|f_{r_0}\|_\infty>M.
$$
Since $M>0$ is arbitrary, we deduce that
$$
\sup_{0\leq r<1}\|f(rS_1,\ldots, rS_n)\|=\infty,
$$
which is a contradiction.
Consequently, (ii)$\Longleftrightarrow$(iii).
Now, let us prove that (ii)$\implies$(iv).
Assume (ii) and define the map $\tilde\varphi:[0,1]\to B(F^2(H_n)$ by setting
$$
\tilde\varphi(r):= \begin{cases}
F(rS_1,\ldots, rS_n) &\quad  \text{if } 0\leq r<1\\
L_f &\quad  \text{if } r=1.
\end{cases}
$$
Since $f(rS_1,\ldots, rS_n)=F(rS_1,\ldots, rS_n)$, \ $0\leq r<1$,  the SOT-continuity of $\tilde\varphi$ at $r=1$ is due to relation \eqref{So}, while the continuity  of $\tilde\varphi$ on $[0,1)$ is a consequence of Theorem \ref{continuous}.
Therefore, the item (iv) holds.

Assume  now that (iv) holds.
For each $x\in F^2(H_n)$, the map $[0,1)\ni \mapsto \|\varphi(r)x\|\in \RR^+$ is  bounded, i.e.,
$\sup\limits_{0\leq r<1}\|\varphi(r)x\|<\infty$.  Due  to the principle of uniform  boundedness, we deduce condition (iii).

The implication (i)$\implies$(iii) is obvious, and
the implication (iii)$\implies$(i) is due to Theorem
\ref{Abel} and the noncommutative von Neumann inequality.
Indeed,
if $[X_1,\ldots, X_n]\in [B(\cH)^n]_1$, $\cH$ is an arbitrary Hilbert space,  and
$\|[X_1,\ldots, X_n]\|=r<1$, then
$$
\left\|\sum_{k=0}^m \sum_{|\alpha|=k} a_\alpha X_\alpha \right\|\leq \left\|\sum_{k=0}^m \sum_{|\alpha|=k} r^{|\alpha|}a_\alpha S_\alpha \right\|, \quad m=1,2,\ldots.
$$
Hence, and taking into account Theorem \ref{Abel}, we deduce
that
$$\|F(X_1,\ldots, X_n)\|\leq \|F(rS_1, \ldots, rS_n)\|,\quad \text{ for any } \ r\in [0,1).
$$
Consequently,
\begin{equation}\label{supsup}
\sup_{[X_1,\ldots, X_n]\in [B(\cH)^n]_1}
\|F(X_1,\ldots, X_n)\|\leq \sup_{0\leq r<1} \|F(rS_1, \ldots, rS_n)\|<\infty,
\end{equation}
whence (i) holds.

We  prove  now the last part of the theorem. If $f\in F_n^\infty$ and  $\epsilon>0$, then there exists a polynomial $q\in \cP$ with $\|q\|_2=1$ such that
$$
\|f\otimes q\|_2>\|f\|_\infty-\epsilon.
$$
Due to relation \eqref{So}, there exists $r_0\in (0,1)$ such that
$\|f_{r_0}(S_1,\ldots, S_n)q\|>\|f\|_\infty-\epsilon$.
Using now relation \eqref{f-inf}, we deduce that
$$
\sup_{0\leq r<1}\|f(rS_1,\ldots, rS_n)\|=\|f\|_\infty.
$$

Now,  let $r_1,r_2\in [0,1)$ with $r_1<r_2$ and let
$f:=\sum_{\alpha\in \FF_n^+}a_\alpha e_\alpha $.
Since
$g:=\sum_{\alpha\in \FF_n^+} r_2^{|\alpha|}a_\alpha  e_\alpha $ is in the noncommutative disc algebra
$\cA_n$, we
have  $\|g_r\|_\infty\leq \|g\|_\infty$ for any $0\leq r<1$.
In particular, when    $r:=\frac {r_1}{r_2}$,   we deduce that
$$
\|f_{r_1}(S_1,\ldots, S_n)\|\leq \|f_{r_2}(S_1,\ldots, S_n)\|.
$$
Consequently,  the   function $[0,1]\ni r\to \|f(rS_1,\ldots,
rS_n)\|\in \RR^+$ is increasing. Hence, and using  relation
\eqref{supsup}, we deduce \eqref{many eq}. Using the same
techniques, one can prove  a  matrix form of
 relation \eqref{many eq}. In particular, we have
 $\|[F_{ij}]_m\|_m=\|[L_{f_{ij}}]_m\|$ for any  $[F_{ij}]_m\in
 M_m\left(H^\infty(B(\cX)^n_1)\right)$ and $m=1,2,\ldots$.
 Hence, we deduce that $\Phi$ is a complete isometry of
 $ H^\infty(B(\cX)^n_1)$ onto $F_n^\infty$.
The proof is complete.
\end{proof}

\begin{theorem}\label{A-infty} Let $F:=\sum_{\alpha\in \FF_n^+} a_\alpha Z_\alpha$ be a free holomorphic function on the open  operatorial unit $n$-ball. Then the following statements are equivalent:
\begin{enumerate}
\item[(i)] $F$ is in $A(B(\cX)^n_1)$;
\item[(ii)] $f:=\sum_{\alpha\in \FF_n^+} a_\alpha e_\alpha$ is in $\cA_n$;
\item[(iii)]
 The map $\varphi:[0,1)\to B(F^2(H_n))$ defined by
$$
\varphi(r):=F(rS_1,\ldots, rS_n)
$$
has a  continuous  extension to $[0,1]$,
 with respect to the operator norm topology of $B(F^2(H_n))$.
\end{enumerate}
Moreover,
 the map
$$\Psi:A((B(\cX)^n_1)\to \cA_n\quad \text{ defined by } \quad \Psi(F):=f
$$
 is a completely isometric isomorphism of operator algebras.
\end{theorem}

\begin{proof}
The implication (i)$\implies$(iii) is due to the definition of $A(B(\cX)^n_1)$. Assume that  item (ii) holds, i.e.,
 $f\in \cA_n$.
The norm continuity of $\varphi$ on [0,1) is due to Theorem
\ref{continuous}, while the continuity of $\varphi$  at $r=1$  is due to the fact that
$\lim_{r\to 1} f_r(S_1,\dots, S_n)=L_f$ in the operator norm for any $f\in \cA_n$ (see the remarks preceeding  this theorem). Therefore,
the implication (ii)$\implies$(iii) is true.
Conversely, assume item (iii) holds. Then
$\lim_{r\to \infty} F(rS_1,\ldots, rS_n) $ exists in the operator norm.
Since $F(rS_1,\ldots, rS_n)\in \cA_n$ and $\cA_n$ is a Banach algebra, there exists $g\in \cA_n$  such that $L_g=\lim_{r\to \infty} F(rS_1,\ldots, rS_n)$ in the operator norm.
On the other hand, due to Theorem \ref{f-infty}, we deduce that
$f:=\sum_{\alpha\in \FF_n} a_\alpha e_\alpha\in F_n^\infty$.
Since $f(rS_1,\ldots, rS_n)=F(rS_1,\ldots, rS_n)$, \ $0\leq r<1$, and
$L_f=\text{\rm SOT}-\lim_{r\to \infty} f(rS_1,\ldots, rS_n)$, we conclude that $L_f=L_g$, i.e., $f=g$. Therefore, condition (ii) holds.

It remains to prove that (ii)$\implies$(i).
According to \cite{Po-funct} (see also \cite{Po-poisson}), if  $f\in\cA_n$ then,
for any  $n$-tuple $[Y_1,\ldots, Y_n]\in [B(\cH)^n]_1^-$,
$$
\tilde F(Y_1,\ldots, Y_n):= \lim_{r\to 1} f(rY_1,\ldots, rY_n),
$$
exists in the operator norm, and
$$
\|\tilde F(Y_1,\ldots, Y_n)\|\leq \|f\|_\infty\quad \text{ for any } \ [Y_1,\ldots, Y_n]\in [B(\cH)^n]_1^-.
$$
 Notice also that $\tilde F$ is an extension of the free holomorphic function $F$ on $[B(\cH)^n]_1$.
Indeed, if $ [X_1,\ldots, X_n]\in [B(\cH)^n]_1$, then
\begin{equation*}
\begin{split}
\tilde F(X_1,\ldots, X_n)&=\lim_{r\to 1}f(rX_1,\ldots, rX_n)\\
&=\lim_{r\to 1}F(rX_1,\ldots, rX_n)=F(X_1,\ldots, X_n).
\end{split}
\end{equation*}
The last equality is due to Theorem \ref{continuous}.

Let us prove that $\tilde F:[B(\cH)^n]_1^-\to B(\cH)$ is continuous.
Since $f\in \cA_n$,  for any $\epsilon>0$ there exists $r_0\in [0,1)$ such that
$\|L_f-f(r_0S_1,\ldots, r_0 S_n)\|<\epsilon$.
Applying the above mentioned result from \cite{Po-poisson} to $ f-\ f_{r_0}\in \cA_n$, we deduce that
\begin{equation}
\label{tild-f}
\|\tilde F(T_1,\ldots, T_n)-f_{r_0}(T_1,\ldots, T_n)\|\leq \|L_f-L_{f_{r_0}} \|< \frac{\epsilon}{3}
\end{equation}
for any $[T_1,\ldots, T_n]\in [B(\cH)^n]_1^-$.
Due to Theorem \ref{continuous}, $F$ is a continuous function on $[B(\cH)^n]_1$. Therefore, there exists $\delta>0$ such that
$$
\|F_{r_0}(T_1,\ldots, T_n)-F_{r_0}(Y_1,\ldots, Y_n)\|<\frac{\epsilon}{3}
$$
 for any $n$-tuples $[T_1,\ldots, T_n]$ and  $[Y_1,\ldots, Y_n]$   in $[B(\cH)^n]_1^-$     such that  $\|[T_1-Y_1,\ldots, T_n-Y_n]\|<\delta$.
Hence, and using \eqref{tild-f}, we have
\begin{equation*}
\begin{split}
\|\tilde F(T_1,\ldots, T_n)-\tilde F(Y_1,\ldots, Y_n)\|
&\leq \|\tilde F(T_1,\ldots, T_n)-f_{r_0}(T_1,\ldots, T_n)\|\\
&\qquad + \| f_{r_0}(T_1,\ldots, T_n)- f_{r_0}(Y_1,\ldots, Y_n)\|\\
&\qquad + \|f_{r_0}(Y_1,\ldots, Y_n)-\tilde F(Y_1,\ldots, Y_n)\|
<\epsilon,
\end{split}
\end{equation*}
whenever $\|[T_1-Y_1,\ldots, T_n-Y_n]\|<\delta$.
This proves the continuity of $\tilde F$ on $[B(\cH)^n]_1^-$.
Therefore, $F\in A(B(\cX)^n_1)$.

To prove the last part of the theorem, notice that if
$f\in\cA_n\subset F_n^\infty$, then by Theorem \ref{f-infty} (see
relation \eqref{many eq} and its matrix form), we have
$\|[L_{f_{ij}}]_m\|=\|[F_{ij}]_m\|_m$. Since $\cA_n\subset
B(F^2(H_n))$ is an operator algebra,  we deduce that $\Psi$ is a
completely  isometric isomorphism of operator  algebras.
 This completes the proof.
  \end{proof}

 Here is our version of the maximum principle for free holomorphic functions.

\begin{theorem}\label{max-mod1}
Let  $\cH$ be an infinite dimensional Hilbert space. Assume that $f:[B(\cH)^n]_1^-\to B(\cH)$ is  a continuous  function in the operator norm,  and  it is free holomorphic on $[B(\cH)^n]_1$. Then
\begin{equation*}
\begin{split}
\max\{\|f(X_1,\ldots, X_n)\|&:\ \|[X_1,\ldots, X_n]\|\leq 1\}\\
&=
\max\{\|f(X_1,\ldots, X_n)\|:\ \|[X_1,\ldots, X_n]\|= 1\}.
\end{split}
\end{equation*}
\end{theorem}

\begin{proof}
Due to the continuity of $f$, for any $[X_1,\ldots, X_n]\in [B(\cH)^n]_1^-$,
$$
\|f(X_1,\ldots, X_n)\|=\lim_{r\to 1} \|f(rX_1,\ldots, rX_n)\|.
$$
On the other hand, the noncommutative von Neumann inequality implies
$$
\|f(rX_1,\ldots, rX_n)\|\leq \|f(rS_1,\ldots, rS_n)\|
\quad \text{
for } \ 0\leq r<1.
$$
By Theorem \ref{A-infty}, $f\in \cA_n$ and, consequently,
$$
\lim_{r\to 1} \|f(rS_1,\ldots, rS_n)\|=\|L_f\|=\|f\|_\infty.
$$
Combining these relations, we deduce that
\begin{equation}
\label{ff}
\|f(X_1,\ldots, X_n)\|\leq \|f\|_\infty\quad \text{ for any } \ [X_1,\ldots, X_n]\in [B(\cH)^n]_1.
\end{equation}
 Since $\cH$ is infinite dimensional, there exists a subspace  $\cK\subset \cH$ and a unitary operator $U:F^2(H_n)\to \cK$.
Define the operators
$$
V_i:=\left(\begin{matrix}
US_iU^*&0\\
0&0
\end{matrix}\right), \quad i=1,\ldots,n,
$$
with respect to the orthogonal decomposition $\cH=\cK\oplus \cK^{\perp}$,
where $S_1,\ldots, S_n$  are the left creation operators on the full Fock space $F^2(H_n)$.
Notice that $\|[V_1,\ldots, V_n]\|=1$ and
$$
f(V_1,\ldots, V_n)=\lim_{r\to 1} \left(\begin{matrix}
Uf_r(S_1,\ldots, S_n)U^*&0\\
0&0
\end{matrix}\right)
$$
in the operator norm. Consequently,
$$
\|f(V_1,\ldots, V_n)\|=\lim_{r\to 1}\|f_r(S_1,\ldots, S_n)\|=\|f\|_\infty.
$$
Hence, and using  inequality \eqref{ff}, we deduce that
\begin{equation*}
\begin{split}
\max\{\|f(X_1,\ldots, X_n)\|&:\ \|[X_1,\ldots, X_n]\|\leq 1\}\\
&=
\max\{\|f(X_1,\ldots, X_n)\|:\ \|[X_1,\ldots, X_n]\|= 1\}\\
&=\|f\|_\infty.
\end{split}
\end{equation*}
This completes the proof.
 \end{proof}

\begin{corollary}\label{max-mod2}
Let $f$ be a free  holomorphic function on $[B(\cH)^n]_1$, where $\cH$ is an infinite dimensional Hilbert space, and let $r\in [0,1)$. Then
\begin{equation*}
\begin{split}
\max\{\|f(X_1,\ldots, X_n)\|&:\ \|[X_1,\ldots, X_n]\|\leq r\}\\
&=
\max\{\|f(X_1,\ldots, X_n)\|:\ \|[X_1,\ldots, X_n]\|= r\}\\
&=\|f(rS_1,\ldots, rS_n)\|.
\end{split}
\end{equation*}
\end{corollary}

 In a forthcoming paper \cite{Po-Bohr},  we obtain operator-valued
multivariable
Bohr type inequalities for   free holomorphic functions on the open  operatorial unit $n$-ball.
As consequences,  we  obtain
 operator-valued
Bohr   inequalities for
the noncommutative disc algebra $\cA_n$ and  the noncommutative analytic Toeplitz algebra $F_n^\infty$.

\bigskip

\section{Free analytic functional calculus  and noncommutative Cauchy transforms
}
\label{free analytic}

In this section,  we introduce  a  free analytic functional calculus
for $n$-tuples $T:=[T_1,\ldots, T_n]\in B(\cH)^n$ of operators
 with joint spectral radius
$r(T_1,\ldots, T_n)<1$.  We   introduce a noncommutative Cauchy transform
$\cC_T:B(F^2(H_n))\to B(\cH)$ associated with any such $n$-tuple of operators   and  prove  that
$$
f(T_1,\ldots, T_n)=C_T(f(S_1,\ldots, S_n)),\quad f\in H^\infty (B(\cX)^n_1),
$$
where  $f(S_1,\ldots, S_n)$ is the boundary function of $f$.
 Similar Cauchy representations are obtained for the
$k$-order Hausdorff derivations of $f$. Finally, we show  that the noncommutative Cauchy transform commutes with the action of the unitary group $\cU(\CC^n)$.

\begin{theorem}\label{abel}
Let $F:=\sum\limits_{\alpha\in \FF_n^+} A_{(\alpha)}\otimes Z_\alpha$ be a free holomorphic
function on the open  operatorial $n$-ball of radius $\gamma$. Then,  for any Hilbert space  $\cH$ and any $n$-tuple  of operators $[X_1,\ldots, X_n]\in B(\cH)^n$ with $r(X_1,\ldots,X_n)<\gamma$, the series
$$
F(X_1,\ldots, X_n)=\sum\limits_{k=0}^\infty \sum\limits_{|\alpha|=k} A_{(\alpha)}\otimes X_\alpha
$$ is convergent in the operator norm  of $B(\cK\otimes \cH)$.
Moreover, if \  $0<r<1$, then
\begin{equation} \label{lim-Fr}
\lim_{r\to 1}F_r(X_1,\ldots, X_n)=F(X_1,\ldots, X_n)
\end{equation}
and
\begin{equation} \label{lim-PFr}
\lim_{r\to 1}\left(\frac{\partial^k F_r}{\partial Z_{i_1}\cdots Z_{i_k}}\right)(X_1,\ldots, X_n)=\left(\frac{\partial^k F}{\partial Z_{i_1}\cdots \partial Z_{i_k}}\right)(X_1,\ldots, X_n)
 \end{equation}
for $i_1,\ldots, i_k\in \{1,\ldots, n\}$,
where the limits are in the operator  norm.
\end{theorem}

\begin{proof}.
 Assume that  $[X_1,\ldots, X_n]$ is an $n$-tuple of operators on $\cH$ such that $r(X_1,\ldots,X_n)<R$, where $R$ is the radius of convergence of $F$.
Let $\rho',\rho>0$ be such that
$r(X_1,\ldots,X_n)<\rho'<\rho<R$. Due to the definition of $r(X_1,\ldots,X_n)$, there exists $k_0\in \NN$ such that
\begin{equation}\label{ro'}
\left\|\sum\limits_{|\alpha|=k}X_{\alpha}X_\alpha^*\right\|^{1/2k}< \rho'\quad  \text{ for any }\ k\geq k_0.
\end{equation}
Since $\frac{1}{\rho}> \frac{1}{R}$, we can find $m_0$ such that
\begin{equation}\label{ro}
\left\|\sum\limits_{|\alpha|=k}A_{(\alpha)}A_{(\alpha)}^*\right\|^{1/2k}< \frac{1}{\rho}\quad  \text{ for any }\ k\geq m_0.
\end{equation}
If $k\geq \max\{k_0,m_0\}$, then relations \eqref{ro'} and \eqref{ro} imply
\begin{equation*}
\begin{split}
\left\|\sum\limits_{|\alpha|=k}A_{(\alpha)}\otimes X_\alpha \right\|&=
\left\|\left[ I\otimes X_\alpha:\ |\alpha|=k\right]
\left[\begin{matrix}
A_{(\alpha)}\otimes I\\
:\\|\alpha|=k
\end{matrix}\right]\right\|\\
&=\left\|\sum\limits_{|\alpha|=k}X_{\alpha}X_\alpha^*\right\|^{1/2}\left\|\sum\limits_{|\alpha|=k}A_{(\alpha)}A_{(\alpha)}^*\right\|^{1/2}\\
&\leq \left(\frac{\rho'}{\rho}\right)^k.
\end{split}
\end{equation*}
This proves the convergence of the series
$\sum\limits_{k=0}^\infty \left( \sum\limits_{|\alpha|=k}A_{(\alpha)}\otimes X_\alpha\right)$
  in the operator norm.
Now,  using the above inequalities, we obtain
\begin{equation*}
\begin{split}
\left\|\sum_{k=0}^\infty \sum_{|\alpha|=k} (r^{|\alpha|}-1)A_\alpha\otimes X_\alpha\right\|&\leq \sum_{k=1}^\infty (r^k-1)\left\|\sum_{|\alpha |=k}A_\alpha \otimes X_\alpha\right\|\\
&\leq \sum_{k=1}^\infty (r^k-1) \left(\frac{\rho'}{\rho}\right)^k\\
&\leq (r-1)\sum_{k=1}^\infty k \left(\frac{\rho'}{\rho}\right)^k.
\end{split}
\end{equation*}
Since $\rho'<\rho$, the latter series is convergent and therefore relation \eqref{lim-Fr} holds.
Due to Theorem \ref{derivation}, $\frac{\partial F}{\partial Z_i}$ is a free holomorphic
function on the open  operatorial $n$-ball of radius $\gamma$,
and
$$\frac{\partial^k F_r}{\partial Z_{i_1}\cdots \partial Z_{i_k}}(X_1,\ldots, X_n)=r^k
\frac{\partial^k F}{\partial Z_{i_1}\cdots \partial Z_{i_k}}(rX_1,\ldots, rX_n),\quad 0<r<1.
$$
 Applying relation \eqref{lim-Fr} to $\frac{\partial^k F}{\partial Z_{i_1}\cdots \partial  Z_{i_k}}$,  we  deduce \eqref{lim-PFr}.
The proof is complete.
\end{proof}

Let $T:=[T_1,\ldots, T_n]\in B(\cH)^n$ be an   $n$-tuple of operators with  joint spectral radius
$r(T_1,\ldots, T_n)<1$.
 We introduce the {\it Cauchy kernel} associated  with $T$  to be  the operator
$C_T(S_1,\ldots, S_n)\in B(F^2(H_n)\otimes \cH)$ defined by
\begin{equation}
\label{Cauc}
C_T(S_1,\ldots, S_n):=\sum_{k=0}^\infty\sum_{|\alpha|=k} S_\alpha\otimes T_{\tilde\alpha}^*,
\end{equation}
where $S_1,\ldots, S_n$ are the left creation operators on the full Fock space $F^2(H_n)$, and $\tilde\alpha$ is the reverse of $\alpha$, i.e.,  $\tilde \alpha= g_{i_k}\cdots g_{i_k}$ if $\alpha=g_{i_1}\cdots g_{i_k}$.
Applying Theorem \ref{Abel}, when $A_{(\alpha)}:=T_\alpha^*$, $\alpha\in \FF_n^+$ and $X_i:=S_i$, $i=1,\ldots, n$, we
deduce that
$$\frac{1}{R}=\lim_{k\to\infty}
\left\|\sum_{|\alpha|=k} T_\alpha T_\alpha^*\right\|^{1/2k}=r(T_1,\ldots, T_n)<1
$$
and $\|[S_1,\ldots, S_n]\|=1<R$. Consequently, the series in \eqref{Cauc} is convergent in the operator norm and  $C_T(S_1,\ldots, S_n)\in \cA_n\bar\otimes B(\cH)\subset  B(F^2(H_n)\otimes \cH)$.
Now, one can  easily see  that

\begin{equation}
\label{Cauc-inv}
C_T(S_1,\ldots, S_n)=\left( I-S_1\otimes T_1^*-\cdots -S_n\otimes T_n^*\right)^{-1}.
\end{equation}
We  call the operator
$$S_1\otimes T_1^*+\cdots +S_n\otimes T_n^*$$
 the {\it  reconstruction operator} associated with the $n$-tuple $[T_1,\ldots, T_n]$.
We should mention that this operator plays an important role
in  noncommutative  multivariable operator theory (see
\cite{Po-varieties}, \cite{Po-unitary}).
We remark that if $1$ is not in the spectrum of the reconstruction operator,
  then the Cauchy kernel defined by \eqref{Cauc-inv} makes sense.  In this case,  $C_T(S_1,\ldots, S_n)$ is in $F_n^\infty\bar\otimes B(\cH)$, the $WOT$-closed operator algebra generated by the spatial tensor product,  and not necessarily in $\cA_n\bar\otimes B(\cH)$. Morever, we can think of  the series $\sum_{k=0}^\infty\sum_{|\alpha|=k} S_\alpha\otimes T_{\tilde \alpha}^*$ as the Fourier representation of the Cauchy kernel.

In what follows we also use the notation  $C_T:=C_T(S_1,\ldots, S_n)$.

\begin{proposition}\label{Prop-Cauc}
Let $T:=[T_1,\ldots, T_n]\in B(\cH)^n$ be an   $n$-tuple of operators with  joint spectral radius
$r(T_1,\ldots, T_n)<1$. Then:
\begin{enumerate}
\item[(i)]
$\|C_T\|\leq \sum\limits_{k=0}^\infty\left\|\sum\limits_{|\alpha|=k} T_\alpha T_\alpha^*\right\|^{1/2}$. In particular, if  $T:=[T_1,\ldots, T_n]\in[B(\cH)^n]_1$, then
$\|C_T\|\leq \frac{1}{1-\|T\|}$.
\item[(ii)]
$C_T-C_X=C_T\left[\sum\limits_{i=1}^n S_i\otimes (T_i^*-X_i^*)\right] C_X$
and
$$
\|C_T-C_X\|\leq \|C_T\| \|C_X\|\|[T_1-X_1,\ldots, T_n-X_n]\|
$$
for any $n$-tuple   $X:=[X_1,\ldots, X_n]\in B(\cH)^n$   with  joint spectral radius
$r(X_1,\ldots, X_n)<1$.
\end{enumerate}
 \end{proposition}

\begin{proof} Since $S_1,\ldots, S_n$ are isometries with orthogonal ranges, we have
$$
\|C_T\|\leq \sum_{k=0}^\infty\left\| \sum_{|\alpha|=k}
S_\alpha\otimes T_{\tilde \alpha}^*\right\|
=
\sum_{k=0}^\infty\left\| \sum_{|\alpha|=k}
T_\alpha T_{\alpha}^*\right\|^{1/2}.
$$
If $\|[T_1,\ldots, T_n]\|<1$, then
$$\sum_{k=0}^\infty\left\| \sum_{|\alpha|=k}
T_\alpha T_{\alpha}^*\right\|^{1/2}
\leq
\sum_{k=0}^\infty\left\| \sum_{i=1}^n T_iT_i^*\right\|^{k/2}
=\frac{1}{1-\|T\|}.
$$
To prove (ii), notice that
\begin{equation*}
\begin{split}
C_T-C_X&=
\left(I-\sum_{i=1}^n S_i\otimes T_i^*\right)^{-1}
\left[ I-\sum_{i=1}^n S_i\otimes X_i^*-\left(I-\sum_{i=1}^n S_i\otimes T_i^*\right)\right]\left(I-\sum_{i=1}^n S_i\otimes X_i^*\right)^{-1}\\
&=C_T\left[\sum\limits_{i=1}^n S_i\otimes (T_i^*-X_i^*)\right] C_X,
\end{split}
\end{equation*}
and

\begin{equation*}
\begin{split}
\|C_T-C_X\|&\leq
\|C_T\| \|C_X\|\left\|
\sum_{i=1}^n S_i\otimes (T_i^*-X_i^*)\right\|\\
&=\|C_T\| \|C_X\|\left\|\sum_{i=1}^n
(T_i-X_i)(T_i-X_i)^*\right\|^{1/2},
\end{split}
\end{equation*}
which completes the proof.
\end{proof}


The {\it Cauchy transform} at $T:=[T_1,\ldots, T_n]\in B(\cH)^n$ is the mapping
$$
\cC_T:B(F^2(H_n))\to B(\cH)
$$
defined by
$$
\left< \cC_T(A)x,y\right>:=
\left<(A\otimes I_\cH)(1\otimes x), C_T(R_1,\ldots, R_n)(1\otimes y)\right>
$$
for any $x,y\in \cH$, where $R_1,\ldots, R_n$ are the right creation operators on the full Fock space $F^2(H_n)$. The operator $\cC_T(A)$ is called the Cauchy transform of $A$ at $T$. Given $A\in B(F^2(H_n))$, the Cauchy transform generates a function (the Cauchy transform of $A$)
$$
\cC[A]:[B(\cH)^n]_1\to B(\cH)
$$
by setting
$$
\cC[A](X_1,\ldots, X_n):=\cC_X(A)
\quad \text{
for any } \ X:=[X_1,\ldots, X_n]\in [B(\cH)^n]_1.
$$
Indeed, it is enough to see   that $r(X_1,\ldots, X_n)\leq
\|[X_1,\ldots, X_n]\|<1$, and therefore $\cC_X(A)$ is well-defined.
This gives rise to  an important question:  when is  $\cC[A]$  a free holomorphic function on $[B(\cH)^n]_1$.

Due to Theorem \ref{abel}, if $f=\sum_{k=0}^\infty \sum_{|\alpha|=k} a_\alpha Z_\alpha$ is a free holomorphic function on the open  operatorial unit $n$-ball and  $[T_1,\ldots, T_n]\in B(\cH)^n$ is any $n$-tuple of operators with
$r(T_1,\ldots, T_n)<1$ then,  we can define a bounded linear operator
$$
f(T_1,\ldots, T_n):=\sum_{k=0}^\infty\sum_{|\alpha|=k}
a_\alpha T_\alpha,
$$
where the series converges in norm. This provides the {\it free analytic functional calculus}.

If $F=\sum_{k=0}^\infty \sum_{|\alpha|=k} a_\alpha Z_\alpha $ is
in the Hardy algebra  $H^\infty(B(\cX)^n_1)$, we denote by
$F(S_1,\ldots, S_n)$ the boundary function of $F$, i.e.,
$F(S_1,\ldots, S_n):=L_f\in B(F^2(H_n))$, where $f:= \sum_{\alpha\in \FF_n^+} a_\alpha e_\alpha$.

\begin{theorem}\label{an=cauch}
Let $T:=[T_1,\ldots, T_n]\in B(\cH)^n$ be an   $n$-tuple of operators with  joint spectral radius
$r(T_1,\ldots, T_n)<1$. Then, for any $f\in H^\infty(B(\cX)^n_1)$,
$$
f(T_1,\ldots, T_n)=\cC_T(f(S_1,\ldots, S_n)),
$$
where $f(T_1,\ldots, T_n)$ is defined by the free analytic functional calculus, and  $f(S_1,\ldots, S_n)$ is the boundary function of $f$. Moreover,
$$
\|f(T_1,\ldots, T_n)\|\leq \left(\sum\limits_{k=0}^\infty\left\|\sum\limits_{|\alpha|=k} T_\alpha T_\alpha^*\right\|^{1/2}\right) \|f\|_\infty.
$$
\end{theorem}

\begin{proof}
 First, we prove the above equality for monomials.
Notice that
\begin{equation*}
\begin{split}
\left<\cC_T(S_\alpha)x,y\right>&=
\left< (S_\alpha\otimes I_\cH)(1\otimes x), C_T(R_1,\ldots, R_n)(1\otimes y)\right>\\
&=\left< e_\alpha\otimes x, \left(\sum_{\beta\in \FF_n^+} R_\beta\otimes T_{\tilde \beta}^*\right)(1\otimes y)\right>\\
&=
\left< e_\alpha\otimes x, \sum_{\beta\in \FF_n^+} e_{\tilde \beta}\otimes
T_{\tilde \beta}^*y\right>\\
&=\left<T_\alpha x,y\right>
\end{split}
\end{equation*}
for any $ x,y\in \cH$.
Now, assume that $f:=\sum_{k=0}^\infty\sum_{|\alpha|=k} a_\alpha Z_\alpha$ is in $H^\infty(B(\cX)^n_1)$ and   $0<r<1$. Then, due to Theorem \ref{abel}, we have
$$
\lim_{m\to\infty}\sum_{k=0}^m r^k \sum_{|\alpha|=k} a_\alpha S_\alpha=f_r(S_1,\ldots, S_n)\in \cA_n
$$
in the operator norm of $B(F^2(H_n))$,
and
$$
\lim_{m\to\infty}\sum_{k=0}^m r^k \sum_{|\alpha|=k} a_\alpha T_\alpha=f_r(T_1,\ldots, T_n)
$$
in the operator norm of $B(\cH)$.
Now,  due to the continuity of the  noncommutative Cauchy transform
in the operator norm, we deduce that
\begin{equation}\label{f_r-C}
f_r(T_1,\ldots, T_n)=\cC_T(f_r(S_1,\ldots, S_n)).
\end{equation}
Since $f(S_1,\ldots, S_n)\in F_n^\infty$, we know that
$\lim\limits_{r\to 1} f_r(S_1,\ldots, S_n)=f(S_1,\ldots, S_n)$ in the strong operator topology.
Since $\|f_r(S_1,\ldots, S_n)\|\leq \|f\|_\infty$, we deduce that
$$\text{\rm SOT}-\lim\limits_{r\to 1}f_r(S_1,\ldots, S_n)\otimes I_\cH=f(S_1,\ldots, S_n)\otimes I_\cH.
$$
On the other hand, by  Theorem \ref{abel},
$\lim\limits_{r\to 1} f_r(T_1,\ldots, T_n)=f(T_1,\ldots, T_n)$ in the operator norm.
 Passing to the limit, as $r\to 1$, in the equality
$$
\left<f_r(T_1,\ldots, T_n)x,y\right>=\left<(f_r(S_1,\ldots, S_n)\otimes I_\cH)(1\otimes x), C_T(R_1,\ldots, R_n)(1\otimes y)\right>,
\quad x, y\in \cH,
$$
we obtain $f(T_1,\ldots, T_n)=\cC_T(f(S_1,\ldots, S_n))$, which proves the first part of the theorem.
 Now, we can deduce the second part of the theorem using Proposition
\ref{Prop-Cauc}.
 This completes the proof.
   \end{proof}

Using the Cauchy representation provided by Theorem \ref{an=cauch}, one can deduce the following result.

\begin{corollary}\label{conv-u-w*}
Let $[T_1,\ldots, T_n]\in B(\cH)^n$ be an $n$-tuple of operators with  $r(T_1,\ldots, T_n)<1$.
\begin{enumerate}
\item[(i)] If $\{f_k\}_{k=1}^\infty$ and $f$ are free holomorphic functions  in $Hol(B(\cX)^n_1)$ such that $\|f_k-f\|_\infty\to 0$, as $k\to \infty$, then
$f_k(T_1,\ldots, T_n)\to f(T_1,\ldots, T_n)$ in the operator norm of $B(\cH)$.
\item[(ii)] If $\{f_k\}_{k=1}^\infty$ and $f$ are   in  the algebra $H^\infty (B(\cX)^n_1)$ such that
$f_k(S_1,\ldots, S_n)\to f(S_1,\ldots, S_n)$ in the $w^*$-topology (or strong operator topology)  and $\|f_k\|_\infty\leq M$ for any $k=1,2,\ldots$, then
$f_k(T_1,\ldots, T_n)\to f(T_1,\ldots, T_n)$ in the weak operator topology.
\end{enumerate}
\end{corollary}

We can extend Theorem \ref{an=cauch}  and obtain Cauchy representations  for the $k$-order Hausdorff derivations of
bounded  free holomorphic  functions.

\begin{theorem}\label{cauc-dif}
Let $T:=[T_1,\ldots, T_n]\in B(\cH)^n$  be an $n$-tuple of operators with  the joint  spectral radius
$r(T_1,\ldots, T_n)<1$ and  let $f \in H^\infty(B(\cX)^n_1)$.
Then
  \begin{equation}\label{deriv-Cau}
\begin{split}
\Bigl<\left(\frac{\partial^k f}{\partial Z_{i_1}\cdots \partial Z_{i_k}}\right)&(T_1,\ldots, T_n) x,y\Bigr>\\
&=
\left<\left[\frac{\partial^k \left(C_T(R_1,\ldots, R_n)^*\right)}{\partial T_{i_1}\cdots \partial T_{i_k}}\right](f(S_1,\ldots, S_n)\otimes I_\cH)(1\otimes x), 1\otimes y\right>
\end{split}
\end{equation}
for any \ $i_1,\ldots, i_k\in \{ 1,\ldots, n\}$  and $x,y\in \cH$, where  $f(S_1,\ldots, S_n)$ is the boundary function of $f$.
Moreover,
\begin{equation}\label{deriv-est}
\left\| \left(\frac{\partial f}{\partial Z_i}\right)(T_1,\ldots, T_n)\right\|\leq \|f\|_\infty
 \sum_{k=1}^\infty k^{3/2}\left\|\sum_{|\beta|=k-1}T_\beta T_\beta\right\|^{1/2}, \quad i=1,\ldots,n.
\end{equation}
\end{theorem}

\begin{proof}
First, notice that
$$
C_X(R_1,\ldots, R_n)^*=\sum_{k=0}^\infty \sum_{|\alpha|=k} R_{\tilde \alpha}^*\otimes X_\alpha,
$$
where the series is convergent in norm for each $n$-tuple
$[X_1,\ldots, X_n] $ with $r(X_1,\ldots, X_n)<1$.
Therefore,
$$G:=\sum_{k=0}^\infty \sum_{|\alpha|=k} R_{\tilde \alpha}^*\otimes Z_\alpha
$$
is a free holomorphic function on the open  operatorial unit $n$-ball. Due to Theorem \ref{derivation},  $\frac{\partial^k G}{\partial Z_{i_1}\cdots \partial Z_{i_k}}$ is also a free holomorphic function. By Theorem \ref{abel},
$\frac{\partial^k G}{\partial Z_{i_1}\cdots \partial Z_{i_k}}(X_1,\ldots, X_n)$ is a bounded operator for any $n$-tuple $[X_1,\ldots, X_n]$ with  spectral radius $r(X_1,\ldots, X_n)<1$.

 Now, notice that, for each $\alpha\in \FF_n^+$, $i=1,\ldots, n$, and $ x,y\in \cH$, we have
\begin{equation*}
\begin{split}
\Bigl<\Bigl[\frac{\partial \left(C_T(R_1,\ldots, R_n)^*\right)}{\partial T_{i}}\Bigr]&(S_\alpha\otimes I_\cH) (1\otimes x), 1\otimes y\Bigr>\\
&=
\left<
\left(\sum_{k=0}^\infty \sum_{|\beta|=k}R_\beta^*\otimes \frac{\partial T_{\tilde \beta}}{\partial T_i}\right)(S_\alpha\otimes I_\cH) (1\otimes x), 1\otimes y\right>\\
&=
\left<e_\alpha\otimes x, \sum_{k=0}^\infty \sum_{|\beta|=k}
e_{\tilde \beta}\otimes \left(\frac{\partial T_{\tilde \beta}}{\partial T_i}\right)^*y\right>\\
&=
\left< \frac{\partial T_\alpha}{\partial T_i}x,y\right>
=\left<\left(\frac{\partial Z_\alpha}{\partial Z_i}\right)(T_1,\ldots, T_n) x,y\right>.
 \end{split}
\end{equation*}
Hence, we deduce  relation \eqref{deriv-Cau} for polynomials.
Let $f=\sum_{k=1}^\infty \sum_{|\alpha|=k} a_\alpha Z_\alpha$
be in $H^\infty(B(\cX)^n_1)$.
Due to Theorem \ref{abel}, we have

$$
\left(\frac{\partial f_r}{\partial Z_i}\right)(T_1,\ldots, T_n)
=\lim_{m\to\infty}\sum_{k=0}^m \sum_{|\alpha|=k} r^{|\alpha|} a_\alpha \left(\frac{\partial Z_\alpha}{\partial Z_i}\right)(T_1,\ldots, T_n),
$$
where the convergence is in the operator  norm of $B(\cH)$, and
$$
f_r(S_1,\ldots, S_n)=\lim_{m\to\infty}\sum_{k=0}^m \sum_{|\alpha|=k} r^{|\alpha|} a_\alpha S_\alpha\in \cA_n
$$
where the convergence is in the operator norm of $B(F^2(H_n))$.
Since \eqref{deriv-Cau} holds for polynomials, the last two relations imply
\begin{equation*} \begin{split}
\Bigl<\left(\frac{\partial f_r}{\partial Z_{i}}\right)&(T_1,\ldots, T_n) x,y\Bigr>\\
&=
\left<\left[\frac{\partial \left(C_T(R_1,\ldots, R_n)^*\right)}{\partial T_{i} }\right](f_r(S_1,\ldots, S_n)\otimes I_\cH)(1\otimes x), 1\otimes y\right>
\end{split}
\end{equation*}
for any $x,y\in \cH$ and $0<r<1$.
Using again Theorem \ref{abel}, we have
$$
\lim_{r\to 1} \left(\frac{\partial f_r}{\partial Z_{i}}\right)(T_1,\ldots, T_n)=
\left(\frac{\partial f}{\partial Z_{i}}\right)(T_1,\ldots, T_n)
$$
in the operator norm.  Since  $f(S_1,\ldots, S_n)\in F_n^\infty$ (see Theorem \ref{f-infty}), as in the proof of Theorem \ref{an=cauch}, we deduce that
$$
\text{\rm SOT}-\lim_{r\to\infty} f_r(S_1,\ldots, S_n)\otimes I_\cH=f(S_1,\ldots, S_n)\otimes I_\cH.
$$
Passing to the limit, as $r\to\infty$, in the above equality, we deduce relation
\eqref{deriv-Cau} in the particular case when $k=1$.   Repeating this argument, one can prove the general case when $\frac{\partial}{\partial T_i}$
is replaced by $\frac{\partial^k}{\partial T_{i_1}\cdots \partial T_{i_k}}$.

Now, we prove the second part of the theorem.
 Notice that
\begin{equation*}
\begin{split}
\left\|\frac{\partial G}{\partial Z_i}(X_1,\ldots, X_n)\right\|
&\leq
\sum_{k=0}^\infty\left\|\sum_{|\alpha|=k}R_{\tilde \alpha}\otimes \left(\frac{\partial X_\alpha}{\partial X_i}\right)^*\right\|\\
&\leq
\sum_{k=0}^\infty\left\|\sum_{|\alpha|=k}
\left(\frac{\partial X_\alpha}{\partial X_i}\right)\left(\frac{\partial X_\alpha}{\partial X_i}\right)^*\right\|^{1/2}.
\end{split}
\end{equation*}
For each $\alpha\in \FF_n^+$, $|\alpha|=k$, we can prove that
\begin{equation}
\label{X_ga}
\left(\frac{\partial X_\alpha}{\partial X_i}\right)\left(\frac{\partial X_\alpha}{\partial X_i}\right)^*\leq k^2\sideset{}{^\alpha_d}\sum_\gamma X_\gamma X_\gamma^*,
\end{equation}
where the sum   is taken over all distinct words $\gamma$ obtained by deleting each occurence of $g_i$ in $\alpha$.
Indeed,   notice first that
$\frac{\partial X_\alpha}{\partial X_i}=\sideset{}{^\alpha}\sum\limits_\beta
X_\beta$, where the sum  is taken over all words $\beta$ obtained by deleting each occurence of $g_i$ in $\alpha$.
Since the above some contains at most $k$ terms,
one can   show that
\begin{equation*}
 \left(\frac{\partial X_\alpha}{\partial X_i}\right)\left(\frac{\partial X_\alpha}{\partial X_i}\right)^*\leq k\sideset{}{^\alpha}\sum_\beta X_\beta X_\beta^*.
\end{equation*}
Indeed,  it enough  to use the following result which is an easy consequence of the classical Cauchy inequality:
if $A_1,\ldots, A_k\in B(\cH)$, then
$$
\left(\sum_{i=1}^k A_i\right)\left(\sum_{i=1}^k A_i^*\right)
\leq k\sum_{i=1}^k A_iA_i^*.
$$
 Now,
the  $X_\beta$'s  in the above sum are not necessarily distinct but each of them can occur at most $k$ times. Consequently,
$$
\sideset{}{^\alpha}\sum\limits_\beta
X_\beta X_\beta^*\leq k  \sideset{}{^\alpha_d}\sum_\gamma X_\gamma X_\gamma^*.
$$
Combining these inequalities, we deduce  \eqref{X_ga}.
(We remark that the inequality  \eqref{X_ga} is sharp and the  equality  occurs, for example,  when $\alpha=g_i^k$.)
Therefore, we have
\begin{equation*}
\begin{split}
 \sum_{k=0}^\infty\left\|\sum_{|\alpha|=k}\left(\frac{\partial X_\alpha}{\partial X_i}\right)\left(\frac{\partial X_\alpha}{\partial X_i}\right)^*\right\|^{1/2}
&\leq
\sum_{k=0}^\infty\left\|\sum_{|\alpha|=k}
k^2\sideset{}{^\alpha_d}\sum_\gamma X_\gamma X_\gamma^*\right\|^{1/2}.
\end{split}
\end{equation*}
We remark that if $\beta\in \FF_n^+$, $|\beta|=k-1$,
then $X_\beta$ can come from free differentiation  with respect to $X_i$ of the monomials  $X_{\chi(g_i,m,\beta)}$, $m=0,1,\ldots, k-1$, where $\chi(g_i,m,\beta)$ is the insertion mapping of $g_i$ on the $m$ position of $\beta$
(see the proof of Theorem \ref{derivation}).
Consequently, we have
$$
\sum_{|\alpha|=k} \sideset{}{^\alpha_d}\sum_\gamma X_\gamma X_\gamma^*\leq k \sum_{|\beta|=k-1} X_\beta X_\beta^*.
$$
Using the above  inequalities, we obtain
$$
\left\|\frac{\partial G}{\partial Z_i}(X_1,\ldots, X_n)\right\|\leq \sum_{k=1}^\infty k^{3/2}\left\|\sum_{|\beta|=k-1}X_\beta X_\beta\right\|^{1/2}.
$$
Hence,  and due to relation \eqref{deriv-Cau}, we deduce inequality \eqref{deriv-est}. The proof is complete.
 \end{proof}

We remark that inequalities of type \eqref{deriv-Cau} can be obtained for $k$-order Housdorff derivations.
On the other hand, a similar result to Corollary \ref{conv-u-w*} can be obtain  for $k$-order Housdorff derivations, if one uses Theorem \ref{cauc-dif}.

\smallskip

In the last part of this section, we show that the noncommutative Cauchy transform commutes with certain classes of automorphisms.
Let $\cU(H_n)$ be the group of all unitaries on $H_n$ and let $U\in \cU(H_n)$.     If $U:=\left[\lambda_{ij}\right]_{i,j=1}^n$
and
 $T:=[T_1,\ldots, T_n]\in B(\cH)^n$,
we define
 $$
\beta_U(T_j):=\sum_{i=1}^n \lambda_{ij} T_i,\quad j=1,\ldots, n,
$$
and the map $\beta_U: B(\cH)^n\to B(\cH)^n$ by setting
$\beta_U(T):=[\beta_U(T_1),\ldots, \beta_U(T_n)].
$

\begin{theorem}\label{auto}
If \, $U\in \cU(H_n)$, \ $U:=\left[\lambda_{ij}\right]_{i,j=1}^n$, then the map
$\beta_U $ is an isometric  automorphism of the open unit ball $[B(\cH)^n]_1$ and also of the ball
$$
\{[T_1,\ldots, T_n]\in B(\cH)^n:\ r(T_1,\ldots, T_n)<1\}.
$$
Moreover, there is a unique  completely isometric automorphism
 of the noncommutative disc algebra $\cA_n$, denoted also by $\beta_U$,  such that
$$
\beta_U(S_j):= \sum_{i=1}^n \lambda_{ij} S_i,\quad j=1,\ldots, n,
$$
where $S_1,\ldots, S_n$ are the left creation operators
on the full Fock space.
\end{theorem}

\begin{proof} For each $j=1,\ldots, n$,
 we define the operators
$$
{\bf U}_j:=\left[\begin{matrix}
\lambda_{1j}I_\cH\\
\vdots\\
\lambda_{nj}I_\cH
\end{matrix}\right]:\cH\to \cH^{(n)},
$$
where $\cH^{(n)}$ is the direct sum of $n$ copies of $\cH$.
Notice that
\begin{equation}\label{Cuntz}
{\bf U}_i^*{\bf U}_j=\delta_{ij} I_{\cH^{(n)}}\quad \text{ and  } \quad \sum_{i=1}^n{\bf U}_i{\bf U}_i^*=I_{\cH^{(n)}}.
\end{equation}
We have $\beta_U(T)=[B_1,\ldots,B_n]$, where $B_i:=T{\bf U}_i$, \ $i=1,\ldots, n$  and \,$T:=[T_1,\ldots, T_n]$.
Now, it is clear that $\sum_{i=1}^n B_iB_i^*=\sum_{i=1}^n T_iT_i$.
If $A\in B(\cH)$ then
$$
{\bf U}_iA=\text{diag}_n (A) {\bf U}_i,\qquad i=1,\ldots, n,
$$
where $\text{diag}_n (A)$ is the $n\times n$ block diagonal operator matrix having $A$ on the diagonal and $0$ otherwise.
Using this relation and \eqref{Cuntz}, we deduce that
\begin{equation*}
\begin{split}
\sum_{|\alpha|=2} B_\alpha B_\alpha^*&=\sum_{i=1}^n B_i\left(\sum_{|\alpha|=1} B_\alpha B_\alpha^*\right) B_i\\
&=
T\left[ \sum_{i=1}^n {\bf U}_i(TT^*){\bf U}_i^*\right] T^*\\
&=
T\text{diag}_n (TT^*) \left( \sum_{i=1}^n{\bf U}_i{\bf U}_i^*
\right) T^*\\
&=T\text{diag}_n (TT^*)T=\sum_{|\alpha|=2} T_\alpha T_\alpha^*.
\end{split}
\end{equation*}
By induction over $k$, one can similarly prove that
\begin{equation}
\label{B-T}
\sum_{|\alpha|=k} B_\alpha B_\alpha^*=\sum_{|\alpha|=k} T_\alpha T_\alpha^*, \quad k=1,2,\ldots.
\end{equation}
Consequently, we have
$$
\|\beta_U(T)\|=\|T\|\quad \text{ and }\quad r(\beta_U(T))=r(T).
$$
Hence,  and since $\beta_U(T)=T{\bf U}$, where ${\bf U}:=[{\bf U}_1,\ldots, {\bf U}_n]$ is a unitary operator, we deduce  that  the map
$\beta_U:[B(\cH)^n]_1\to [B(\cH)^n]_1$  is an isometric authomorphism of the open unit ball of $B(\cH)^n$  and
$$\beta_U^{-1}(Y)=Y{\bf U}^*,\quad Y\in [B(\cH)^n]_1.
$$
Moreover, $\beta_U$ is an isometric automorphism of
the operatorial ball
$$
\{[T_1,\ldots, T_n]\in B(\cH)^n:\ r(T_1,\ldots, T_n)<1\}.
$$

  Now, let us prove the second part of the theorem.
  Using the  same  notation for the unitary  operator ${\bf U}$,
  when $\cH:=F^2(H_n)$, we deduce that
$[\beta_U(S_1),\ldots, \beta_U(S_n)]=S{\bf U}$,
where
$S:=[S_1,\ldots, S_n]$.
Setting $V_i:=\beta_U(S_i)$,\ $i=1,\ldots, n$, one can easily see that $V_1,\ldots, V_n$ are isometries with orthogonal ranges.
For any polynomial $p(S_1,\ldots, S_n)$ in the noncommutative disc algebra $\cA_n$,
we have  $\beta_U(p(S_1,\ldots, S_n))=p(V_1,\ldots, V_n)$.
According to \cite{Po-disc}, we have
$$
\|[p_{ij}(S_1,\ldots, S_n)]_m\|=\|[p_{ij}(V_1,\ldots, V_n)]_m\|.
$$
Since $\cA_n$ is the norm closure of all polynomials in $S_1,\ldots,
S_n$ and the identity, $\beta_U$  can be uniquely extended to  a
completely isometric   homomorphism from $\cA_n$ to $\cA_n$. Define
the $n$-tuple $[X_1,\ldots, X_n]:=[S_1,\ldots,S_n]{\bf U}^*$ and
notice that each  entry $X_i$ is a homogenous polynomial of degree
one in $S_1,\ldots, S_n$. Since
$$[\beta_U(X_1),\ldots, \beta_U(X_n)]=
[X_1,\ldots, X_n] {\bf U}= [S_1,\ldots,S_n],
$$
we deduce that $\beta_U(X_i)=S_i$, \ $i=1,\ldots, n$, and consequently,
$\beta_U(X_\alpha)=S_\alpha$, \ $\alpha\in \FF_n^+$.
Hence, the range of $\beta_U:\cA_n\to \cA_n$ contains all
polynomials in $\cA_n$.  Using again the norm density
of polynomials in $\cA_n$, we conclude that $\beta_U$
 is a completely  isometric automorphism of $\cA_n$.
\end{proof}

In what follows we show that the noncommutative Cauchy transform commutes with the action of the unitary group $\cU(H_n)$.

\begin{theorem}\label{Cau-inv}
Let $T:=[T_1,\ldots, T_n]\in B(\cH)^n$ be an   $n$-tuple of operators with  joint spectral radius
$r(T_1,\ldots, T_n)<1$ and $U\in \cU(H_n)$. Then
$$
\cC_T(\beta_U(f))=\cC_{\beta_U(T)}(f), \quad   f\in   \cA_n,
$$
where  $\beta_U$ is the canonical automorphism generated by $U$.
\end{theorem}

\begin{proof}
Remember that $\cA_n$ is the norm closure of the polynomials in $S_1,\ldots, S_n$ and the identity. Due to the continuity of the noncommutative Cauchy transform in the operator norm,  it is enough to prove the above  relation for $f:=S_\alpha$, $\alpha\in \FF_n^+$.
By Theorem \ref{an=cauch}, we have
\begin{equation*}
\begin{split}
\left<\cC_T(\beta_U(S_\alpha))x,y\right>&=\left<C_T(R_1,\ldots, R_n)^*(\beta_U(S_\alpha)\otimes I_\cH)(1\otimes x),1\otimes y)\right>\\
&=
\left<B_\alpha x,y\right>
\end{split}
\end{equation*}
for any $x,y\in \cH$, where $[B_1,\ldots, B_n]:=\beta_U(T)$.
On the other hand,  due to Theorem \ref{auto}, we have $r(\beta_U(T))<1$.  Applying again Theorem \ref{an=cauch},
we obtain
\begin{equation*}
\begin{split}
\left<\cC_{\beta_U(T)}(S_\alpha) x,y\right>&=\left<C_{\beta_U(T)}(R_1,\ldots, R_n)^*(S_\alpha\otimes I_\cH)(1\otimes x),1\otimes y)\right>\\
&=
\left<B_\alpha x,y\right>.
\end{split}
\end{equation*}
Hence, $\cC_T(\beta_U(S_\alpha))=\cC_{\beta_U(T)}(S_\alpha)$, and
 the result follows.
\end{proof}

The continuity    and the uniqueness of the free analytic functional calculus for
$n$-tuples of operators
with joint spectral radius $r(T_1,\ldots, T_n)<1$ will be proved in the next section.

\bigskip

\section{Weierstrass and Montel theorems for free holomorphic functions} \label{Weierstrass and Montel}

In this section, we obtain Weierstrass and Montel type theorems for the algebra   of free holomorphic functions with scalar coefficients on the open operatorial unit $n$-ball.
This enables us to introduce a metric on $Hol(B(\cX)^n_1)$ with respect to which it becomes a complete metric space, and
the Hausdorff derivations
 are continuous. In the end of this section, we prove  the continuity  and uniqueness of the free functional calculus.
Connections with  the $F_n^\infty$-functional calculus for row contractions \cite{Po-funct} and, in the commutative case,  with Taylor's functional calculus \cite{T2} are  also discussed.

We say that a sequence $\{F_m\}_{m=1}^\infty\subset Hol(B(\cX)^n_1)$  of free holomorphic functions   converges uniformly on the  closed  operatorial $n$-ball of radius $r\in [0,1)$ if it converges uniformly on the closed ball $$
[B(\cH)^n]_{ r}^{-}:=\{ [X_1,\ldots, X_n]\in B(\cH)^n:\ \|X_1X_1^*+\cdots+X_nX_n^*\|\leq r\},
$$
 where $\cH$ is an infinite dimensional Hilbert space.  According to the maximum principle of Corollary \ref{max-mod2}, this  is equivalent to the fact that the sequence $\{F_m(rS_1,\ldots, rS_n)\}_{m=1}^\infty$ is convergent in the operator norm topology of $B(F^2(H_n))$.

The first result of this section
 is a multivariable operatorial version of Weierstrass theorem (\cite{Co}).

\begin{theorem}\label{Weierstrass}
 Let $\{F_m\}_{m=1}^\infty\subset Hol(B(\cX)^n_1)$ be a sequence of free holomorphic functions   which is uniformly convergent on any closed operatorial $n$-ball of radius $r\in [0,1)$. Then there is a free holomorphic function   $F\in Hol(B(\cX)^n_1)$   such that
$F_m$ converges to $F$  on any closed operatorial $n$-ball of radius $r\in [0,1)$.

   Moreover, given $i_1,\ldots, i_k\in \{1,\ldots, n\}$,
the sequence
$\left\{\frac{\partial^k F_m} {\partial Z_{i_1}\cdots \partial Z_{i_k}}\right\}_{m=1}^\infty
$
 is uniformly convergent to $\frac{\partial^k F} {\partial Z_{i_1}\cdots \partial Z_{i_k}}$ on any closed  operatorial $n$-ball of radius $r\in [0,1)$, where $\frac{\partial^k} {\partial Z_{i_1}\cdots \partial Z_{i_k}}$ is the $k$-order Hausdorff derivation.
\end{theorem}

\begin{proof} Let $F_m:=\sum\limits_{k=0}^\infty
\sum\limits_{|\alpha|=k} a_{\alpha}^{(m)} Z_\alpha$ and
fix $r\in (0,1)$.  Then, due to Theorem \ref{caract-shifts},
$$
F_m(rS_1,\ldots, rS_n)=\sum\limits_{k=0}^\infty\sum\limits_{|\alpha|=k} r^{|\alpha|} a_{\alpha}^{(m)} S_\alpha
$$
is in the noncommutative disc algebra $\cA_n$.
Since  $\{F_m\}_{m=1}^\infty$ is uniformly convergent   on the closed operatorial $n$-ball of radius $r$, the sequence
$\{F_m(rS_1,\ldots, rS_n)\}_{m=1}^\infty$ is convergent   in the operator norm of $B(F^2(H_n))$.
 On the other hand, since  the noncommutative disc algebra $\cA_n$  is closed in the operator  norm, there exists $g\in \cA_n$ such that
\begin{equation}\label{Fm-to}
F_m(rS_1,\ldots, rS_n)\to L_g, \quad \text{ as } \ m\to\infty.
\end{equation}
 Assume $g=\sum\limits_{\alpha\in \FF_n^+} b_\alpha(r) e_\alpha$,  and notice also that
$$
b_\alpha(r)=\left< S_\alpha^* L_g(1),1\right>, \quad \alpha\in \FF_n^+.
$$
If $\lambda_{(\beta)}\in \CC$ for    $\beta\in \FF_n^+$ with  $|\beta|=k$,  we have
$$
\left|\left<\sum_{|\beta|=k} \lambda_{(\beta)} S_\beta^*(F_m(rS_1,\ldots, rS_n)-L_g)1,1\right>\right|
\leq \|F_m(rS_1,\ldots, rS_n)-L_g\|
\left\|\sum_{|\beta|=k} \lambda_{(\beta)} S_\beta^*\right\|.
$$
Since $S_1,\ldots, S_n$ are isometries with orthogonal ranges, we deduce that
$$
\left|\sum_{|\beta|=k}(r^ka_\beta^{(m)}-b_\beta(r))\lambda_{(\beta)}\right|\leq \|F_m(rS_1,\ldots, rS_n)-L_g\|
\left(\sum_{|\beta|=k} |\lambda_{(\beta)}|^2\right)^{1/2}.
$$
for any $\lambda_{(\beta)}\in \CC$ with    $|\beta|=k$.
Consequently, we have
$$
\left(\sum_{|\beta|=k}|r^ka_\beta^{(m)}-b_\beta(r)|^2\right)^{1/2}\leq  \|F_m(rS_1,\ldots, rS_n)-L_g\|
$$
for any $k=0,1,\ldots$. Since $\|F_m(rS_1,\ldots, rS_n)-L_g\|
\to 0$, as $m\to\infty$, we deduce that
$r^ka_\beta^{(m)}\to b_\beta(r)$, as $m\to\infty$, for any $|\beta|=k$ and  $k=0,1,\ldots$. Hence, $a_\beta:=\lim\limits_{m\to\infty}a_\beta^{(m)}$ exists and $b_\beta(r)=r^k a_\beta$ for any $\beta\in \FF_n^+$ with $|\beta|=k$ and $k=0,1,\ldots$.
Consider the  formal power series
$F:= \sum_{\alpha\in \FF_n^+} a_\alpha Z_\alpha$.
We show now that $F$ is a free holomorphic function on the open  operatorial unit $n$-ball.
Due to the above calculations, we have
$$
r^k \left|\left( \sum_{|\beta|=k}|a_\beta^{(m)}|^2\right)^{1/2}-
\left(\sum_{|\beta|=k}|a_\beta|^2\right)^{1/2}\right| \leq
\|F_m(rS_1,\ldots, rS_n)-L_g\|.
$$
Therefore,
\begin{equation}
\label{conv-coef}
\sum_{|\beta|=k}|a_\beta^{(m)}|^2 \to \sum_{|\beta|=k}|a_\beta|^2,\quad \text{ as } \ m\to\infty,
\end{equation}
uniformly with respect to $k=0,1,\ldots$.
Let us show that the radius of convergence of $F$ is $\geq 1$. To this end,
assume that $\gamma>1$ and
$$
\limsup_{k\to\infty} \left(\sum_{|\beta|=k}|a_\beta|^2\right)^{1/2k}>\gamma.
$$
Then there is $k\in \NN$ as large as we want such that
\begin{equation}
\label{sup-ga}
\left(\sum_{|\beta|=k}|a_\beta|^2\right)^{1/2}>\gamma^k.
\end{equation}
Choose $\lambda$ such that $1<\lambda< \gamma$ and let $\epsilon>0$ be such that
$\epsilon<\gamma-\lambda$.
Notice that $\epsilon<\gamma^k-\lambda^k$ for any $k=1,2,\ldots$. Now, due to relation \eqref{conv-coef}, there exists $N_\epsilon\in \NN$ such that
$$
\left|\left(\sum_{|\beta|=k}|a_\beta^{(m)}|^2\right)^{1/2} - \left(\sum_{|\beta|=k}|a_\beta|^2\right)^{1/2}\right|<
\epsilon
$$
for any  $m>N_\epsilon$ and any $k=0,1,\ldots$.
Hence, and using inequality \eqref{sup-ga}, we deduce that
$$
\left(\sum_{|\beta|=k}|a_\beta^{(m)}|^2\right)^{1/2}\geq \gamma^k-\epsilon>\lambda^k
$$
for any  $m>N_\epsilon$ and some $k$ as large as we want.
Consequently,  we have
$$
\limsup_{k\to \infty} \left(\sum_{|\beta|=k}|a_\beta^{(m)}|^2\right)^{1/2k}\geq \lambda>1
$$
for $m\geq N_\epsilon$. Due to Theorem \ref{Abel}, this  shows that the radius of convergence of $F_m$ is $<1$, which contradicts the fact that $F_m$ is a free holomorphic function with radius of convergence $\geq 1$. Therefore,
$$
\limsup_{k\to\infty}\left(\sum_{|\beta|=k}|a_\beta|^2\right)^{1/2k}\leq 1
$$ and, consequently, Theorem \ref{Abel} shows that $F$ is a
free holomorphic function on the open  operatorial unit ball.
The same theorem implies that
$F(rS_1,\ldots, rS_n)=\sum\limits_{k=0}^\infty \sum\limits_{|\alpha|=k} r^{|\alpha|} a_\alpha S_\alpha$ is convergent in norm.
Since $L_g$ and $F(rS_1,\ldots, rS_n)$ have the same Fourier coefficients, we must have
$L_g=F(rS_1,\ldots, rS_n)$. Due to relation \eqref{Fm-to}, we have
$$
\|F_m(rS_1,\ldots, rS_n)-F(rS_1,\ldots, rS_n)\|\to 0, \quad \text{ as }\ m\to \infty.
$$
If $[X_1,\ldots, X_n]\in [B(\cH)^n]_{1}$ and
$\|[X_1,\ldots, X_n]\|=r<1$, the noncommutative von Neumann
inequality implies
$$
\|F_m(X_1,\ldots, X_n)-F(X_1,\ldots, X_n)\|\leq
\|F_m(rS_1,\ldots, rS_n)-F(rS_1,\ldots, rS_n)\|.
$$
Taking $m\to\infty$, we deduce that
$F_m$ converges to $F$  on any closed operatorial $n$-ball of radius $r\in [0,1)$.

Now, we show that for each $\gamma\in (0,1)$
\begin{equation}
\label{co-de}
\left(\frac{\partial F_m}{\partial Z_i}\right)(\gamma S_1,\ldots, \gamma S_n)\to \left(\frac{\partial F}{\partial Z_i}\right)(\gamma S_1,\ldots, \gamma S_n)
\end{equation}
 in the operator norm, as $m\to \infty$.
Let $r,r'\in (0,1)$ such that $\gamma=r r'$. Since $(F_m)_r$ and $  F_r\in \cA_n$ are in the noncommutative disc algebra $\cA_n$, we can apply Theorem \ref{cauc-dif} (see inequality \eqref{deriv-est}) and obtain
$$
\left\|\left(\frac{\partial ((F_m)_r-F_r)}{\partial Z_i}\right)(r' S_1,\ldots, r' S_n)\right\|
\leq M\|(F_m)_r-F_r\|_\infty,
$$
where $M$ is an appropriate constant which does not depend on $m$.
Since $\|(F_m)_r-F_r\|_\infty\to 0$ as $m\to\infty$ and

$$
\left(\frac{\partial ((F_m)_r-F_r)}{\partial Z_i}\right)(r' S_1,\ldots, r' S_n)=r\left(\frac{\partial (F_m-F)}{\partial Z_i}\right)(\gamma S_1,\ldots, \gamma S_n),
$$
we deduce relation \eqref{co-de}. Using   the result  for $\frac{\partial}{\partial Z_i}$, one can obtain  the general case for $k$-order Housdorff partial derivations.
The proof is complete.
 \end{proof}

We say that a set $\cF\subset Hol(B(\cX)^n_1)$ is normal if each sequence in $\cF$ has a subsequence which converges to a function in $Hol(B(\cX)^n_1)$ uniformly on any closed operatorial ball of radius $r\in [0,1)$.
 The set $\cF$ is called  locally bounded if, for any $r\in[0,1)$, there exists $M>0$ such that
$\|f(X_1,\ldots, X_n)\|\leq M$ for any $f\in \cF$ and $[X_1,\ldots, X_n]\in [B(\cH)^n]_r$, where $\cH$ is an infinite dimensional Hilbert space.

We can prove now   the following noncommutative version of Montel theorem (see \cite{Co}).

\begin{theorem}\label{Montel}
Let $\cF\subset Hol(B(\cX)^n_1)$  be a family of free holomorphic functions. Then the  following statements are equivalent:
\begin{enumerate}
\item[(i)] $
  \sup_{f\in \cF}\|f(rS_1,\ldots, rS_n)\|<\infty
$
for each
$r\in [0,1)$.

 \item[(ii)] $\cF$ is a  normal set.
\item[(iii)] $\cF$ is  locally bounded.
\end{enumerate}
\end{theorem}

\begin{proof}
Assume that condition (i) holds. For each $f\in \cF$, let $\{a_\alpha(f)\}_{\alpha\in \FF_n^+}$ be the sequence
of coefficients. Due to (i), for each $r\in [0,1)$, there exists $M_r>0$ such that
\begin{equation}
\label{FrMr}
\|f(rS_1,\ldots, rS_n)\|\leq M_r
\quad \text{ for any } \ f\in \cF.
\end{equation}
By the Cauchy type estimate of Theorem \ref{Cauchy-est}, if $r\in(0,1)$, then
\begin{equation}
\label{Cau-est}
\left( \sum_{|\alpha|=k} |a_\alpha(f)|^2\right)^{1/2}\leq \frac{1}{r^k} M_r\quad \text{ for any } \  f\in \cF, k=0,1,\ldots.
\end{equation}
Let $\{F_m\}_{m=1}^\infty$ be a sequence  of elements in $\cF$.
Then, relation \eqref{FrMr} implies
$$
|a_0(F_m)|\leq M_0\quad \text{ for any } \ m=1,2,\ldots.
$$
Due to the classical Bolzano-Weierstrass theorem for bounded sequences of complex numbers, there is a subsequence $\{F_{m_k^{(0)}}\}_{k=1}^\infty$ of $ \{F_m\}_{m=1}^\infty$
such that the scalar sequence $\{a_0(F_{m_k^{(0)}})\}_{k=1}^\infty$ is convergent in $\CC$, as $k\to\infty$.
Inductively, using relation \eqref{Cau-est}, we find, for each $\alpha\in \FF_n^+$, $|\alpha|\geq 1$,  a subsequence
$\{F_{m_k^{(\alpha)}}\}_{k=1}^\infty$ of
$\{F_{m_k^{(\beta)}}\}_{k=1}^\infty$, where $\alpha$ is the succesor of $\beta$ in the lexicographic order of $\FF_n^+$, such that
the sequence $\{a_\alpha(F_{m_k^{(\alpha)}})\}_{k=1}^\infty$
is convergent in $\CC$, as $k\to\infty$.
Using the diagonal process, we find a subsequence $\{F_{p_k}\}_{k=1}^\infty$ of $\{F_m\}_{m=1}^\infty$ such that $\{a_\alpha(F_{p_k})\}_{k=1}^\infty$ converges in $\CC$ as $k\to\infty$, for any $\alpha\in \FF_n^+$.

Now let us prove that, if $\gamma>1$, then
 $\{F_{p_k}(\frac{r}{\gamma}S_1,\ldots, \frac{r}{\gamma}S_n)\}_{k=1}^\infty $  converges in  the norm
topology of $B(F^2(H_n))$.
Indeed, if $N\in \NN$, then  relation  \eqref{Cau-est} implies
\begin{equation*}
\begin{split}
&\Bigl\|F_{p_k}\Bigl(\frac{r}{\gamma}S_1,\ldots, \frac{r}{\gamma}S_n\Bigr)-F_{p_s}\Bigl(\frac{r}{\gamma}S_1,\ldots, \frac{r}{\gamma}S_n\Bigr)\Bigr\|\\
&\leq \sum_{j=1}^N  \frac{r^j}{\gamma^j}\left(\sum_{|\alpha|=j}
 |a_\alpha(F_{p_k})-a_\alpha(F_{p_s})|^2\right)^{1/2}
+
\sum_{j=N+1} \frac{r^j}{\gamma^j}\left(\sum_{|\alpha|=j}
|a_\alpha(F_{p_k})-a_\alpha(F_{p_s})|^2\right)^{1/2}\\
&\leq
\sum_{j=1}^N \frac{r^j}{\gamma^j}\left(\sum_{|\alpha|=j}
|a_\alpha(F_{p_k})-a_\alpha(F_{p_s})|^2\right)^{1/2}
+\sum_{j=N+1}^\infty \frac{r^j}{\gamma^j} \frac{2M_r}{r^j}\\
&\leq \sum_{j=1}^N \frac{r^j}{\gamma^j}\left(\sum_{|\alpha|=j}
|a_\alpha(F_{p_k})-a_\alpha(F_{p_s})|^2\right)^{1/2}
+\frac{2M_r}{\gamma^N(\gamma-1)}.
\end{split}
\end{equation*}
Given $\epsilon>0$, we choose $N\in \NN$ such that
$\frac{2M_r}{\gamma^N}<\frac{\epsilon}{2}$.
On the other hand, since $\{a_\alpha(F_{p_k})\}_{k=1}^\infty$ is a Cauchy sequence in $\CC$, there is $k_0\in \NN$ such that
$$
 \sum_{j=1}^N \frac{r^j}{\gamma^j}\left(\sum_{|\alpha|=j}
|a_\alpha(F_{p_k})-a_\alpha(F_{p_s})|^2\right)^{1/2}<\frac{\epsilon}{2}\quad
\text{ for any } \ k,s\geq k_0.
$$
Summing up the above results, we deduce that
$$
\Bigl\|F_{p_k}\Bigl(\frac{r}{\gamma}S_1,\ldots, \frac{r}{\gamma}S_n\Bigr)-F_{p_s}\Bigl(\frac{r}{\gamma}S_1,\ldots, \frac{r}{\gamma}S_n\Bigr)\Bigr\|
<\epsilon \quad
\text{ for any } \ k,s\geq k_0.
$$
 This proves that the the sequence
$\{F_{p_k}(\frac{r}{\gamma}S_1,\ldots, \frac{r}{\gamma}S_n)\}_{k=1}^\infty $  converges in  the norm
topology of $B(F^2(H_n))$, for any $r\in[0,1)$ and $\gamma>1$. Since the set $A:=\{\frac{r}{\gamma}:\ 0\leq r<1, \gamma>1\}$ is equal to $[0,1)$, one can choose an increasing sequence  $\{t_q\}_{q=1}^\infty$  such that $t_q\in A$ and $t_q\to 1$ as $q\to\infty$.

Now, if $\{F_m\}_{m=1}^\infty\subset \cF$, then, using the above result, there is a subsequence
$\{F_{n_k^{(1)}}\}_{k=1}^\infty$ of $\{F_m\}_{m=1}^\infty$
such that
$\{ F_{n_k^{(1)}}(t_1S_1,\ldots, t_1S_n)\}$ is convergent in  the norm
topology of $B(F^2(H_n))$, as $k\to\infty$.
Inductively, for each $q=2,3,\ldots$, we find a subsequence
$\{F_{n_k^{(q)}}\}_{k=1}^\infty$ of $\{F_{n_k^{(q-1)}}\}_{k=1}^\infty$ such that
$\{ F_{n_k^{(q)}}(t_qS_1,\ldots, t_qS_n)\}$ is convergent in  the norm
topology of $B(F^2(H_n))$, as $k\to\infty$.
Using again the diagonal process, we find a subsequence
$\{F_{m_k}\}_{k=1}^\infty$ of $\{F_m\}_{m=1}^\infty$ such that,  for each $r\in [0,1)$,
the subsequence  $\{F_{m_k}(rS_1,\ldots, rS_n)\}$ is convergent in  the norm
topology of $B(F^2(H_n))$, as $k\to\infty$.
Applying Theorem \ref{Weierstrass}, we deduce that $\cF$ is a normal set.
Therefore,  the implication $(i)\implies (ii)$ is true.

To prove the converse,
assume that there is $r_0\in (0,1)$ such that
$$
\sup_{f\in \cF}\|f(r_0S_1,\ldots, r_0S_n)\|=\infty.
$$
Let $\{f_m\}_{m=1}^\infty\subset \cF$ be such that
\begin{equation}
\label{r_0}
\|f_m(r_0S_1,\ldots, r_0S_n)\|\to\infty \quad \text{ as } \ m\to \infty.
\end{equation}
Since (ii) holds, there exists a subsequence $\{f_{m_k}\}_{k=1}^\infty$ such that   $\{f_{m_k}(rS_1,\ldots, r S_n)\}_{k=1}^\infty$
is convergent for any $r\in [0,1)$. This contradicts relation \eqref{r_0}.
The equivalence (i)$\Longleftrightarrow$(ii) follows from Corollary \ref{max-mod2}.
The proof is complete.
 \end{proof}

Now, we can obtain the following Vitali type result in our setting.
\begin{theorem}\label{Vitali}
Let $\{F_m\}_{m=1}^\infty$ be a sequence of free holomorphic functions on $[B(\cH)^n]_{1}$ with scalar coefficients  such that,
for each $r\in [0,1)$,
$$
 \sup_{m}
\|F_m(rS_1,\ldots, rS_n)\|<\infty.
$$
If there exists $0<\gamma<1$ such that $F_m(\gamma S_1,\ldots, \gamma S_n)$ converges in norm as $m\to \infty$, then $F_m$ converges uniformly  on $[B(\cH)^n]_{ r}^{-}$ for any $r\in [0, 1)$.
\end{theorem}

\begin{proof}
Suppose that $\{F_m\}_{m=1}^\infty$ does not converge uniformly on $[B(\cH)^n]_{r_0}^-$ for some $r_0\in (0,1)$.
Then there exist $\delta>0$,  subsequences
$\{F_{m_k}\}_{k=1}^\infty$  and  $\{F_{n_k}\}_{k=1}^\infty$
of $\{F_{m}\}_{m=1}^\infty$, and  $n$-tuples of operators
$[X_1^{(k)},\ldots, X_n^{(k)}]\in [B(\cH)^n]_{r_0}^-$ such that
\begin{equation}\label{Fnkmk}
\|F_{n_k}(X_1^{(k)},\ldots, X_n^{(k)})-
F_{m_k}(X_1^{(k)},\ldots, X_n^{(k)}\|\geq \delta
\end{equation}
for any $k=1,2,\ldots$. By Theorem \ref{Montel}, we find a subsequence $\{k_p\}_{p=1}^\infty$ of $\{k\}_{k=1}^\infty$
such that $\{F_{m_{k_p}}\}_{k=1}^\infty$  and  $\{F_{n_{k_p}}\}_{k=1}^\infty$ are uniformly convergent to $f$ and $g$, respectively,
 on any closed operatorial $n$-ball  of radius $r\in [0,1)$. Using Theorem \ref{Weierstrass}, we deduce that $f,g$ are free holomorphic functions on $[B(\cH)^n]_1$
 Now, the inequality \eqref{Fnkmk} and the noncommutative von Neumann inequality imply
$$
\|F_{n_{k_p}}(r_0 S_1,\ldots, r_0S_n)-F_{m_{k_p}}(r_0 S_1,\ldots, r_0S_n)\|\geq \delta>0
$$
for any $k=1,2,\ldots$.
Consequently, we have
\begin{equation}
\label{fr0gr0}
\|f(r_0 S_1,\ldots, r_0S_n)-g(r_0 S_1,\ldots, r_0S_n)\|\geq \delta>0.
\end{equation}
On the other hand, since $\{F_m(\gamma S_1,\ldots, \gamma S_n)\}_{m=1}^\infty$ converges in norm as $m\to\infty$,  we must have
$$
f(\gamma S_1,\ldots, \gamma S_n)=g(\gamma S_1,\ldots, \gamma S_n).
$$
Since $0<\gamma<1$ and $f,g$ are free holomorphic functions on $[B(\cH)^n]_1$, we deduce that $f=g$, which contradicts inequality
\eqref{fr0gr0}. The proof is complete.
 \end{proof}

\bigskip

Let $\cH$ be a Hilbert space and let $C(B(\cH)^n_1, B(\cH))$
be the vector space  of all continuous functions from  the open  operatorial unit ball
$[B(\cH)^n]_1$ to $B(\cH)$.
If $f,g\in C(B(\cH)^n_1, B(\cH))$ and $0<r<1$, we define
$$
\rho_r(f,g):=\sup_{[X_1,\ldots, X_n]\in [B(\cH)^n]_r^-}
\|f(X_1,\ldots, X_n)-g(X_1,\ldots, X_n)\|.
$$
 Let $0<r_m<1$ be such that $\{r_m\}_{m=1}^\infty$ is an increasing sequence convergent  to $1$.
For any $f,g\in C(B(\cH)^n_1, B(\cH))$, we define
$$
\rho (f,g):=\sum_{m=1}^\infty \left(\frac{1}{2}\right)^m \frac{\rho_{r_m}(f,g)}{1+\rho_{r_m}(f,g)}.
$$
Based on standard arguments, one can prove that $\rho$ is a metric on $C(B(\cH)^n_1, B(\cH))$.
 Following the corresponding result (see \cite{Co}) for   the set of all continuous functions from a set $G\subset \CC$ to a metric space $\Omega$, one can easily obtain  the following operator version. We leave the proof to the reader.

\begin{lemma}\label{Conway}
If $\epsilon>0$, then  there exists $\delta>0$ and $m\in \NN$ such that for any   $f,g\in C(B(\cH)^n_1, B(\cH))$
$$
\sup_{[X_1,\ldots, X_n]\in [B(\cH)^n]_{r_m}^-}
\|f(X_1,\ldots, X_n)-g(X_1,\ldots, X_n)\|<\delta\implies \rho(f,g)<\epsilon.
$$

Conversely, if $\delta>0$ and $m\in \NN$ are fixed, then there is  $\epsilon>0$ such that for  any  $f,g\in C(B(\cH)^n_1, B(\cH))$
$$
\rho(f,g)<\epsilon \implies \sup_{[X_1,\ldots, X_n]\in [B(\cH)^n]_{r_m}^-}
\|f(X_1,\ldots, X_n)-g(X_1,\ldots, X_n)\|<\delta.
$$
\end{lemma}
An immediate consequence of Lemma \ref{Conway} is the following: if
  $\{f_m\}_{k=1}^\infty$  and $f$  are in  $C(B(\cH)^n_1, B(\cH))$, then $f_k$ is convergent  to  $f$  in the metric $\rho$ if and only if $f_m\to f$ uniformly on any closed ball $[B(\cH)^n]_{r_m}^-$,  \ $m=1,2,\ldots$. This  result is needed to prove the following.

\begin{theorem}
\label{cont-comp}
$\left( C(B(\cH)^n_1, B(\cH)),  \rho\right)$ is a complete
metric space.
\end{theorem}

\begin{proof}
Suppose that $\{f_k\}_{k=1}^\infty$ is a Cauchy sequence in
$\left(C(B(\cH)^n_1, B(\cH)),\rho\right)$.
Due to Lemma \ref{Conway}, the sequence
$\left\{f_k|_{[B(\cH)^n]_{r}^-}\right\}_{k=1}^\infty$ is  Cauchy  in
$C([B(\cH)^n]_r^-, B(\cH))$. Consequently, for any $\epsilon>0$, there exists $N\in \NN$, such that
\begin{equation}\label{unif}
\sup_{[X_1,\ldots, X_n]\in [B(\cH)^n]_{r}^-} \|f_m(X_1,\ldots, X_n)-f_k(X_1,\ldots, X_n)\|<\epsilon\quad \text{ for any } k,m\geq N.
\end{equation}
In particular, $\{f_k(X_1,\ldots, X_n)\}_{k=1}^\infty$ is a Cauchy sequence  in  the operator norm  of $B(\cH)$. Therefore, there is an operator $f(X_1,\ldots, X_n)\in B(\cH)$ such that
\begin{equation}\label{li}
f(X_1,\ldots, X_n)=\lim_{k\to\infty} f_k(X_1,\ldots, X_n)
\end{equation}
in the operator norm. This gives rise to a function
$f:[B(\cH)^n]_1\to B(\cH)$.
We need to show that $\rho(f_k,f)\to 0$, as $k\to\infty$, and that $f$ is continuous.
    If $[X_1,\ldots, X_n]\in [B(\cH)^n]_{r}^-$,  then, due to relations \eqref{unif} and \eqref{li}, there exists $m\geq N$ such that
$$
\|f(X_1,\ldots, X_n)-f_m(X_1,\ldots, X_n)\|<\epsilon
\quad \text{ and } \quad \|f(X_1,\ldots, X_n)-f_k(X_1,\ldots, X_n)\|<\epsilon
$$ for any $k\geq N$. Since $N$ does not depend on $[X_1,\ldots, X_n]\in [B(\cH)^n]_{r}^-$, we deduce that
 $\{f_k\}_{k=1}^\infty$ converges to $f$ uniformly on any closed ball
$[B(\cH)^n]_{r}^-$. Due to Lemma \ref{Conway}, this shows that $\rho(f_k,f)\to 0$, as $k\to\infty$.
The continuity of $f$ can be proved using standard arguments
in the theory of metric spaces. We leave it to the reader.
\end{proof}

Let $\cH$ be an infinite dimensional Hilbert space and
 denote by $Hol(B(\cH)^n_1)$ the algebra of free holomorphic functions on $[B(\cH)^n]_1$.

\begin{theorem}\label{complete-metric}
$\left(Hol(B(\cH)^n_1), \rho\right)$  is a complete
metric space and  the Hausdorff derivations
$$\frac {\partial}{\partial Z_i}: \left(Hol(B(\cH)^n_1), \rho\right) \to \left(Hol(B(\cH)^n_1), \rho\right),\quad i=1,\ldots, n,
$$
are continuous.
 \end{theorem}

\begin{proof}
First, note that Theorem \ref{continuous} implies that
$Hol(B(\cH)^n_1)\subset C(B(\cH)^n_1, B(\cH)$.
Due to Theorem \ref{cont-comp}, it is enough to show that
  $\left(\cH ol(B(\cH)^n_1), \rho\right)$ is
closed in $\left( C(B(\cH)^n_1, B(\cH)),  \rho\right)$.
Let $\{f_m\}_{m=1}^\infty\subset Hol(B(\cH)^n_1)$ and $f\in C(B(\cH)^n_1, B(\cH)$ be  such that $\rho(f_m,f)\to 0$, as $m\to\infty$. Due to Lemma \ref{Conway}, $f_m\to f$ uniformly on any closed ball $[B(\cH)^n]_{r_m}^-$, \ $m=1,2,\ldots$.
Applying now Theorem \ref{Weierstrass}, we deduce that
$f\in  Hol(B(\cH)^n_1)$ and that
$$
\frac{\partial f_m}{\partial Z_i}\to \frac{\partial f}{\partial Z_i}
$$
uniformly on any closed ball $[B(\cH)^n]_{r_m}^-$ and, therefore, in the metric $\rho$. This completes the proof of the theorem.
 \end{proof}

Now, Theorem \ref{Montel} implies the following compactness criterion for subsets of $Hol(B(\cH)^n_1)$.

\begin{corollary}
A subset $\cF$ of  $(Hol(B(\cH)^n_1), \rho)$ is compact if and only if it is closed and locally bounded.
\end{corollary}

\smallskip

We return now to the setting
 of  Section \ref{free analytic}, where we showed  that
 if $f=\sum_{k=0}^\infty \sum_{|\alpha|=k} a_\alpha Z_\alpha$ is a free holomorphic function on the open  operatorial unit $n$-ball and  $[T_1,\ldots, T_n]\in B(\cH)^n$ is any $n$-tuple of operators with
$r(T_1,\ldots, T_n)<1$, then we can define the bounded linear operator
$$
f(T_1,\ldots, T_n):=\sum_{k=0}^\infty\sum_{|\alpha|=k}
a_\alpha T_\alpha,
$$
where the series converges in norm. This provides a {\it free analytic functional calculus}, which now turns out to be continuous and unique.

\begin{theorem}
 If $T:=[T_1,\ldots, T_n]\in B(\cH)^n$ is any $n$-tuple of operators with  joint spectral radius
$r(T_1,\ldots, T_n)<1$ then the mapping
$\Phi_T: Hol(B(\cX)^n_1) \to B(\cH)$ defined by
 $$
\Phi_T(f):=f(T_1,\ldots, T_n)
$$
is a continuous unital algebra homomorphism.
Moreover, the free analytic functional calculus  is uniquely
determined by the mapping
$$
Z_i\mapsto T_i,\qquad i=1,\ldots,n.
$$
\end{theorem}

\begin{proof}
Due to Theorem \ref{abel} and  Theorem \ref{operations},
we deduce that $\Phi_T$ is a  well-defined unital algebra homomorphism. To prove the continuity of $\Phi_T$, let $f_m$ and $f$ be in $Hol(B(\cX)^n_1)$ such that $f_m\to f$ in the metric  $\rho$ of $Hol(B(\cX)^n_1)$, as $m\to\infty$.
Due to Lemma \ref{Conway} and Corollary \ref{max-mod2}, this is equivalent to the fact that, for each $r\in [0,1)$,
\begin{equation}\label{conv-S}
f_m(rS_1,\ldots, rS_n)\to f(rS_1,\ldots, rS_n),\quad \text{ as }\ m\to\infty,
\end{equation}
where the convergence is  in the operator norm of $B(F^2(H_n))$.
We shall prove  that
\begin{equation}
\label{conv-f_m}
\|f_m(T_1,\ldots, T_n)-f(T_1,\ldots, T_n)\|\to 0, \quad \text{ as }\ m\to\infty.
\end{equation}
Let $f:=\sum_{k=0}^\infty \sum_{|\alpha|=k} a_\alpha Z_\alpha$
and $f_m:=\sum_{k=0}^\infty \sum_{|\alpha|=k} a_\alpha^{(m)} Z_\alpha$.
Due to Theorem \ref{abel}, the series defining $f_m(T_1,\ldots, T_n)$ and $f(T_1,\ldots, T_n)$ are  norm convergent.
Notice that
\begin{equation*}
\begin{split}
\|f_m(T_1,\ldots, T_n)-f(T_1,\ldots, T_n)\|&=
\left\|\sum_{k=0}^\infty \sum_{|\alpha|=k}(a_\alpha^{(m)}-a_\alpha)T_\alpha\right\| \\
&\leq \sum_{k=0}^\infty\left\|\sum_{|\alpha|=k}(a_\alpha^{(m)}-a_\alpha)T_\alpha\right\| \\
&\leq \sum_{k=0}^\infty
\left\|\sum_{|\alpha|=k} T_\alpha T_\alpha^*\right\|^{1/2}\left(\sum_{|\alpha|=k}|a_\alpha^{m)}-a_\alpha|^2\right)^{1/2}.
\end{split}
\end{equation*}
If  $r(T_1,\ldots, T_n)<\rho<r<1$, then there exists $k_0\in \NN$
such that
$$
\left\|\sum_{|\alpha|=k} T_\alpha T_\alpha^*\right\|^{1/2}\leq \rho^k \quad \text{ for any }\ k\geq k_0.
$$
According to Theorem \ref{Cauchy-est}, we have
$$
\left(\sum_{|\alpha|=k}|a_\alpha^{m)}-a_\alpha|^2\right)^{1/2}\leq
\frac{1}{r^k}\|f_m(rS_1,\ldots, rS_n)-f(rS_1,\ldots, rS_n)\|.
$$
Combining this with the above inequalities, we obtain
\begin{equation*}
\begin{split}
\|f_m(T_1,\ldots, T_n)-f(T_1,\ldots, T_n)\|&\leq M(T,\rho,r)
\|f_m(rS_1,\ldots, rS_n)-f(rS_1,\ldots, rS_n)\|,
\end{split}
\end{equation*}
where
$$
M(T,\rho,r):=\sum_{k=0}^{k_0}\left\|\sum_{|\alpha|=k} T_\alpha T_\alpha^*\right\|^{1/2}\frac{1}{r^k}+\sum_{k=k_0+1}^\infty\left(\frac{\rho}{r}\right)^{k}.
$$
Now, using relation \eqref{conv-S}, we deduce \eqref{conv-f_m}, which proves the continuity of $\Phi_T$.

To prove the uniqueness of the free analytic functional calculus, let $\Phi:Hol(B(\cX)_1^n)\to B(\cH)$
be a continuous unital algebra homomorphism such that
$\Phi(Z_i)=T_i$, \ $i=1,\ldots, n$. Hence, we deduce that
\begin{equation}\label{pol2}
\Phi_T(p(Z_1,\ldots, Z_n))=\Phi(p(Z_1,\ldots, Z_n))
\end{equation}
for any polynomial $p(Z_1,\ldots, Z_n)$ in $Hol(B(\cX)_1^n)$.   Let $f=\sum_{k=0}^\infty \sum_{|\alpha|=k} a_\alpha Z_\alpha$ be an element in
$Hol(B(\cX)_1^n)$ and let
$p_m:=\sum_{k=0}^m\sum_{|\alpha|=k} a_\alpha Z_\alpha$, \  $m=1,2,\ldots$.
Since
$$
f(rS_1,\ldots, rS_n)=\sum_{k=0}^\infty \sum_{|\alpha|=k} r^k a_\alpha S_\alpha$$
and the series $\sum_{k=0}^\infty
r^k\left\|\sum_{|\alpha|=k} a_\alpha S_\alpha\right\|$ converges due to Theorem \ref{caract-shifts},
we deduce that
$$
p_m(rS_1,\ldots, rS_n)\to f(rS_1,\ldots, rS_n)
$$
in the operator norm, as $m\to\infty$.
Therefore, $p_m\to f$ in the metric $\rho$ of $Hol(B(\cX)^n_1)$. Hence, using \eqref{pol2} and the continuity of $\Phi$ and $\Phi_T$, we deduce that $\Phi=\Phi_T$.
This completes the proof.
 \end{proof}

Using Theorem \ref{f-infty}, Theorem \ref{abel}, and the results from  \cite{Po-funct} concerning the $F_n^\infty$ functional calculus for row contractions, one can make  the following observation.

\begin{remark}
For strict row contractions, i.e. $\|[T_1,\ldots, T_n]\|<1$,
and $F\in H^\infty(B(\cX)^n_1)$,
the free analytic functional calculus  $F(T_1,\ldots, T_n)$ coincides with the $F_n^\infty$-functional calculus for row contractions.
\end{remark}

Let $\{F_m\}_{m=1}^\infty$ and $F$ be in $Hol(B(\cX)^n_1)$ and let $\{f_m\}_{m=1}^\infty$ and $f$  be the corresponding representations on $\CC$, respectively (see Corollary \ref{part-case}).
Due to the noncommuting von Neumann inequality, we have
$$
\sup_{|\lambda_1|^2+\cdots +|\lambda_n|^2\leq r^2} |f_m(\lambda_1,\ldots, \lambda_n)-f(\lambda_1,\ldots, \lambda_n)|\leq \|F_m(rS_1,\ldots, rS_n)-F(rS_1,\ldots, rS_n)\|
$$
for any $r\in [0,1)$. Hence, we deduce that if $F_m\to F$ in the metric $\rho$ of $Hol(B(\cX)^n_1)$, then $f_m\to f$ uniformly on compact subsets of $\BB_n$.
Since there is a sequence of polynomials $\{p_m\}_{m=1}^\infty$ such that $p_m\to F$ in the metric $\rho$, one can use the continuity of Taylor's functional calculus and the continuity of the free analytic functional calculus as well as  the fact that they coincide on polynomials, to deduce the following result.

\begin{remark}  If  $f$ is the representation of  a free holomorphic function
$F\in Hol(B(\cX)^n_1)$ on $\CC$ and $[T_1,\ldots, T_n]\in B(\cH)^n$ is an $n$-tuple of commuting operators with Taylor spectrum $\sigma(T_1,\ldots, T_n)\subset \BB_n$,
then the  free analytic calculus $F(T_1,\ldots, T_n)$ coincides with  Taylor's functional calculus $f(T_1,\ldots, T_n)$.
\end{remark}

\bigskip

\section{Free pluriharmonic functions and  noncommutative Poisson transforms}
\label{free harmonic}

Given an operator $A\in B(F^2(H_n))$, the noncommutative Poisson transform \cite{Po-poisson} generates a function
$$
P[A]: [B(\cH)^n]_1\to B(\cH).
$$
In this section,  we provide classes of operators $A\in B(F^2(H_n))$
such that $P[A]$ is a free holomorphic (resp. pluriharmonic)
function on $[B(\cH)^n]_1$.  We characterize the free holomorphic
functions $u$ on $[B(\cH)^n]_1$ such that $u=P[f]$ for some boundary
function $f$ in the noncommutative analytic Toeplitz algebra
$F_n^\infty$, or the noncommutative disc algebra $\cA_n$. We also
obtain noncommutative multivariable versions of Herglotz theorem and
Dirichlet extension problem (see \cite{Co}, \cite{H}), for free
pluriharmonic functions.

We define the  operator $K_T(S_1,\ldots, S_n)\in B(F^2(H_n)\otimes \cH)$ associated with a row contraction
$T:=[T_1,\ldots, T_n]\in B(\cH)^n$
by setting
$$
K_T(S_1,\ldots, S_n):=\sum_{k=0}^\infty\sum_{|\alpha|=k} S_\alpha \otimes \Delta_T T_\alpha^*,
$$
where $\Delta_T:=(I_\cH-\sum_{i=1}^n T_iT_i^*)^{1/2}$.
Due to Theorem \ref{Abel}, when $A_{(\alpha)}:=\Delta_T T_\alpha^*$ and $X_i:=S_i$, $i=1,\ldots, n$,
the above series is  convergent in the operator norm  if
\begin{equation}
\label{cond-conv}
\limsup_{k\to\infty} \left\|\sum_{|\alpha|=k} T_\alpha T_\alpha^*-\sum_{|\alpha|=k+1} T_\alpha T_\alpha^*
\right\|^{1/2k}<1.
\end{equation}
In particular, if $\|[T_1,\ldots, T_n]\|<1$, then  relation
\eqref{cond-conv} holds and the operator $K_T(S_1,\ldots, S_n)$ is in $\cA_n\bar \otimes B(\cH)$. Notice  also  that
$$
(S_\alpha^*\otimes I_\cH)K_T(S_1,\ldots, S_n)=K_T(S_1,\ldots, S_n) (I_{F^2(H_n)} \otimes T_\alpha^*),\qquad \alpha\in \FF_n^+.
$$
 Introduced in  \cite{Po-poisson},
the noncommutative Poisson transform at $T:=[T_1,\ldots, T_n]$ is the map
 $P_T:B(F^2(H_n))\to B(\cH)$ defined by
\begin{equation*}
\begin{split}
\left<P_T(A)x,y\right>&:=
\left<K_T(S_1,\ldots, S_n)^* ( A\otimes I_\cH) K_T(S_1,\ldots, S_n) (1\otimes x),1\otimes y\right>\\
&:=\left< K_T^*( A\otimes I_\cH) K_Tx,y\right>
\end{split}
\end{equation*}
for any $x,y\in B(\cH)$, where
$K_T:=K_T(S_1,\ldots, S_n)|_{1\otimes \cH}:\cH\to F^2(H_n)\otimes \cH$. We recall  that
the Poisson kernel $K_T$ is an isometry   if $\|T\|<1$, and
\begin{equation} \label{pol}
p(T_1,\ldots, T_n)=K_T^*(p(S_1,\ldots, S_n)\otimes I_\cH)K_T
\end{equation}
for any polynomial $p$. We refer to \cite{Po-poisson}, \cite{Po-curvature}, \cite{Po-similarity}, and \cite{Po-unitary} for more on noncommutative Poisson transforms on $C^*$-algebras generated by isometries.

Given an operator
$A\in B(F^2(H_n))$, the noncommutative Poisson transform
generates a function
$$
P[A]:[B(\cH)^n]_1\to B(\cH)
$$
by setting
$$
P[A](X_1,\ldots, X_n):=P_X(A)\quad  \text{ for }\ X:=[X_1,\ldots, X_n]\in
[B(\cH)^n]_1.
$$
In what follows, we provide   classes of operators     $A\in B(F^2(H_n))$  such that the mapping   $P[A]$ is a free holomorphic function  on $[B(\cH)^n]_1$. In this case, the operator $A$ can be seen as the boundary function of the Poisson transform $P[A]$.

As in the previous sections,  we  identify $f\in F_n^\infty$ with the multiplication operator \,$L_f\in B(F^2(H_n))$.

\begin{theorem}\label{behave} Let $\cH$ be a  Hilbert space and  $u$ be  a free holomorphic function on $[B(\cH)^n]_1$.
\begin{enumerate}
\item[(i)]
There exists $f\in F_n^\infty$ with $u=P[f]$ if and only if
$\sup\limits_{0\leq r<1}\|u(rS_1,\ldots, rS_n)\|<\infty$.
In this case, $u(rS_1,\ldots, rS_n)\to f$,  as $ r\to 1$,
  in the $w^*$-topology (or strong operator topology).

\item[(ii)]
There exists $f\in \cA_n$ with $u=P[f]$ if and only if
$\{u(rS_1,\ldots, rS_n)\}_{0\leq r<1}$ is
convergent in norm as, $r\to 1$. In this case, $u(rS_1,\ldots, rS_n)\to f$ in the operator norm, as $r\to 1$.
\end{enumerate}
\end{theorem}

\begin{proof}
To prove (i), assume that $f\in F_n^\infty$ and $u=P[f]$, where $f$ is identified with the multiplication operator $L_f\in B(F^2(H_n)$. Then
$$u(X_1,\ldots, X_n)=K_X^*(L_f\otimes I_\cH)K_X,\quad [X_1,\ldots, X_n]\in [B(\cH)^n]_1$$
and  $\|u(X_1,\ldots, X_n)\|\leq \|L_f\|=\|f\|_\infty$
for any $[X_1,\ldots, X_n]\in [B(\cH)^n]_1$.
In particular,
\begin{equation}\label{sup-u}\sup\limits_{0\leq r<1}\|u(rS_1,\ldots, rS_n)\|\leq \|f\|_\infty<\infty.
\end{equation}
Conversely, assume that $u(X_1,\ldots, X_n):=\sum_{k=0}\sum_{|\alpha|=k} a_\alpha X_\alpha$ is a free holomorphic  function on $[B(\cH)^n]_1$ such that \eqref{sup-u} holds. By Theorem \ref{f-infty}, $f:=\sum_{\alpha\in \FF_n^+} a_\alpha e_\alpha$ is in $F_n^\infty$.
Due to Theorem  \ref{Abel}, we have that
$u_r(X_1,\ldots, X_n):=\sum_{k=0}^\infty \sum_{|\alpha|=k} r^{|\alpha|} a_\alpha X_\alpha$ is convergent in norm for any $[X_1,\ldots, X_n]\in [B(\cH)^n]_1$ and $r\in [0,1]$.
Similarly,  we have that $f_r(S_1,\ldots, S_n):=\sum_{k=0}^\infty \sum_{|\alpha|=k} r^{|\alpha|} a_\alpha S_\alpha$ is convergent in norm for any $r\in [0,1)$.
 Using relation \eqref{pol}, we deduce that
$$
\sum_{k=0}^m \sum_{|\alpha|=k} r^{|\alpha|} a_\alpha X_\alpha
=K_X^*\left(\sum_{k=0}^m \sum_{|\alpha|=k} r^{|\alpha|} a_\alpha S_\alpha\otimes I_\cH\right)K_X.
$$
Taking $m\to \infty$ and using the above convergences, we get
\begin{equation}\label{u_r-f_r}
u_r(X_1,\ldots, X_n)=K_X^*(f_r(S_1,\ldots, S_n)\otimes I_\cH) K_X,\qquad r\in [0,1).
\end{equation}
By Theorem \ref{continuous}, we have
$$
\lim_{r\to 1} u_r(X_1,\ldots, X_n)=u(X_1,\ldots, X_n)
$$
in the operator norm. On the other hand,  due to
relation
\eqref{So}, we have
\begin{equation}
\label{So2}
\text{\rm SOT-}\lim_{r\to 1} f_r(S_1,\ldots, S_n)=L_f.
\end{equation}
Since $\|f_r(S_1,\ldots, S_n)\|\leq \|f\|_\infty$ and the map $A\mapsto A\otimes I_\cH$ is SOT-continuous on bounded subsets of $B(F^2(H_n))$, we take $r\to 1$ in relation \eqref{u_r-f_r}
and deduce that
$u(X_1,\ldots, X_n)=P_X(f)$ for any $[X_1,\ldots, X_n]\in [B(\cH)^n]_1$. Since $u_r(S_1,\ldots, S_n)=f(rS_1,\ldots, rS_n)$ and the strong operator topology coincides with the $w^*$-topology on $F_n^\infty$ (see \cite{DP1}), one can use \eqref{So2} to complete the proof of part (i).

To prove (ii), assume that $f=\sum_{\alpha\in \FF_n^+} a_\alpha e_\alpha$ is in $\cA_n$ and $u=P[f]$, i.e.,
$$u(X_1,\ldots, X_n)= K_X^* (L_f\otimes I_\cH) K_X
$$ for any
$X=[X_1,\ldots, X_n]\in [B(\cH)^n]_1$.
Due to Theorem \ref{A-infty}, we have $\lim\limits_{r\to 1} f_r(S_1,\ldots, S_n)=L_f$ in the operator norm.
Hence,  using relation \eqref{pol} and Theorem \ref{continuous}, we deduce that
\begin{equation*}
\begin{split}
K_X^* (L_f\otimes I_\cH) K_X&=\lim_{r\to 1}K_X (f_r(S_1,\ldots, S_n)\otimes I_\cH) K_X\\
&=\lim_{r\to 1} f(rX_1,\ldots, rX_n)=f(X_1,\ldots, X_n).
\end{split}
\end{equation*}
This proves that $u(X_1,\ldots, X_n)=f(X_1,\ldots, X_n)$
for any $[X_1,\ldots, X_n]\in [B(\cH)^n]_1$.
In particular, we deduce that
$$
u(rS_1,\ldots, rS_n)=f_r(S_1,\ldots, S_n)\to L_f, \quad \text{ as } \ r\to1,
$$
in the operator norm.

Conversely, assume that $u:=\sum_{k=0}\sum_{|\alpha|=k} a_\alpha Z_\alpha$ is a free holomorphic function on the open  operatorial unit $n$-ball, such that
$\{u(rS_1,\ldots, rS_n)\}_{0\leq r<1}$ is
convergent in norm, as $r\to 1$. By Theorem \ref{caract-shifts}, we have that $u(rS_1,\ldots, rS_n)\in \cA_n$. Since $\cA_n$ is a Banach algebra, there exists $f\in \cA_n$ such that
$
u(rS_1,\ldots, rS_n)\to f$ in norm, as $r\to 1$. Due to Theorem \ref{A-infty}, we must have  $f=\sum_{k=0}^\infty\sum_{|\alpha|=k} a_\alpha e_\alpha$.
As in the proof of part (i), we have
$$
u(X_1,\ldots, X_n)=\lim_{r\to 1} f_r(X_1,\ldots, X_n)=
\lim_{r\to 1}  K_X^*(f_r(S_1,\ldots, S_n)\otimes I_\cH)K_X
$$
for any $[X_1,\ldots, X_n]\in [B(\cH)^n]_1$.
Now, since
$\lim\limits_{r\to 1} f_r(S_1,\ldots, S_n)=L_f$ in norm, we deduce that
$u=P[f]$.
This completes the proof.
  \end{proof}

We now turn our attention to a noncommutative generalization of the
harmonic functions on the open unit disc $\DD$.
 We say that $G$ is a self-adjoint  free pluriharmonic function on $[B(\cH)^n]_1$ if there exists a free holomorphic function $F$ on $[B(\cH)^n]_1$ such that
$$
G(X_1,\ldots, X_n)=\text{\rm Re}\,F(X_1,\ldots, X_n):=\frac{1}{2}\left(F(X_1,\ldots, X_n)+ F(X_1,\ldots, X_n)^*\right)
$$
We remark that if $\cH$ be an infinite dimensional Hilbert space,
then $G$  determines $F$ up to an imaginary complex number. Indeed,
if we assume that $\text{\rm Re}\, F=0$ and take the representation
on the full Fock space $F^2(H_n)$, we obtain $F(rS_1,\ldots,
rS_n)=-F(rS_1,\ldots, rS_n)^*$,\ $0<r<1$. If $F(rS_1,\ldots, rS_n)$
has the representation $\sum_{k=0}^\infty\sum_{|\alpha|=k}
r^{|\alpha|} a_\alpha S_\alpha$,\ $a_\alpha\in \CC$, the above
relation implies
$$
\sum_{k=0}^\infty\sum_{|\alpha|=k} r^{|\alpha|} a_\alpha e_\alpha=F(rS_1,\ldots, rS_n)1=-F(rS_1,\ldots, rS_n)^*1=-\overline{a}_0.
$$
Hence, $a_\alpha=0$ if $|\alpha|\geq 1$ and $a_0+\overline{a}_0=0$. Therefore, $F=a_0$, where $a_0$ is an imaginary complex number.
This proves our assertion.
  Due to Theorem \ref{Abel},
$$G(X_1,\ldots, X_n):=\sum_{k=1}^\infty \sum_{|\alpha|=k} \overline{a}_\alpha X_\alpha^* +a_0 I+ \sum_{k=1}^\infty \sum_{|\alpha|=k} {a_\alpha} X_\alpha
$$
represents a self-adjoint free pluriharmonic function on
$[B(\cH)^n]_1$ if and only if
$$
\limsup_{k\to\infty}\left(\sum_{|\alpha|=k}|a_\alpha|^2\right)^{1/2k}\leq 1.
$$
If $H_1$ and $H_2$ are self-adjoint free pluriharmonic functions on
$[B(\cH)^n]_1$, we say that $H:=H_1+iH_2$ is a free pluriharmonic
function on $[B(\cH)^n]_1$. Notice that any free holomorphic
function on $[B(\cH)^n]_1$ is a free pluriharmonic function. This is
due to the fact that $f=\frac{f+f^*}{2}+i\frac{f-f^*}{2i}$.

\begin{proposition}
Let $g$ be a free pluriharmonic function on  the open operatorial
$n$-ball of radius $1+\epsilon$, $\epsilon>0$.  Then
$$
g(X_1,\ldots, X_n)=P_X(g(S_1,\ldots, S_n)),\quad X:=[X_1,\ldots, X_n]\in [B(\cH)^n]_1,
$$
where $P_X$ is the noncommutative Poisson transform at $X$.
Moreover, if $\cH$ is an infinite dimensional Hilbert space, then
 $g(S_1,\ldots, S_n)\geq 0$ if and only if
$g(X_1,\ldots, X_n)\geq 0$ for any  $[X_1,\ldots, X_n]\in [B(\cH)^n]_1$.
\end{proposition}

\begin{proof} Without loss of generality,
we can assume that  $g$ is a self-adjoint  free pluriharmonic
function and $g(X_1,\ldots, X_n)=f(X_1,\ldots, X_n)+f(X_1,\ldots,
X_n)^*$ for any $[X_1,\ldots, X_n)]\in [B(\cH)^n]_{1+\epsilon}$,
where the function $ f(X_1,\ldots, X_n)=\sum_{k=0}^\infty
\sum_{|\alpha|=k} a_\alpha X_\alpha$ is  free holomorphic on
$[B(\cH)^n]_{1+\epsilon}$. According to Theorem \ref{caract-shifts},
the series $\sum_{k=0}^\infty \sum_{|\alpha|=k} r^{|\alpha|}
a_\alpha S_\alpha$ converges in the operator norm for any $r\in
[0,1+\epsilon )$. Due to   relation \ref{pol} and taking limits in
the operator norm, we have
\begin{equation*}
\begin{split}
f(X_1,\ldots, X_n)&=\sum_{k=0}^\infty \sum_{|\alpha|=k} a_\alpha X_\alpha=
P_X[f(S_1,\ldots, S_n)] \ \text{ and}\\
f(X_1,\ldots, X_n)^*&=\sum_{k=0}^\infty \sum_{|\alpha|=k} \overline{a}_\alpha X_\alpha^*=
P_X[f(S_1,\ldots, S_n)^*].
\end{split}
\end{equation*}
Consequently,
$$g(X_1,\ldots, X_n)=P_X[g(S_1,\ldots,S_n)], \quad
[X_1,\ldots, X_n]\in [B(\cH)^n]_1.
$$

We prove now the last part of the proposition. One implication is obvious due to  the above
relation.   Conversely, assume that $g(X_1,\ldots, X_n)\geq 0$ for any  $[X_1,\ldots, X_n]\in [B(\cH)^n]_1$.
Then, since $\cH$ is infinite dimensional,  we deduce that
$g(rS_1,\ldots, rS_n)\geq 0$ for any $r\in [0,1)$. On the other hand, due to Theorem \ref{continuous}, $\lim\limits_{r\to 1} g(rS_1,\ldots, rS_n)=g(S_1,\ldots, S_n)$ in the operator norm.
Hence, $g(S_1,\ldots, S_n)\geq 0$, and the proof is complete.
 \end{proof}

Now, we obtain a  noncommutative  multivariable version of
Herglotz theorem (see \cite{H}).

\begin{theorem}
\label{Herglotz} Let $f\in (F_n^\infty)^*+ F_n^\infty$ and let
$u=P[f]$ be its  noncommutative Poisson transform. Then $u$ is a
free   pluriharmonic function on $[B(\cH)^n]_1$, where $\cH$ is a
Hilbert space. Moreover,  $u\geq 0$ on $[B(\cH)^n]_1$, where $\cH$
is an  infinite dimensional Hilbert space, if and only if $f\geq 0$.
\end{theorem}

\begin{proof} First, notice that, without loss of generality,  we can assume that $f=f^*$.
Then, one  can prove that $f=g^*+g$ for some $g\in F_n^\infty$.
Indeed, if $f=h^*+g$ for some $h,g\in F_n^\infty$, the we must have $(g-h)^*=g-h$. Hence, $(g-h)^*1=(g-h)1$   and one can easily   deduce that $g-h$ is a constant, which proves our assertion.
According to Theorem \ref{behave},  $P[g]$ is a free holomorphic function on the open  operatorial unit $n$-ball.
On the other hand, due to \cite{Po-varieties}, we have
$$
\text{\rm SOT}-\lim_{r\to 1} g_r(S_1,\ldots, S_n)^*=L_g^*.
$$
Hence, using the properties of the Poisson tranform and Theorem \ref{continuous}, we deduce that \begin{equation*}
\begin{split}
\left< P[g^*] x,y\right>&=
\lim_{r\to 1}\left< K_X(g_r(S_1,\ldots, S_n)^*\otimes I_\cH)K_X x,y\right>\\
&=
\lim_{r\to 1}\left< g_r(X_1,\ldots, X_n)^*x,y\right>\\
&=
\left< g(X_1,\ldots, X_n)^* x,y\right>\\
&=\left< P[g]^*x,y\right>.
\end{split}
\end{equation*}
Hence,  we have $P[g]^*=P[g^*]$.
Consequently,
$$
u=P[f]=P[g^*]+P[g]=P[g]^*+P[g],
$$
which proves that $u$ is a  self-adjoint free   pluriharmonic
function on $[B(\cH)^n]_1$.

Now, it is clear that if $f\geq 0$ then $u=P[f]\geq 0$. Conversely, assume that $u(X_1,\ldots, X_n)\geq 0$ for any
$[X_1,\ldots, X_n]\in [B(\cH)^n]_1$. Since $\cH$  is an infinite dimensional Hilbert space $\cH$, we deduce that
$$
u(rS_1,\ldots, rS_n)=g(rS_1,\ldots, rS_n)^*+ g(rS_1,\ldots, rS_n)\geq 0,\qquad r\in [0,1).
$$
Due to Theorem  \ref{behave}, we have
$$
\text{\rm WOT}-\lim_{r\to 1} [g(rS_1,\ldots, rS_n)^*+ g(rS_1,\ldots, rS_n)]=L_g^*+L_g\geq 0.
$$
Under the identification of  $g$ with $L_g$, we deduce
 $f=g^*+g\geq 0$, and complete the proof.
 \end{proof}

Here again, we  remark that $ f$ plays the role of the boundary function from the classical complex analysis.

Our   version of the classical Dirichlet extension problem for the unit disc (see \cite{Co}, \cite{H}) is the following extension of Theorem \ref{A-infty}.

\begin{theorem}\label{Dirichlet}
If $f\in \cA_n^*+\cA_n$, then  $u:=P[f]$ is  a  free pluriharmonic
function on the open  operatorial unit $n$-ball such that
\begin{enumerate}
\item[(i)] $u$ has a continuous extension $\tilde u$ to $[B(\cH)^n]_1^-$ for  any Hilbert space $\cH$, in the operator norm;
\item[(ii)]
$\tilde u(S_1,\ldots, S_n)=f$.
\end{enumerate}
\end{theorem}

\begin{proof} Without loss of generality, we can assume that $f$ is self-adjoint.
As in the proof of Theorem \ref{Herglotz}, one can prove that
$f=g^*+g$ for some $g\in \cA_n$ and $u:=P[f]=P[g]^*+P[g]$ is a
self-adjoint pluriharmonic function on the open  operatorial unit
$n$-ball. Since $g\in \cA_n$, we know that $g_r(S_1,\ldots, S_n)\to
L_g$ in norm, as $r\to 1$. Consequently,
$$f_r(S_1,\ldots, S_n):=g_r(S_1,\ldots, S_n)^*+g_r(S_1,\ldots, S_n)\to L_f^*+ L_f,\quad \text{ as } \ r\to1,
$$ in norm.
As in the proof  of Theorem \ref{Herglotz}, we have
$$
u(X_1,\ldots, X_n)=f(X_1,\ldots, X_n),\quad \text{ for } \ [X_1,\ldots, X_n]\in [B(\cH)^n]_1.
$$
 Moreover,   $v:=P[g]$ is a free holomorphic function such that $v(X_1,\ldots, X_n)=g(X_1,\ldots, X_n)$,\  for any $ [X_1,\ldots, X_n]\in [B(\cH)^n]_1$.

 For each $n$-tuple $[Y_1,\ldots, Y_n]\in [B(\cH)^n]_1^-$, we define
$$
\tilde v(Y_1,\ldots, Y_n):=\lim_{r\to 1} P_{rY}[g],
$$
where $rY:=[rY_1,\ldots, rY_n]$.
Hence, we have  $\tilde v(Y_1,\ldots, Y_n)=\lim_{r\to 1} g(rY_1,\ldots, rY_n)$.
Now, as in the proof of Theorem \ref{A-infty},  we deduce that
the map $\tilde v:[B(\cH)^n]_1^-\to B(\cH)$ is a continuous
extension of $v$.
 Therefore, the map  $\tilde u:={\tilde v}^*+\tilde v$
is a continuous extension of $u$ to $[B(\cH)^n]_1^-$.
To prove (ii), apply part (i) when $\cH=F^2(H_n)$ and take into account Theorem \ref{A-infty}.
We obtain
$$
\tilde v(S_1,\ldots, S_n)=\lim_{r\to 1} g(rS_1,\ldots, rS_n)=g,
$$
where we used the identification of $g$ with $L_g$, and the limit is in the operator norm.
Therefore,
$$
\tilde u(S_1,\ldots, S_n)=\tilde v(S_1,\ldots, S_n)^*+\tilde v(S_1,\ldots, S_n)=g^*+g=f.
$$
This completes the proof.
\end{proof}

Let $u$  and $v$ be  two  self-adjoint free pluriharmonic functions
on $[B(\cH)^n]_1$.  We say that $v$ is the pluriharmonic conjugate
of $u$ if  $u+iv$ is a free holomorphic function on $[B(\cH)^n]_1$.

\begin{remark}
The pluriharmonic conjugate  of a self-adjoint free pluriharmonic
function on $[B(\cH)^n]_1$ is unique up to an additive real
constant.
\end{remark}

\begin{proof} Let $f$ be a free holomorphic function on $[B(\cH)^n]_1$ and $u=\text{\rm Re}\, f$.
 Assume that $v$ is a selfadjoint free pluriharmonic function such that $u+iv=g$ is a
free holomorphic function on $[B(\cH)^n]_1$. Hence, we have
\begin{equation}\label{v}
v=\frac{2g-f-f^*}{2i}.
\end{equation}
Since $v=v^*$, we must have $(g-f=(g-f)^*$, i.e.,
$\text{\rm Re}\, (g-f)=0$. Based on the remarks following Theorem \ref{behave}, we have $g-f=w$, where $w$ is an imaginary complex number. Consequently,
relation \eqref{v}, implies
$v=\frac{f-f^*}{2i}-iw$.
This proves the assertion.
\end{proof}

We remark that if $u=\text{\rm Re}\, f$ and $f(0)$ is real then
$v=\frac{f-f^*}{2i}$ is the unique pluriharmonic conjugate of $u$
such that $v(0)=0$.

\begin{theorem}\label{cauch-conj}
Let $T:=[T_1,\ldots, T_n]\in B(\cH)^n$ be an   $n$-tuple of operators with  joint spectral radius
$r(T_1,\ldots, T_n)<1$. If $f\in H^\infty (B(\cX)^n_1) $, \
$u=\text{\rm Re}\,f$, and $f(0)$ is real, then
\begin{equation*}
  \left<f(T_1,\ldots, T_n)x,y\right>=\left<(u(S_1,\ldots, S_n)\otimes I_\cH)(1\otimes x),
[2C_T(R_1,\ldots, R_n)-I](1\otimes y)\right>
\end{equation*}
for any $x,y\in \cH$, where $u(S_1,\ldots, S_n)$ is the boundary function of $u$.
 \end{theorem}

\begin{proof}
Due to Theorem \ref{an=cauch}, we have
\begin{equation*}
\begin{split}
\left<(f(S_1,\ldots, S_n)\otimes I_\cH)(1\otimes x)\right. &,
\left.
[2C_T(R_1,\ldots, R_n)-I](1\otimes y)\right>\\
&=
2\left<(f(S_1,\ldots, S_n)\otimes I_\cH)(1\otimes x), [C_T(R_1,\ldots, R_n)](1\otimes y)\right>\\
 &\qquad-\left<(f(S_1,\ldots, S_n)\otimes I_\cH)(1\otimes x), 1\otimes y\right>\\
&=
2\left< f(T_1,\ldots, T_n)x,y\right>-f(0)\left< x,y\right>.
\end{split}
\end{equation*}
On the other hand, it is easy to see that
\begin{equation*}
\begin{split}
\left<(f(S_1,\ldots, S_n)^*\otimes I_\cH)(1\otimes x)\right.&,\left.
 [2C_T(R_1,\ldots, R_n)-I](1\otimes y)\right>\\
&=
\left<(\overline{(f(0)}\otimes I_\cH)(1\otimes x),
[2C_T(R_1,\ldots, R_n)-I](1\otimes y)\right>\\
 &=\overline{(f(0)}\left<x,y\right>.
\end{split}
\end{equation*}
If $f(0)\in \RR$, then adding up the above relations, we
complete the proof.
\end{proof}

We remark that under the conditions of Theorem \ref{cauch-conj} and  using the noncommutative Cauchy transform, one can  express
  the pluriharmonic conjugate of $u$ in terms of $u$.

In  a forthcoming paper \cite{Po-Bohr},
we will consider
operator-valued Bohr type inequalities for
 classes of free  pluriharmonic
functions on the open  operatorial unit $n$-ball   with operator-valued coefficients.

\bigskip

\section{ Hardy spaces of free holomorphic functions }
\label{Banach}

In this section, we define the radial  maximal Hardy space
$H^p(B(\cX)^n_1)$, $p\geq 1$, and   the symmetrized Hardy space
$H^\infty_{\text{\rm sym}}(\BB_n)$,    and prove that they are Banach spaces with respect to  some appropriate norms. In this setting, we  obtain  von Neumann type inequalities for $n$-tuples of operators.

Let $F$ be a free holomorphic function on the open  operatorial unit $n$-ball. The map $\varphi:[0,1)\to B(F^2(H_n))$ defined by $\varphi(r):= F(rS_1,\ldots, rS_n)$ is called
the {\it radial boundary function} associated with $F$.
Due to Theorem \ref{continuous}, $\varphi$ is continuous
 with respect to the operator norm topology of $B(F^2(H_n))$.
When $\lim\limits_{r\to 1} \varphi(r)$ exists, in one of the classical topologies of $B(F^2(H_n))$, we  call it the {\it boundary function}  of  $F$.

Due to the  maximum principle for free holomorphic functions (see  Theorem \ref{max-mod1}), we have
$$
\|\varphi(r)\|=\sup \|F(X_1,\ldots, X_n)\|,\quad 0\leq r<1,
$$
where the supremum is taken over all $n$ tuples of operators $[X_1,\ldots, X_n]$ in either one of the following sets
$[B(\cH)^n]_r,\ [B(\cH)^n]_r^-$,
or
$$\{[X_1,\ldots, X_n]\in B(\cH)^n: \  \|[X_1,\ldots, X_n]\|=r\},
$$
where $\cH$ is an arbitrary infinite dimensional Hilbert space.
The {\it radial maximal function} $M_F:[0,1)\to [0,\infty)$ associated with a free holomorphic function $F\in Hol(B(\cX)^n_1)$ is defined by
$$
M_F(r):=\|\varphi(r)\|=\|F(rS_1,\ldots, rS_n)\|.
$$
  $M_F$ is an increasing continuous function (see the proof of Theorem \ref{f-infty}).
We define the {\it radial maximal  Hardy space} $H^p(B(\cX)^n_1)$, \ $p\geq 1$, as the set of all free holomorphic functions
$F\in Hol(B(\cX)^n_1)$  such that $M_F$ is in the  Lebesque space $ L^p[0,1]$.
Setting
$$
\|F\|_p:=\|M_F\|_p:=\left(\int_0^1\|F(rS_1,\ldots, rS_n)\|^p dr \right)^{1/p},
$$
it is easy to see that $\|\cdot \|_p$ is a norm on the linear space $H^p(B(\cX)^n_1)$.

\begin{theorem}\label{radial-Banach}
If $p\geq 1$, then the radial maximal Hardy space $H^p(B(\cX)^n_1)$ is a Banach space.
\end{theorem}

\begin{proof}
First we prove the result for $p=1$.
Let $\{F_k\}_{k=1}^\infty \subset H^1(B(\cX)^n_1)$ be a sequence such that
\begin{equation}
\label{ser-conv-1}
\sum_{k=1}^\infty \|F_k\|_1\leq M<\infty.
\end{equation}
We need to prove that $\sum_{k=1}^\infty F_k$ converges in $\|\cdot\|_1$. By \eqref{ser-conv-1}, we have
$$
\sum_{k=1}^m \int_0^1 \|F_k(rS_1,\ldots, rS_n)\| dr\leq M,\quad \text{ for any } \ m\in \NN.
$$
Using Fatou's lemma, we deduce that the function
$\psi(r):=\sum_{k=1}^\infty \|F_k(rS_1,\ldots, rS_n)\|$
is integrable on $[0,1]$.
Notice that the series $\sum_{k=1}^\infty \|F_k(rS_1,\ldots, rS_n)\|<\infty$ for any $r\in [0,1)$.
Indeed, assume that there exists $r_0\in [0,1)$ such that
$\sum_{k=1}^\infty \|F_k(r_0S_1,\ldots, r_0S_n)\|=\infty$.
Since the radial maximal function is increasing,  we have
$$
\sum_{k=1}^\infty \|F_k(rS_1,\ldots, rS_n)\|\geq \sum_{k=1}^\infty \|F_k(r_0S_1,\ldots, r_0S_n)\|=\infty
$$
for any $r\in [r_0,1)$.
Hence, we deduce that
$$
\int_0^1\sum_{k=1}^\infty  \|F_k(rS_1,\ldots, rS_n)\| dr\geq (1-r_0) \sum_{k=1}^\infty \|F_k(r_0S_1,\ldots, r_0S_n)\|=\infty,
$$
which contradicts the fact that $\psi$ is integrable on $[0,1]$. Therefore, we deduce that
$\sum\limits_{k=1}^\infty \|F_k(rS_1,\ldots, rS_n)\|$ is convergent for any $r\in [0,1)$.
Hence, the series $\sum\limits_{k=1}^\infty F_k(rS_1,\ldots, rS_n)$ is convergent in the operator norm of $B(F^2(H_n))$ for each $r\in [0,1)$.
For each $m\geq 1$, define $g_m:=\sum_{k=1}^m F_k$.
Since  $\{g_m\}_{m=1}^\infty$  is a  sequence of free holomorphic  functions  such that
$\{g_m(rS_1,\ldots, rS_n)\}_{m=1}^\infty$ is convergent in norm for each $r\in [0,1)$, we deduce that
$\{g_m\}_{m=1}^\infty$ is uniformly convergent on any closed operatorial ball $[B(\cX)^n]_r^-$, $r\in [0,1)$.
According to our noncommutative Weierstrass type result,  Theorem \ref{Weierstrass}, there is a free holomorphic function $g$ on the open  operatorial unit $n$-ball such that
$\|g_m(rS_1,\ldots, rS_n)-g(rS_1,\ldots, rS_n)\|\to 0$, as $m\to\infty,$ and
therefore
$$g(rS_1,\ldots, rS_n)=\sum_{k=1}^\infty F_k(rS_1,\ldots, rS_n)\quad \text{   for any } \ r\in [0,1).
$$
Moreover, due to  the fact that $\psi$ is integrable, we have
$$
\int_0^1 \|g(rS_1,\ldots, rS_n)\| dr \leq \int_0^1
\sum_{k=1}^\infty \|F_k(rS_1,\ldots, rS_n)\| dr <\infty,
$$
which shows that $g\in H^1(B(\cX)^n_1)$.
Now, notice that
\begin{equation*}
\begin{split}
\|g-g_m\|_1&=\int_0^1 \|g(rS_1,\ldots, rS_n)-g_m(rS_1,\ldots, rS_n)\| dr\\
&=\int_0^1\left\|\sum_{k=m+1}^\infty F_k(rS_1,\ldots, rS_n)\right\| dr\\
&\leq  \int_0^1 \sum_{k=m+1}^\infty \|F_k(rS_1,\ldots, rS_n)\| dr
\end{split}
\end{equation*}
Since $\sum_{k=1}^\infty \|F_k(rS_1,\ldots, rS_n)\|<\infty$, we have
$$
\lim_{m\to\infty} \sum_{k=m+1}^\infty \|F_k(rS_1,\ldots, rS_n)\|=0 \quad \text{ for any} \ r\in [0,1).
$$
On the other hand, $\sum_{k=m+1}^\infty \|F_k(rS_1,\ldots, rS_n)\|\leq \psi(r)$ for any $m\in \NN$.
Since $\psi$ is integrable on $[0,1]$, we can apply Lebesgue's dominated convergence theorem and deduce that
$$
\lim_{m\to\infty} \int_0^1\sum_{k=m+1}^\infty \|F_k(rS_1,\ldots, rS_n)\| dr=0.
$$
Now,  we deduce that $\|g-g_m\|_1\to 0$, as $m\to\infty$, which shows that the series $\sum_{k=1}^\infty F_k$ is convergent in $\|\cdot\|_1$. This completes the proof when $p=1$.

Assume now that $p>1$ and let $\{F_k\}_{k=1}^\infty\subset H^p((B(\cX)^n_1)$ be a sequence such that $\sum_{k=1}^\infty
\|F\|_p\leq M<\infty$.
Since $\|F_k\|_1\leq \|F_k\|_p$, we have $\sum_{k=1}^\infty
\|F\|_1\leq M$.
Applying the first part of the proof, we find $g\in H^1(B(\cX)^n_1)$ such that, for each $r\in [0,1)$,
$$
g(rS_1,\ldots, rS_n)=\sum_{k=1}^\infty F_k(rS_1,\ldots, rS_n),
$$
where the convergence is in the operator norm of $B(F^2(H_n))$.
Moreover, we have
\begin{equation*}
\begin{split}
\int_0^1\left\|\sum_{k=1}^m F_k(rS_1,\ldots, rS_n)\right\|^p dr
&\leq  \int_0^1\left( \sum_{k=1}^m \|F_k(rS_1,\ldots, rS_n)\|\right)^{p} dr\\
&\leq
\left[\sum_{k=1}^m\left(\int_0^1 \|F_k(rS_1,\ldots, rS_n)\|^p\right)^{1/p}\right]^p\\
&=
\left(\sum_{k=1}^m \|F_k\|_p\right)^p\leq M^p.
\end{split}
\end{equation*}
Using Fatou's lemma, we deduce that the function $r\mapsto \left\| \sum\limits_{k=1}^\infty F_k(rS_1,\ldots, rS_n)\right\|^p$ is integrable on $[0,1]$ and therefore
$g\in H^p((B(\cX)^n_1)$.
Notice also that
\begin{equation}
\label{norm-int}
\|g-g_m\|_p\leq
\left[\int_0^1 \left( \sum_{k=m+1}^\infty \|F_k(rS_1,\ldots, rS_n)\|\right)^p\right]^{1/p}.
\end{equation}
Since
$ \sum_{k=m+1}^\infty \|F_k(rS_1,\ldots, rS_n)\|\leq \psi$ for any $m\in \NN$, and
$$
\lim\limits_{m\to\infty} \sum_{k=m+1}^\infty \|F_k(rS_1,\ldots, rS_n)\|=0\quad \text{ for any }\ r\in [0,1),
$$
 we can apply again Lebesgue's  dominated convergence theorem and deduce that
$$
\lim_{m\to\infty}
\left[\int_0^1 \left( \sum_{k=m+1}^\infty \|F_k(rS_1,\ldots, rS_n)\|\right)^p\right]^{1/p}=0.
$$
Hence and using inequality \eqref{norm-int}, we deduce that
$\|g-g_m\|_p\to0$ as $m\to\infty$.
Consequently,
the series $\sum_{k=1}^\infty F_k$ converges in the norm $\|\cdot\|_p$.
This completes the proof.
 \end{proof}

\begin{proposition}\label{prop-Hp}
Let  $p\geq 1$.
\begin{enumerate}
\item[(i)]
If  $f\in H^\infty(B(\cX)^n_1)$, then
$\|f\|_1\leq \|f\|_p\leq \|f\|_\infty$. Moreover,
$$H^\infty(B(\cX)^n_1)\subset H^p(B(\cX)^n_1)\subset H^1(B(\cX)^n_1)\subset Hol(B(\cX)^n_1).$$
\item[(ii)] If $f\in H^\infty(B(\cX)^n_1)$, then
$$
\|f\|_\infty=\lim_{p\to\infty}\left(\int_0^1\|f(rS_1,\ldots, rS_n)\|^p dr\right)^{1/p}.
$$
\item[(iii)] If $f=\sum\limits_{k=0}^\infty \sum\limits_{|\alpha|=k} a_\alpha Z_\alpha$ is in $H^p(B(\cX)^n_1)$, then
$$
\left(\sum_{|\alpha|=k}|a_\alpha|^2\right)^{1/2}\leq (pk+1)^{1/p} \|f\|_p.
$$
\end{enumerate}
\end{proposition}

\begin{proof} Part (i) follows as in the classical  theory of $L^p$ spaces.
To prove (ii), define the function $G:[0,1]\to [0,\infty)$ by setting
$G(r):=\|f(rS_1,\ldots, rS_n)\|$ if $r\in [0,1)$ and
$G(1):=\lim\limits_{r\to 1}\|f(rS_1,\ldots, rS_n)\|$. Due to Theorem \ref{f-infty}, $G$ is an increasing  continuous  function and $G(1)=\|f\|_\infty$.
Therefore,
\begin{equation*}
\begin{split}
\lim_{p\to\infty}\left(\int_0^1\|f(rS_1,\ldots, rS_n)\|^p dr\right)^{1/p}&= \lim_{p\to\infty}
\left( \int_0^1 G(r)^p\right)^{1/p}\\
&=\max_{r\in [0,1]} G(r)=G(1)=\|f\|_\infty.
\end{split}
\end{equation*}
To prove (iii),   notice that Theorem \ref{Cauchy-est}
implies
$$
r^k \left(\sum_{|\alpha|=k}|a_\alpha|^2\right)^{1/2}\leq
\|f(rS_1,\ldots, rS_n)\|,\quad r\in [0,1).
$$
Integrating over $[0,1]$, we complete the proof of (iii).
\end{proof}

The next result extends the noncommutative von Neumann inequality from $H^\infty(B(\cX)^n_1)$ to the radial maximal  Hardy space $H^p(B(\cX)^n_1)$,\ $p\geq 1$.

\begin{theorem}\label{vN-Hp} If $T:=[T_1,\ldots, T_n]\in [B(\cH)^n]_1$ and $p\geq 1$, then the mapping $$
\Psi_T:H^p(B(\cX)^n_1)\to B(\cH)\quad \text{ defined by  } \ \Psi_T(f):=f(T_1,\ldots, T_n)
$$ is continuous, where $f(T_1,\ldots, T_n)$  is defined by the free analytic functional calculus and $B(\cH)$ is considered with the operator norm  topology. Moreover,
\begin{equation*}
 \|f(T_1,\ldots, T_n)\|\leq\frac{1}{(1-\|[T_1,\ldots, T_n]\|)^{1/p}} \|f\|_p
\end{equation*}
for any
 $f\in H^p(B(\cX)^n_1)$.
\end{theorem}

\begin{proof}
 Assume that $\|[T_1,\ldots, T_n]\|=r_0<1$ and let $f\in H^p(B(\cX)^n_1)$. Since the radial maximal function is increasing and
 and due to    Corollary \ref{max-mod2}, we have
\begin{equation*}
\begin{split}
\|f\|_p&\geq
\left(\int_{r_0}^1 \|f(rS_1,\ldots, rS_n)\|^p dr\right)^{1/p}\\
&\geq (1-r_0)^{1/p} \|f(r_0S_1,\ldots, r_0S_n)\|\\
&\geq (1-r_0)^{1/p}\|f(T_1,\ldots, T_n)\|.
\end{split}
\end{equation*}
Hence, we deduce the above von Neumann type inequality, which can be used to prove the continuity of $\Psi_T$.
\end{proof}

We remark that if $f\in H^\infty(B(\cX)^n_1)$, then  one can  recover the noncommutative von
Neumann  inequality \cite{Po-von} for strict row contractions, i.e.,
$\|f(T_1,\ldots, T_n)\|\leq \|f\|_\infty$. Indeed, take $p\to\infty$ in the above inequality   and use part (ii) of Proposition \ref{prop-Hp}.

\smallskip

In the last part of this paper,  we introduce  a  Banach space of analytic functions on the open unit ball of $\CC^n$ and obtain a von Neumann type inequality in this setting.
We use the standard multi-index notation.
Let ${\bf p}:=(p_1,\ldots, p_n)$ be a multi-index in $\ZZ_+^n$. We denote
$|{\bf p}|:=p_1+\cdots + p_n$ and ${\bf p} !:={ p}_1 !\cdots { p}_n !$.
If $\lambda:=(\lambda_1,\ldots,\lambda_n)$, then we set
$\lambda^{\bf p}:=\lambda_1^{p_1}\cdots \lambda_n^{p_n}$ and
define the symmetrized functional calculus
$$
(\lambda^{\bf p})_{\text{\rm sym}} (S_1,\ldots, S_n):=\frac {{\bf p}!} {|{\bf p}|! }\sum_{\alpha\in \Lambda_{\bf p}} S_\alpha,
$$
where
$$
\Lambda_{\bf p}:=\{\alpha\in \FF_n^+: \lambda_\alpha= \lambda^{\bf p} \text{ for any } \lambda\in \BB_n\}
$$
and $S_1,\ldots, S_n$ are the left creation operators on the Fock space $F^2(H_n)$.
Notice that card\,$\Lambda_{\bf p}=\frac {|{\bf p}|!}
{{\bf p}!}$.
  Denote by $H_{\text{\rm sym}}(\BB_n)$ the set of all
analytic functions on $\BB_n$ with scalar coefficients
 $$
f(\lambda_1,\ldots,\lambda_n):=\sum\limits_{\bf p\in \ZZ_+^n}
\lambda^{\bf p} a_{\bf p}, \quad  a_{\bf p}\in \CC,
$$
such that
\begin{equation}\label{sup-AA}
\limsup_{k\to \infty}\left(
\sum\limits_{{\bf p}\in \ZZ_+^n,|{\bf p}|=k}
 \frac
{|{\bf p}|!}{{\bf p}!}
|a_{\bf p}|^2\right)^{1/2k}\leq 1.
\end{equation}
Then
 \begin{equation*}
\begin{split}
f_{\text{\rm sym}}(rS_1,\ldots, rS_n)
&:=\sum_{k=0}^\infty
\sum\limits_{{\bf p}\in \ZZ_+^n,|{\bf p}|=k}
r^k a_{{\bf p}}[(\lambda^{\bf p})_{\text{\rm sym}} (S_1,\ldots, S_n)]\\
&=\sum_{k=0}^\infty \sum_{|\alpha|=k} r^{|\alpha|} c_{\alpha}S_\alpha,
\end{split}
\end{equation*}
where $c_{0}:=a_{0}$   and $c_{\alpha}:=
 \frac {{\bf p}!}{|{\bf p}|!}a_{{\bf p}}$ for ${\bf p}\in \ZZ_+^n$, ${\bf p}\neq (0,\ldots, 0)$, and $\alpha\in \Lambda_{\bf p}$.
It is clear that, for each $k=1,2,\ldots, $ we have
\begin{equation*}
\begin{split}
\sum_{|\alpha|=k}|c_{\alpha}|^2&=
\sum\limits_{{\bf p}\in \ZZ_+^n,|{\bf p}|=k}\left(\sum_{\alpha\in \Lambda_{\bf p}} |c_{\alpha}|^2\right)\\
&=\sum\limits_{{\bf p}\in \ZZ_+^n,|{\bf p}|=k}
 \frac
{{\bf p}!}{|{\bf p}|!}
 |a_{\alpha}|^2.
 \end{split}
\end{equation*}
Due to Theorem \ref{Abel},  condition \eqref{sup-AA} implies that $f_{\text{\rm sym}}(rS_1,\ldots, rS_n)$ is norm convergent for each $r\in [0,1)$, and  $f_{\text{\rm sym}}(Z_1,\ldots, Z_n)$ is a free holomorphic function on the open  operatorial unit $n$-ball.
We define $H_{\text{\rm sym}}^\infty(\BB_n) $ as the set of all functions $f\in H_{\text{\rm sym}}(\BB_n)$ such that
$$\|f\|_{\text{\rm sym}}:=\sup_{0\leq r<1}\left\| f_{\text{\rm sym}}(rS_1,\ldots, rS_n)\right\|<\infty.
$$
\begin{theorem}\label{sym}
$\left(H_{\text{\rm sym}}^\infty(\BB_n), \|\cdot\|_{\text{\rm sym}}\right)$ is a Banach space.
\end{theorem}

\begin{proof}
First notice that if $f\in H_{\text{\rm sym}}^\infty(\BB_n)$ then $f_{\text{\rm sym}}(rS_1,\ldots, rS_n)$ is norm convergent and
$f_{\text{\rm sym}}(Z_1,\ldots, Z_n)$ is a free holomorphic function on the open  operatorial unit $n$-ball. Using Theorem \ref{operations},
it is easy to see that
$H_{\text{\rm sym}}^\infty(\BB_n)$ is a vector space and
$\|\cdot \|_{\text{\rm sym}}$ is a norm.
Let $\{f_m\}_{m=1}^\infty$ be a Cauchy sequence of functions in
$H_{\text{\rm sym}}^\infty(\BB_n)$.
According to
Theorem \ref{f-infty}, $(f_m)_{\text{\rm sym}}\in F_n^\infty$ and
$\{(f_m)_{\text{\rm sym}}\}_{m=1}^\infty$  is a Cauchy sequence in $\|\cdot\|_\infty$, the norm of the Banach algebra $F_n^\infty$.
Therefore, there exists $g\in F_n^\infty$ such that
$\|(f_m)_{\text{\rm sym}}-L_g\|_\infty\to 0$, as $m\to\infty$.
If $
f(\lambda_1,\ldots,\lambda_n)=\sum\limits_{\bf p\in \ZZ_+^n}
a_{\bf p}^{(m)}\lambda^{\bf p}, \quad  a_{\bf p}\in \CC,
$
then $(f_m)_{\text{\rm sym}}(S_1,\ldots, S_n)=\sum_{k=0}^\infty\sum_{|\alpha|=k} c_\alpha^{(m)} S_\alpha$,
where $c_\alpha^{(m)}:=\frac
{|{\bf p}|!}{{\bf p}!} a_{\bf p}^{(m)}$
for ${\bf p}\in \ZZ_+^n$, ${\bf p}\neq (0,\ldots,0)$ and $\alpha\in \Lambda_{\bf p}$.
If $g=\sum_{\alpha\in \FF_n^+} b_\alpha e_\alpha$ is  the Fourier representation of $g$ as an element of $F^2(H_n)$, then we have
\begin{equation*}
\begin{split}
|c_\alpha^{(m)}-b_\alpha|&=
\left|\left<[(f_m)_{\text{\rm sym}}(S_1,\ldots, S_n)-L_g]1,1\right>\right|\\
&\leq \|(f_m)_{\text{\rm sym}}-L_g\|_\infty.
\end{split}
\end{equation*}
Taking $m\to \infty$,  we deduce that $c_\alpha^{(m)}\to b_\alpha$ for each $\alpha\in \FF_n^+$. Since $c_\alpha^{(m)}=c_\beta^{(m)}$ for any $\alpha,\beta \in \Lambda_{\bf p}$, we get $b_\alpha=b_\beta$.
Setting $h(\lambda_1,\ldots, \lambda_n):= \sum_{k=0}^\infty \sum_{|\alpha|=k} b_\alpha \lambda_\alpha$, one can see that $h$ is holomorphic in $\BB_n$ and $h_{\text{\rm sym}}=L_g$. Moreover, $\|h\|_{\text{\rm sym}}=\|g\|_\infty<\infty$.
This shows that
$H_{\text{\rm sym}}^\infty(\BB_n)$ is a Banach space.
\end{proof}

Now, using Theorem \ref{abel} in the scalar case, we can deduce the following.

\begin{proposition}
  If  $T:=[T_1,\ldots, T_n]\in B(\cH)^n$  is  a commuting $n$-tuple  of operators with  the joint spectral radius $r(T_1,\ldots, T_n)<1$  and
$f(\lambda_1,\ldots,\lambda_n):=\sum\limits_{\bf p\in \ZZ_+^n} a_{\bf p}
\lambda^{\bf p} $ is in
$ H_{\text{\rm sym}}(\BB_n)$, then
$$f(T_1,\ldots, T_n):=
\sum_{k=0}^\infty
\sum\limits_{{\bf p}\in \ZZ_+^n,|{\bf p}|=k}a_{\bf p} T^{\bf p}
$$
is a well-defined operator in $B(\cH)$, where the series is convergent in  the operator norm topology.
Moreover,
the map
$$
\Psi_T: H_{\text{\rm sym}}(\BB_n)\to B(\cH)\qquad \Psi_T(f)=f(T_1,\ldots, T_n)
$$
is continuous and
$$
\|f(T_1,\ldots, T_n)
\|\leq M\|f\|_{\text{\rm sym}},
$$
where $M=\sum_{k=0}^\infty \left\|\sum_{|\alpha|=k}T_\alpha T_\alpha^*\right\|^{1/2}$.
\end{proposition}

In a forthcoming paper \cite{Po-Bohr}, we obtain operator-valued Bohr type inequalities for the Banach space $H_{\text{\rm sym}}^\infty(\BB_n)$.

      \bigskip

       %

      \end{document}